\documentstyle[11pt]{article}
\makeatletter
\newcommand{\l@a}[2]{\hbox to\textwidth{#1\leaders\hbox to
0.78em{\hfil .\hfil}\hfil #2}}

\renewcommand{\l@section}{\@dottedtocline{1}{0em}{2em}}
\makeatother
\newcounter{lm}[section]
\newcounter{thm}[section]
\newcounter{prop}[section]
\newcounter{rem}[section]
\newcounter{cor}[section]
\newcounter{nexr}[section]
\newcounter{ex}[section]

\newcounter{defin}[section]

\begin{document}

\def\contentsname{\centerline {\bf Contents}}
\def\refname{\centerline {\bf References}}

\newcommand{\thm}{\refstepcounter{thm} {\bf Theorem \arabic{section}.%
\arabic{thm}.} }

\renewcommand{\thethm}{\thesection.\arabic{thm}}

\newcommand{\nthm} [1] {\refstepcounter{thm} {\bf Theorem
\arabic{section}.%
\arabic{thm}} #1}

\newcommand{\prop}{\refstepcounter{prop} {\bf Proposition
\arabic{section}.%
\arabic{prop}.} }

\newcommand{\nprop} [1] {\refstepcounter{prop} {\bf Proposition
\arabic{section}.%
\arabic{prop}} #1}

\renewcommand{\theprop}{\thesection.\arabic{prop}}

\newcommand{\cor}{\refstepcounter{cor} {\bf Corollary
\arabic{section}.%
\arabic{cor}.} }

\newcommand{\ncor} [1] {\refstepcounter{cor} {\bf Corollary
\arabic{section}.%
\arabic{cor}} #1}

\renewcommand{\thecor}{\thesection.\arabic{cor}}

\newcommand{\lemma}{\refstepcounter{lm} {\bf Lemma
\arabic{section}.%
\arabic{lm}.} }

\newcommand{\nlemma} [1] {\refstepcounter{lm} {\bf Lemma
\arabic{section}.%
\arabic{lm}} #1}

\renewcommand{\thelm}{\thesection.\arabic{lm} }

\newcommand{\ex}{\refstepcounter{ex} {\bf Example \arabic{section}.%
\arabic{ex}.} }

\newcommand{\nex} [1] {\refstepcounter{ex} {\bf Example
\arabic{section}.%
\arabic{ex}} #1}

\renewcommand{\theex}{\thesection.\arabic{ex} }

\newcommand{\defin}{\refstepcounter{defin} {\bf Definition
\arabic{section}.%
\arabic{defin}.} }

\newcommand{\ndefin} [1] {\refstepcounter{defin} {\bf Definition
\arabic{section}.%
\arabic{defin}} #1}

\renewcommand{\thedefin}{\thesection.\arabic{defin}}

\newcommand{\rem}{\refstepcounter{rem} {\bf Remark
\arabic{section}.%
\arabic{rem}.} }

\newcommand{\nrem}{\refstepcounter{rem} {\bf \arabic{section}.%
\arabic{rem}.} }

\renewcommand{\therem}{\thesection.\arabic{rem} }

\newcommand{\rems} {{\bf Remarks.}}

\newcommand{\nexr}{\refstepcounter{nexr}{\bf \arabic{section}.%
\arabic{nexr}.}}

\renewcommand{\thenexr}{\thesection.\arabic{nexr}}

\newcommand{\exr}{{\bf Exercise} }

\newcommand{\exrs} {{\bf Exercises.} }

\def\qed{\hfill $\Box$}
\newcommand{\vi}{{\varphi}}
\newcommand{\C}{{\bf C}}
\newcommand{\Pp}{{\bf P}}
\newcommand{\R}{{\bf R}}
\newcommand{\Z}{{\bf Z}}
\newcommand{\Q}{{\bf Q}}
\newcommand{\N}{{\bf N}}
\newcommand{\tf}{{\tilde {f}}}
\newcommand{\tg}{{\tilde {g}}}
\newcommand{\th}{{\tilde {h}}}
\newcommand{\tp}{{\tilde {p}}}
\newcommand{\tx}{{\tilde {x}}}
\newcommand{\tvar}{{\tilde {\varphi }}}
\newcommand{\tpsi}{{\tilde {\psi }}}
\newcommand{\tchi}{{\tilde {\chi }}}
\newcommand{\G}{{\Gamma}}
\newcommand{\bh}{{\bar {H}}}
\newcommand{\be}{{\bar {E}}}
\newcommand{\bx}{{\bar {x}}}
\newcommand{\by}{{\bar {y}}}
\newcommand{\bu}{{\bar {u}}}
\newcommand{\bA}{{\widehat {A}}}
\newcommand{\bI}{{\widehat {I}}}
\newcommand{\bG}{{\widehat {G}}}
\newcommand{\bg}{{\widehat {g}}}
\newcommand{\baf}{{\widehat {f}}}
\newcommand{\bP}{{\widehat {P}}}
\newcommand{\hvar}{{\widehat {\varphi }}}
\newcommand{\hpsi}{{\widehat {\psi }}}
\newcommand{\hx}{{\widehat {x}}}
\newcommand{\hz}{{\widehat {z}}}
\newcommand{\hy}{{\widehat {y}}}
\newcommand{\hh}{{\widehat {h}}}
\newcommand{\hf}{{\widehat {f}}}
\newcommand{\hg}{{\widehat {g}}}
\newcommand{\hr}{{\widehat {r}}}
\newcommand{\hp}{{\widehat {p}}}
\newcommand{\ha}{{\widehat {a}}}
\newcommand{\hc}{{\widehat {C}}}
\newcommand{\hE}{{\widehat {E}}}
\newcommand{\hM}{{\widehat {M}}}
\newcommand{\htt}{{\widehat {t}}}
\newcommand{\hR}{{\widehat {R}}}
\newcommand{\hX}{{\widehat {X}}}
\newcommand{\si}{{\sigma}}
\newcommand{\te}{{\tilde {E}}}
\newcommand{\gi}{{G_{\infty }}}
\newcommand{\k}{{\bar {k}}}
\newcommand{\br}{{\bar {r}}}
\newcommand{\bR}{{\bar {R}}}

\newcommand{\dX}{{\dot {X}}}
\newcommand{\dM}{{\dot {M}}}
\newcommand{\dD}{{\dot {D}}}
\newcommand{\dE}{{\dot {E}}}

\newcommand{\hA}{{\widehat {A}}}
\newcommand{\hI}{{\widehat {I}}}
\newcommand{\bF}{{\bar {F}}}
\newcommand{\dimc}{{{\dim}_{\C}\,}}
\newcommand{\De}{{Deg \,}}
\newcommand{\tA}{{\tilde {A}}}
\newcommand{\tH}{{\tilde {H}}}
\newcommand{\is}{{\simeq}}

\def\Hom{{\rm Hom\,}}
\def\Im{{\rm Im\,}}
\def\Aut{{\rm Aut\,}}
\def\det{{\rm det\,}}

\newcommand{\fr} [1] {{F^{#1}A}}
\newcommand{\po}{{\pi_0}}
\newcommand{\id}{{\bf 1}}
\newcommand{\pe}{{\pi_1}}
\newcommand{\supp}{{\rm supp ~}}
\newcommand{\const}{{\rm const}}
\newcommand{\tors}{{\rm Tors ~}}
\newcommand{\Tor}{{\rm Tor ~}}
\newcommand{\Pic}{{\rm Pic ~}}
\newcommand{\ML}{{\rm ML}}
\newcommand{\Dk}{{\rm Dk}}
\newcommand{\LND}{{\rm LND}}
\newcommand{\LNDG}{{\rm LND_{gr}}}
\newcommand{\GL}{{\rm GL}}
\newcommand{\Gr}{{\rm Gr ~}}
\newcommand{\gr}{{\rm gr ~}}
\newcommand{\Div}{{\rm Div ~}}
\newcommand{\widehatpar}{{\widehat {\partial}}}
\newcommand{\proof} {{\bf Proof.} }

\newcommand{\mod} {{\rm mod ~}}
\newcommand{\chr} {{\rm char ~}}
\newcommand{\grad} {{\rm grad ~}}
\newcommand{\spec} {{\rm spec ~}}
\newcommand{\reg} {{\rm reg ~}}
\newcommand{\dt} {{\rm det ~}}
\newcommand{\sing} {{\rm sing\,}}

\newcommand{\dg} {{\rm deg ~}}
\newcommand{\Ker} {{\rm Ker ~}}
\newcommand{\trdg} {{\rm tr.deg ~}}
\newcommand{\dgv} {{\rm deg}_{\varphi} ~}
\newcommand{\dgy} {{\rm deg}_y ~}
\newcommand{\dge} {{\rm deg}_1 ~}
\newcommand{\dgp} {{\rm deg}_\partial ~}

\newcommand{\ldm} {{\rm L-dim ~}}
\newcommand{\kdm} {{\rm K-dim ~}}
\newcommand{\dm} {{\rm dim ~}}
\newcommand{\cdm} {{\rm codim ~}}

\newcommand{\ok} {{\overline k}}
\newcommand{\ox} {{\overline x}}
\newcommand{\oy} {{\overline y}}
\newcommand{\om} {{\omega}}

\newcommand{\ff}{{\bf F}}
\newcommand{\T}{{\bf T}}
\newcommand{\sph}{{\bf S}}

\title {LECTURES on EXOTIC ALGEBRAIC STRUCTURES on AFFINE SPACES}

\vspace{.5cm}

\author{M. Zaidenberg}

\date{}
\maketitle

\begin{abstract}
\noindent These notes are based on the lecture courses given at
the Ruhr-Universit{\"a}t-Bochum (03--08.02.1997) and
at the Universit{\'e} Paul Sabatier (Toulouse, 08-12.01.1996).
\end{abstract}

\tableofcontents

\newpage

\section{Introduction}

Within the traditional
algebraic geometry of quasi-projective varieties,
the affine geometry occupies a special place, being
known as a source of difficult problems.
Let us recall the most famous ones. Hereafter we restrict the
consideration to varieties defined over $\C$.

\bigskip

\noindent {\bf 1. The Zariski Cancellation Problem}:

\smallskip

\noindent {\it Is it true that
an isomorphism $X \times \C^n \approx \C^{n+k}$
is only possible if $X \approx \C^{k}$?}

\bigskip

\noindent {\bf 2. The Structure of the Automorphism Group (Nagata)}:

\smallskip

\noindent {\it Given a polynomial
automorphism of ${\C}^n,$ can it be presented as a product
of linear and triangular ones?}

\bigskip

\noindent {\bf 3. The Linearization Problem}:

\smallskip

\noindent {\it Is any regular ${\C}^*$-action on
${\C}^n$ conjugate with a linear one in the automorphism group $\Aut \C^n?$}

\bigskip

\noindent {\bf 4. The Embedding Problem (Abhyankar, Sathaye)}:

\smallskip

\noindent {\it Is any regular embedding ${\C}^k \hookrightarrow {\C}^n$
equivalent to a linear one up to the actions of the groups
$\Aut \C^n$ on $\C^n$ resp. $\Aut \C^k$ on $\C^k?$}

\bigskip

\noindent {\bf 5. The Jacobian Problem (Keller): }

\smallskip

\noindent {\it Given a regular mapping ${\C}^n \hookrightarrow {\C}^n$
with a constant non-zero Jacobian, is it necessarily an automorphism of
${\C}^n?$}

\bigskip

\noindent To clarify the present day situation,
some comments are in order.

\medskip

\noindent 1. The affirmative answer to the question (1) for $k=2$
was the result of a series of papers by Miyanishi, Sugie and Fujita
\cite{MiySug, Fu 1} (see also \cite{Kam 1}). In higher dimensions $k \ge 3$
there is no significant progress.

An analog of the Zariski Cancellation Problem in birational setting was
answered in negative by Beauville, Colliot-Thelene, Sansuc and Swinnerton-Dyer
\cite{BCTSSD}.

There is the following more general Cancellation Problem
(see e.g. \cite{AEH, EH, Ho}):

\smallskip

\noindent {\it Given
an isomorphism of polynomial rings $A[x] \simeq B[x]$ over two rings $A$ and
$B,$ when does it follow that $A \simeq B?$}

\smallskip

\noindent
In the corresponding counterexample by Danielewski \cite{Dan}
(see also Fieseler \cite{Fi}, tom Dieck \cite{tD 3})
$A$ and $B$ are the rings of
regular functions on any two of the smooth affine surfaces
$\{x^n y + z^2 =1\} \subset \C^3, \quad n \in \N.$

\smallskip

\noindent 2. The structure of the automorphism group
Aut${\C}^n$ for $n=2$ is classically known
(Jung \cite{Ju}, van der Kulk \cite{vdK}).
Starting with $n = 3$ it is completely mysterious
(see e.g. \cite{AAS, Na 1, Na 2, Wr 1}).
For instance, it is still unknown whether or not the Nagata automorphism
$g \in \Aut \C^3,\,\,\,g\,:\,(x,\,y,\,z) \longmapsto
(x+z\Delta,\,y+2x\Delta+z\Delta^2,\,z)$ where $\Delta:=x^2-yz,$ is `tame',
i.e. decomposable as it's needed in the problem (2) \cite[(6.8)]{Na 1, Co}.

\smallskip

\noindent 3. To answer (2) it would be useful to describe the
one-parameter subgroups of the automorphism group
Aut${\C}^n,$ that is, the regular ${\C}_+$-actions and
${\C}^*$-actions on ${\C}^n$ where ${\C}_+$ resp. ${\C}^*$
denotes the additive resp. the multiplicative group of the complex number
field. It seems rather natural to expect that any ${\C}_+$-action resp. any
${\C}^*$-action on ${\C}^n$ is conjugate with a triangular one resp.
with a linear one. The former one was shown to be false starting with $n=3$
(Bass \cite{Ba}), while it is true for $n=2$ (Rentschler \cite{Re}).
It is worthwhile noticing that giving a ${\C}_+$-action on an affine
variety $X$ is the same as giving
a locally nilpotent derivation (LND for short) of the algebra
${\C}[X]$ of regular functions on $X$ \cite{Re}.

As for the latter one, i.e. for the Linearization Problem,
the positive answer is known for $n = 2$ (Gutwirth \cite{Gut}).
Below we say more about the
recent positive solution for
$n=3$ and the role of exotic $\C^3$-s
in this solution (see Koras--Russell \cite{KoRu 2, KoRu 3}, Makar-Limanov
\cite{ML 2}, Kaliman and Makar-Limanov \cite{KaML 3},  Kaliman, Koras,
Makar-Limanov, Russell \cite{KaKoMLRu}).
The positive answer for $n = 3$ implies
linearizability of any connected, reductive group action on the affine 3-space
(Kraft, Popov \cite{KrPo, Po}).
An example of a non-linearizable action of the
semi-direct product of ${\bf C}^*$ and ${\bf Z}/2{\bf Z}$ on ${\bf C}^4$ was
constructed by G. Schwarz \cite{Sch}.
It is known that the Linearization Problem restricted to semi-stable
${\bf C}^*-$actions, that is,
the actions of the multiplicative semigroup $\C,$
is equivalent to the Zariski Cancellation Problem (Kambayashi-Russell
\cite{KamRu}; see also Bass-Haboush \cite{BaHa}).

\medskip

\noindent 4. Any regular embedding
${\C} \hookrightarrow {\C}^2$ is equivalent to a linear one (Abhyankar-Moh,
Suzuki \cite{AM, Suz 1}). However, already for the
embeddings ${\C} \hookrightarrow {\C}^3$ and
${\C}^2 \hookrightarrow {\C}^3$ it is unknown whether the analogous fact
is true.
But the embeddings ${\C}^k \hookrightarrow
{\C}^n$ are linearizable as soon as $n \ge 2(k+1)$
(Jelonek \cite{Je}, Kaliman \cite{Ka 4, Ka 5}, Nori, Srinivas \cite{Sr}).

\smallskip

\noindent 5. There is a number of equivalent formulations of
the Jacobian Problem and partial results; see e.g.
\cite{An, AAS, BCW, Dr, Kam 2, Or 1, Wa}; see also \cite{Pi}
for the Jacobian Problem for $n=2$ over $\R.$
False proofs which appear regularly
definitely certify the difficulty of the problem.

\bigskip

Not so far ago, new unusual objects appeared in the affine geometry;
they were called {\it exotic ${\C}^n.$} These are
smooth affine varieties diffeomorphic, but non-isomorphic to the affine spaces.
Actually, the first example for $n=3$ was constructed in
a deep paper of Ramanujam \cite{Ram} where also the
non--existence of exotic ${\C}^2$ was proven.
Later on, many examples of smooth contractible affine $n-$folds
for all $n \ge 3$ were found (Choudary--Dimca
\cite{ChoDi}, tom Dieck \cite{tD 1, tD 2},
tom Dieck--Petrie \cite{tDP 2}, Dimca \cite{Di 1},
Kaliman \cite{Ka 1, Ka 2}, Kaliman-Zaidenberg
\cite{KaZa}, Koras--Russell \cite{KoRu 2}, Russell
\cite{Ru 1}, the author \cite{Za 2, Za 3, Za 4}).
Some of them occurred to be exotic $\C^n-$s, for
the other ones this is still unknown.
The main difficulty here is the absence of
effective tools for recognition of exotics, or, in other words,
for recognition of the affine spaces.

In the work of Koras and Russell on the Linearization Problem for $n = 3$
(see \cite{KoRu 1, KoRu 2})
it was reduced to classification problems
for a series of smooth contractible threefolds $X \subset {\C}^4$
(the {\it Koras--Russell threefolds}), on one hand, and for a series of affine
singular quotient surfaces, on the other one.
As for the latter one, it was recently settled completely \cite{KoRu 3}.

As for the former one, it consists of clarification
whether or not all the Koras-Russel threefolds $X \subset {\C}^4$ are
exotic ${\C}^3$-s.
The first partial results were obtained  by Kaliman and Russell
\cite{Ka 1, Ru 1}. They succeeded to show that
the logarithmic Kodaira dimension is non-negative for at least some of these
threefolds.

A methods suggested by Kaliman and Makar-Limanov \cite{KaML 1} allowed them to
enlarge this class.
Namely, it was shown that under certain restrictions on $X$ there is no
dominant regular mapping ${\C}^3 \to X.$

But all the above
methods failed to distinguish from ${\C}^3$ a certain subseries of the
Koras-Russel threefolds.
The Russell cubic threefold $X \subset {\C}^4$
given by the equation $x + x^2y + z^2 +t^3 =0$ is one of them.
It looks especially simple, but in fact, this one is
the most difficult to analyze. Its geometric structure
can be described as follows. It contains `the book-surface'
$B := \{x = 0\}\subset X$ which is isomorphic to the product
${\C} \times \Gamma_{2,\,3}$ where
$\Gamma_{2,\,3} \subset {\C}^2$ is the affine cuspidal cubic $z^2 +t^3 =0.$
The complement $X \setminus B$ is isomorphic to ${\C}^* \times {\C}^2.$
Thus, $X$ is obtained from ${\C}^3$ after replacing ${\C}^2
\subset {\C}^3$ by the book-surface $B.$ Notice also that there exists a
dominant morphism $\C^3 \to X.$ Using the fact that $B$ is contractible,
one can show that $X$ is contractible, too.
It follows from the Smale h-Cobordism
Theorem that, actually, $X$ is diffeomorphic to ${\bf R}^6.$

\medskip

Finally, Makar-Limanov \cite{ML 2} succeeded to prove that the Russell cubic
is an exotic ${\C}^3.$
Soon after, Kaliman and Makar-Limanov \cite{KaML 3}, along the same approach,
showed that all the  Koras-Russel threefolds
are exotic ${\C}^3$-s. Thus, the Linearization Problem for ${\C}^3$
was answered in positive \cite{KoRu 2, KoRu 3}.

The proof of Makar-Limanov \cite{ML 2} is based on the use of
locally nilpotent derivations (LND, for short).
The principal new ingredients suggested in
\cite{ML 2} consist in

\smallskip

\noindent $*$ using Jacobian derivations; in particular,

\smallskip

\noindent $*$ reducing the study of general LND-s to study of Jacobian
LND-s;

\smallskip

\noindent $*$ introducing and systematically using generalized
degree functions, and then

\smallskip

\noindent $*$ reducing the study of the LND-s of a filtered ring to
those of the associated graded ring.

\bigskip

In section 7 below we present a simplified proof of
the Makar-Limanov Theorem due to  Derksen \cite{De}. In section 2 we
deal with contractible and more general acyclic surfaces. They serve
as a base for constructing exotic $\C^n$--s, but certainly merit being studied
on their own right. Sections 3--6 are devoted to constructions
of exotic $\C^n$--s. Besides, in section 6 examples of computations of the
logarithmic Kodaira dimension are given.

To simplify presentation we often restrict it to
particularly interesting examples. We do not address at all,
or say very little on closely related subjects such as analytically exotic
structures (see e.g.
\cite{Ka 2, Za 3, Za 5}), deformations of exotic structures (see \cite{FlZa 1,
Za 5}), $\bf Q$-acyclic surfaces (see e.g. \cite{FlZa 1, Fu 2, Miy 2, Or 2}),
the positive characteristic case, etc.
The interested reader can find additional information and open
problems in \cite{OPOV, Za 5}.

It is my pleasure to thank Profs. H. Flenner and G. Schumacher,
who suggested to give a lecture course on exotic structures at the
Graduiertenkolleg of the Ruhr-Universit{\"a}t-Bochum, 03--07.02.1997,
as well as the organizers of the school
`Structures exotiques de ${\C}^n$' at the Universit{\'e} Paul Sabatier,
Toulouse, 08-12.01.1996, and especially, Mme Laurence Fourrier,
for an analogous suggestion.
The author is grateful to Shulim Kaliman and Yuli Rudyak, who looked through
the text and made many useful comments; to Konstantin Sonin for his help in
editing the LaTeX-version of these notes.

\section{Acyclic surfaces}
\subsection{The first acquaintance}

By the Hironaka Resolution of Singularities Theorem, any
smooth quasi-projective variety $X$ admits
a smooth projective completion $V$ by a divisor $D$ with simple
normal crossings\footnote{i.e. all the irreducible
components of $D$ are smooth hypersurfaces in $V,$
and for any point $p \in D,$ the divisor $D$ can be given in appropriate local
coordinates $(z_1,\dots,z_n)$ on $V$ with center at $p$ by an equation of the
form $z_1\cdot\dots\cdot z_k=0$ where $k \le n.$}$X=V\setminus D.$
We call $(V,D)$ an {\it SNC-completion}
of $X$ or an {\it SNC-pair}. A variety $X$ is {\it acyclic} if
$H_{*}(X,\Z)\simeq \Z.$

\bigskip

\lemma (Fujita \cite{Fu 2}) {\it Let $X$ be a
smooth quasi-projective surface. If $X$
is acyclic, then it is affine.}

\label{lmf1}
\medskip

\proof
Assume that the surface $X$ is acyclic. Let $V$ be a smooth completion of
$X$ by a reduced divisor $D$ (not necessarily SNC).
Let $D=\sum_{i=1}^kD_i$ where each $D_i$ is an irreducible component
of $X.$ We will show that there exists an effective ample divisor
$A=\sum_{i=1}^ka_iD_i$ supported by $D,$ i.e. such that $a_i>0\,\forall i =
1,\dots,k.$
Thus, $mA$ for $m$ large enough is a hyperplane section (for the
embedding $\Phi_{\;\vert\;mA\;\vert\;}:V\hookrightarrow \Pp^N$). Hence,
$$X=V\setminus D=V\setminus \supp(mA)\hookrightarrow \Pp^N\setminus H\simeq
\C^N$$
is affine. By the Nakai-Moi\v sezon criterion, it suffices to choose
any $A$ as above such that $A^2>0$ and $AC>0$ for any irreducible curve $C$ in
$V.$

In view of acyclicity of $X$, from the standard topological dualities
(see the proof of Proposition \ref{pro1} below) it follows that
the natural homomorphism $H_2(D)\longrightarrow H_2(V)$ is surjective, and
$D$ is connected. (In fact, to prove that $X$ is affine we use only
these two conditions. It is well known that the boundary of an irreducible
affine variety is connected, so, the second one is
necessary.) Set $\sum=\{A=\sum_{i \in I}a_iD_i\;\vert\;a_i>0 \quad
\forall i\in I, AD_i>0 \quad
\forall i\in I\},\,\,\,I \subset \{1,\dots,k\}$. First, we show that $\sum$ is
non-empty.
Indeed, let $H\in \Div V$ be any ample divisor.
The classes of $D_i,$ $i=1,\ldots,k,$ (which we denote by the same
letters) generate the group $H_2(V,\Z),$ and so $H=\sum_{i=1}^k h_iD_i=
\sum_{i\in I} a_iD_i-\sum_{j\in J} a_jD_j = A_0 - B_0$,
where $I,J\subseteq \{1,\ldots,k\},\,\,I \neq \emptyset,\,\,I\bigcap
J=\emptyset$, and $a_i > 0 \quad \forall i \in I \cup J$.
For any irreducible curve $C$ in $V,$ we have $A_0C-B_0C=HC>0,$
whence $A_0C>B_0C.$ Given $C=D_i,$ $i\in I,$ this implies
$A_0C>B_0C\ge 0.$ Therefore, $A_0\in \sum.$

Suppose that $A\in \sum,$ $\supp A\ne D,$
$D_j$ is not contained in $\supp A,$ and
$D_jA>0.$ Then $m_jA+D_j\in \sum$ for some $m_j>0.$
Indeed, $(m_jA+D_j)D_i>0$ for all $D_i\subset \supp A,$
and $(m_jA+D_j)D_j>0$ when $m_j>-D^2_j/D_jA.$

Recall that the divisor $D$ is connected. Therefore,
starting with $A_0$ and applying the procedure as above, in
a finite number of steps one can find a divisor $A\in \sum$ with $\supp A=D$.
Clearly, $A^2>0,$ $AD_i>0$ for all $i=1,\ldots, k,$
and $AC\ge DC$ for any irreducible curve $C$ such that $C$ is not
contained in $D.$ Since $A_0C\ge HC>0,$ we have $DC>0,$ whence also $AC>0.$
Thus,
$A$ is an ample divisor, and $\supp A=D.$
$\Box$

\bigskip

\rem In higher dimensions the analogous statement is not true, in
general, as an example of Winkelmann \cite{Win} shows. In this example
$X = Q \setminus E$ is a contractible non-affine (and even non-Stein)
quasi-projective variety where $Q$ is a smooth projective quadric of dimension
$4,$ and $E \subset Q$ is a codimension $2$ smooth subvariety.

\bigskip

\prop {\it Let $X=V\setminus D$ where $V$ is a smooth projective
surface, and $D$ is a curve in $V.$ Then $X$ is acyclic if and only if
the following conditions hold:

$(i)$ $\po (D)=\pe (V)= \pe (D)=\id.$

$(ii)$ $i_{*}: H_2(D,\Z)\longrightarrow H_2(V,\Z)$ is an isomorphism.}

\label{pro1}
\medskip

\proof
By the Lefschetz duality \cite{Do} we have
$$H^i(V,D) \simeq  H_{n-i}(X)\,, \quad H_i(V,D) \simeq H^{n-i}(X)\,,
\,\,\,i=0,\dots,4\,.$$
Assume that $X$ is acyclic; then the above groups are zero for $i=0,\dots,3$.
{}From the standard exact sequences of a pair (all the
homology groups have coefficients in $\Z$):

\begin{picture}(400,90)
\put(140,26){\vector(-1,1){30}}
\put(130,41){$i_{*}$}
\put(5,10){$H_{*}(V,D)$}
\put(76,56){\vector(-1,-1){30}}
\put(40,41){$r_{*}$}
\put(74,60){$H_{*}(V)$}
\put(65,12){\vector(1,0){60}}
\put(88,18){{\small $\partial$}}
\put(130,10){$H_{*}(D)\,\,,$}

\put(314,53){\vector(1,-1){30}}
\put(333,41){$i^{*}$}
\put(205,10){$H^{*}(V,D)$}
\put(249,23){\vector(1,1){30}}
\put(246,41){{\small $r^{*}$}}
\put(278,59){$H^{*}(V)$}
\put(325,12){\vector(-1,0){58}}
\put(288,18){{\small $\partial^{*}$}}
\put(331,10){$H^{*}(D)\,\,,$}
\end{picture}

\noindent where $\dg i_{*}=\dg i^{*}=\dg r_{*}=\dg r^{*}=0,$ $\dg \partial=-1,$
and $\dg \partial^{*}=1,$ it follows that
$H_i(D)\simeq H_i(V),$ $H^i(D)\simeq H^i(V)$ for all
$i,$ $0\le i\le 3.$
In particular, $H_2(D)\simeq H_2(V),$ which proves $(ii)$. Also,
$H^3(V)\simeq H^3(D) = 0.$ By the Poincar{\'e} duality, $H_1(V)\simeq
H_1(D)=0,$
whence $\pe(D)=\id.$ Since $H_0(D)\simeq H_0(V)\simeq \Z,$ $\po(D)=\id.$

The proof of Lemma \ref{lmf1} shows that $D$ is a
hyperplane section. By the Lefschetz Theorem on hyperplane sections
\cite{Lef, AF}, \cite[Thm.  7.4]{Mil 2}, the homomorphism
$i_{*}:\pe(D)\longrightarrow \pe(V)$
is a surjection. This implies that $\pe(V)=\id,$ as claimed.

Conversely, assume that the conditions $(i)$ and $(ii)$ are satisfied. Then
$H_1(V)=H_1(D)=0,$ whence $H^3(V)=0$, by the Poincar{\'e} duality. Furthermore,
the group $H_2(V)\simeq H^2(V)\simeq H^2(D)\simeq \Z^d$ is free (here $d$
stands for the number of irreducible components of $D$).
Since $\tors H^j=\tors H_{j-1},$ $H^1(V)$ is also a free group. Hence
$H^1(V)\simeq H_1(V)=0,$ and so $H_3(V)=0.$ Then $H_i(D)\simeq H_i(V),$
$i=0,\ldots,3,$ whence $H_i(V,D)=0,$ $i=0,\ldots,3.$ Also, we have the
same equalities for cohomologies.

By the Lefschetz duality, $0=H^i(V,D)=H_{n-i}(X),$ $i=0,\ldots,3.$
Therefore, the surface $X$ is acyclic.
$\Box$

\medskip

\cor {\it All the irreducible components of the boundary divisor
$D$ are rational curves without
self-intersections, and are arranged as a tree.}

\bigskip

\defin
Let $D=\sum_{i=1}^dD_i$ be an SNC-curve on
a projective surface $V.$ {\it The dual graph} $\G_D$ of $D$
is the graph which possesses the irreducible components $\{D_i\}$ of $D$
as vertices, and $[D_i,D_j]$ ($i\ne j$) is an edge of $\G_D$ iff $D_iD_j>0.$
Each vertex $D_i$ of $\G_D$ is weighted by $D^2_i.$

If $X=V\setminus D$ is acyclic, then, by Proposition \ref{pro1}, $\G_D$ is a
tree.

\bigskip

\nthm{({\bf Gurjar-Shastri} \cite{GuSha}, {\bf Gurjar, Pradeep, Shastri}
\cite{GuPraSha}).}
{\it Every smooth acyclic surface is rational. Moreover, the same is true for
the
$\Q-$acyclic surfaces $X$ (i.e. such that $H_*(X;\,\Q)=0$)
with at most quotient singularities. }

\medskip

\rem The following general {\bf Van de Ven Problem} \cite{VdV}:

\smallskip

\noindent {\it Is it true that every smooth
contractible affine (or even quasi-projective) variety is rational?}

\smallskip

\noindent is still open (for this and the related Hirzebruch problem on
the description of compactifications
of $\C^n$ see e.g. \cite{MS, Fur, Pro}).

\bigskip

\cor {\it An SNC-pair $(V,D)$ is a completion of an acyclic surface
$X=V\setminus D$
iff $D$ is a rational tree on a smooth rational surface $V$ such that the
Picard group $\Pic V$ is freely generated over $\Z$
by the irreducible components of $D,$ i.e. $\Pic V\simeq G(D)\simeq \Z^d.$}

\medskip

\proof
Indeed, in the case of a rational surface $V$ we have $\Pic V\simeq H_2(V).$
\qed

\bigskip

\defin
An SNC-pair $(V,D)$ is called {\it minimal} if no contraction
of a component of $D$ leads to a new SNC-pair.
(Equivalently, $\G_D$ has neither linear nor end vertices weighted by
$-1;$ recall the Castelnuovo criterion.)

\bigskip

Let $S$ be a smooth compact real manifold with the connected boundary
$\partial S$ and with the interior $X = S \setminus \partial S.$
Then the fundamental group at infinity $\pi_1^{\infty}(X)$ of $X$
can be identified with the group $\pi_1(\partial S).$
This group is known to be a topological invariant of the open manifold $X.$
Notice that smooth affine surface $X$ can be presented as the interior of a
compact real manifold with a connected boundary.

\bigskip

\nthm{({\bf Ramanujam} \cite{Ram}).} {\it $(a)$ Assume that $(V,D)$ is a
minimal
SNC-completion of a smooth acyclic surface $X=V\setminus D.$
Then $X\simeq \C^2$ iff the dual graph $\G_D$ is linear.

\smallskip

\noindent $(b)$ Furthermore, a smooth contractible surface $X$ is isomorphic to
$\C^2$ iff it is simply connected at infinity, i.e.
$\pi_1^{\infty}(S)={\bf 1}.$}

\label{ram}
\bigskip

\ex {\bf The Hirzebruch surface} $\sum_n$ where $n > 0$ is a $\Pp^1-$bundle
over
$\Pp^1:$ $\sum_n\stackrel{p}\longrightarrow \Pp^1$ such that
there exists a unique section $E_n\subset \sum_n$ with $E^2_n=-n;$
for $n= 0$ one takes $\sum_0=\Pp^1 \times \Pp^1$ with a constant section as
$E_0.$
If $F_\infty$ is a fiber over the point $\infty=(1:0)\in \Pp^1,$
then $\sum_n\setminus (E_n\bigcup F_\infty)\simeq \C^2,$ and
the dual graph of this completion of $X = \C^2$ looks like

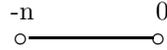
\begin{picture}(300,40)
\put(159,10){$\circ$}
\put(161,20){{\small -n}}
\put(168,13){\line(1,0){47}}
\put(211,10){$\circ$}
\put(216,20){{\small 0}}
\end{picture}

\noindent Note that the standard completion $(\Pp^2,\Pp^1)$ of $\C^2$
has the dual graph

\begin{picture}(300,40)
\put(180,10){$\circ$}
\put(183,20){{\small 1}}
\end{picture}

\ex {\bf The Ramanujam surface} \cite{Ram}. There exists an arrangement
of a smooth conic $C_2$ and a cuspidal cubic $C_1,$  $\{zx^2-y^3=0\}$ say,
in $\Pp^2$ such that $(C_1C_2)_A=1,$ and $(C_1C_2)_B=5$ where
$A,B$ are the smooth intersection points.

Let $\sigma_a:V\longrightarrow \Pp^2$ be the blow-up of $\Pp^2$
at $A$ with the exceptional (-1)-curve $E\subset V,$ and let
$C_1^\prime,C_2^\prime \subset V$ be the proper transforms of $C_1,C_2.$ Set
$D=C_1\bigcup C_2.$
We have $H_2(V)\simeq \Pic V\simeq \Z H'+\Z E$ where $H^\prime$ is the proper
transform of a generic line $H$ in $\Pp^2.$ Since $C_1\sim 3H,$
and $C_2\sim 2H$ in $\Pic \Pp^2\simeq \Z,$ we get $(C_1^\prime,C_2^\prime)=
T(H^\prime,E)$ where $T$ is the unimodular matrix

$$\left(\begin{array}{cc}
3 & 2  \\
-1 & -1
\end{array}
\right)
$$

\noindent Thus, $H_2(V)\simeq H_2(D),$ and so it follows
from Proposition \ref{pro1} that the surface $X=V\setminus D$ is acyclic.
By Fujita's Lemma \ref{lmf1}, $X$ is affine. The resolution
graph of $D\subset V$ looks as follows:

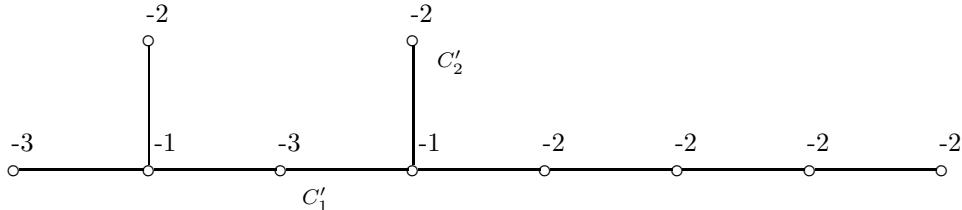
\begin{picture}(450,80)
\put(7,20){{\small -3}}

\put(2,10){$\circ$}
\put(10,13){\line(1,0){47}}

\put(53,10){$\circ$}
\put(61,20){{\small -1}}
\put(59,15){\line(0,1){45}}
\put(53,59){$\circ$}
\put(58,69){{\small -2}}

\put(61,13){\line(1,0){46}}
\put(103,10){$\circ$}
\put(108,20){{\small -3}}
\put(114,0){{\scriptsize $C_1^\prime$}}

\put(111,13){\line(1,0){46}}
\put(153,10){$\circ$}
\put(161,20){{\small -1}}
\put(159,15){\line(0,1){45}}
\put(153,59){$\circ$}
\put(158,69){{\small -2}}
\put(165,52){{\scriptsize $C_2^\prime$}}

\put(161,13){\line(1,0){46}}
\put(203,10){$\circ$}
\put(208,20){{\small -2}}

\put(211,13){\line(1,0){46}}
\put(253,10){$\circ$}
\put(258,20){{\small -2}}

\put(261,13){\line(1,0){46}}
\put(303,10){$\circ$}
\put(308,20){{\small -2}}

\put(311,13){\line(1,0){46}}
\put(353,10){$\circ$}
\put(358,20){{\small -2}}
\end{picture}

\noindent This graph is minimal and non-linear, so the Ramanujam
Theorem \ref{ram}$(a)$ yields that $X$ is not isomorphic to $\C^2.$

\label{exR}
\bigskip

\exrs (\nexr) Show that $\pe(X)=\id,$ and so $X$ is contractible.

(\nexr) The boundary $S$ of a ('tubular') neighborhood
of $D$ in $V$ (= an attached boundary of $X$) is not simply connected;
what is $\pe(S)?$ (see e.g. \cite{Mu, Hir} for an algorithm of computing
$\pe(S)$). Hence, the surface $X$ is not homeomorphic to $\R^4.$

\bigskip

Next we give some more examples of contractible surfaces,
following  \cite{Za 1, Za 4}.

\medskip

\ex \label{ex3} Let $T$ be a matrix
of non-negative integers of the form

$$T = \left(
\begin{array}{cccc}
m_{00} & 0 & n_{00} & 0 \\
m_{10} & 0 & 0 & n_{10} \\
0 & m_{01} & n_{01} & 0 \\
0 & m_{11} & 0 & n_{11}
\end{array}
\right)$$

\noindent Consider the lines $l_{i,j}\simeq \Pp^1$ in the quadric $Q :=
\Pp^1\times \Pp^1$
where $l_{i,j}=\{ix+(1-i)y=j\}_{i,j=0,1}$ (i.e.
$x=0,\,\,x=1,\,\,\,y=0,\,\,y=1$).
Blow up over the points $z_{i,j}=(i,j),$ $i,j=0,1$, until
the  four rational functions
$$\frac{(x-i)^{m_{ij}}}{(y-j)^{n_{ij}}}\,,\,\,\,i,j=0,1\,,$$
become regular; denote the resulting surface by
$V_T\stackrel{\pi}\longrightarrow Q.$ Set $D_0=e_0\bigcup e_1 \bigcup
\{l_{ij}\}\subset Q$ where
$e_0 = \{y = \infty \}, \,\,e_1 = \{x = \infty \}\subset Q$ (see Figure 1).

\begin{picture}(500,120)

\put(170,81){\line(1,0){79}}
\put(170,52){\line(1,0){79}}
\put(170,21){\line(1,0){79}}

\put(185,13){\line(0,1){78}}
\put(215,13){\line(0,1){78}}
\put(245,13){\line(0,1){78}}

\put(252,18){{\small $l_{00}$}}
\put(178,0){{\small $l_{10}$}}
\put(252,48){{\small $l_{01}$}}
\put(208,0){{\small $l_{11}$}}
\put(252,78){{\small $e_0$}}
\put(238,0){{\small $e_1$}}

\put(163,25){{\small $z_{00}$}}
\put(217,25){{\small $z_{10}$}}
\put(163,55){{\small $z_{01}$}}
\put(217,55){{\small $z_{11}$}}

\put(180,19){{\small $\bullet$}}
\put(210,19){{\small $\bullet$}}
\put(180,50){{\small $\bullet$}}
\put(210,50){{\small $\bullet$}}

\end{picture}

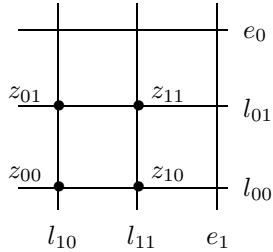
\centerline{\bf Figure 1}

\bigskip

\noindent Set $\pi^{-1}(D_0)=D_T\bigcup \{v_{ij}\}$
where $v_{ij}$ are the only $(-1)$ curves in the
exceptional divisor of $\pi.$ Then $D_T$ is a rational tree
(deleting $v_{ij}$ is called `cutting cycles'
by tom Dieck and Petrie \cite{tDP 1, tDP 3}).
The proper transforms $e_0',\,e_1'$ of
$e_0,e_1$ and the components of the exceptional divisor
of the blowing up $\pi$ form a natural basis in $\Pic V_T=H_2(V_T,\,\Z).$
A new one is given by the components of $D_T$ provided that
the decomposition matrix
$$\left(\begin{array}{cc}
T & 0  \\
B & I
\end{array}
\right)
$$
is unimodular, i.e. that $T$ is unimodular: $\dt T=\pm 1.$
Now Proposition \ref{pro1} asserts that under this condition the surface
$X_T:=V_T\setminus D_T$ is acyclic.

\bigskip

\rem Tom Dieck and Petrie \cite{tDP 1, tDP 3}
found all the basic line
arrangements in $\Pp^2$ that lead to acyclic surfaces in
a blowing up and cutting cycles process as above;
there are seven of them.
The first one (see Figure 2) depends on discrete and
continuous parameters
(cf.  also the Classification Theorem 3.3(d) below);
the other six are projectively rigid.

\begin{picture}(500,120)

\put(30,15){\line(1,0){300}}
\put(348,13){{\small $l_{s+1}$}}
\put(100,0){\line(1,1){100}}
\put(135,55){{\small $p_0$}}
\put(147,50){{\small $\bullet$}}
\put(88,-5){{\small $l_0$}}
\put(130,0){\line(2,3){68}}
\put(142,20){{\small $p_1$}}
\put(134,13){{\small $\bullet$}}
\put(132,-5){{\small $l_1$}}
\put(160,0){\line(1,3){34}}
\put(166,20){{\small $p_2$}}
\put(159,13){{\small $\bullet$}}
\put(162,-5){{\small $l_2$}}
\put(200,20){{$\cdots$}}
\put(280,0){\line(-1,1){100}}
\put(266,20){{\small $p_s$}}
\put(260,13){{\small $\bullet$}}
\put(283,-5){{\small $l_s$}}

\end{picture}

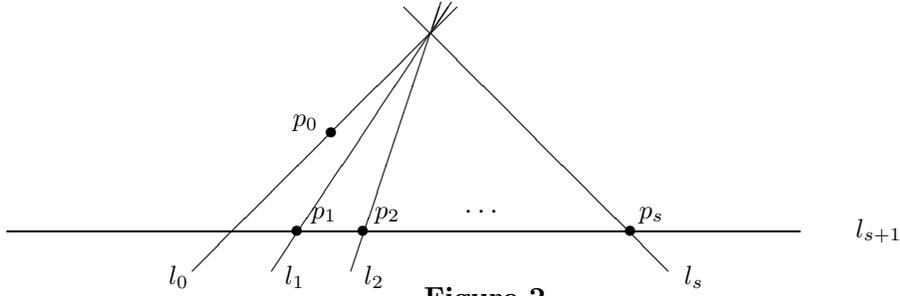
\centerline{\bf Figure 2}

\smallskip

\noindent However, not every smooth acyclic surface arises from
a line arrangement in $\Pp^2;$ see \cite{tD 4}.

\bigskip

Recall the following notion.

\medskip

\defin Let $D$ be a closed analytic hypersurface in a
connected complex manifold $M.$ By a {\it vanishing loop} of $D$ at
a smooth point $e \in D$
we mean any loop $\delta_D = \delta_{D,\,e_0} (= \alpha^{-1}\beta\alpha)$ in
$M \setminus D$ consisting of a path
$\alpha$ which
joins a base point $e_0 \in M \setminus D$ with a point
$e' \in \omega \setminus D$ of a small complex disc $\omega \subset M$
transversal to $D$ at $e$ and a simple loop $\beta$ in positive direction in
$\omega \setminus D$
with the base point $e'.$ We also denote by $\delta_D$
the corresponding element of the fundamental group $\pi_1(M \setminus
D,\,e_0).$
It is uniquely defined up to conjugation by the corresponding irreducible
component of
the hypersurface $D.$

\medskip

The next lemma completes Example \ref{ex3}.

\medskip

\lemma {\it If T is unimodular then $X_T$ is a contractible surface.}

\medskip

\proof The unimodularity condition $\dt T=\pm 1$ implies that $X$ is an acyclic
surface.
Thus, in virtue of the Hurewicz and Whitehead Theorems it is enough to
show that it also implies simply connectedness of $X.$
Denote
$$X_0:= X_T \setminus \{v_{i,j}\}\simeq Q\setminus D_0$$
where, as before, $D_0$ is the union of six fixed
generators of the quadric $Q.$
The vanishing loops $a_i,\,b_j \in \pi_1(X_0),\,\,i,j=0,1,$ of the lines
$l_{i,j} \subset \hookrightarrow Q$ provide generators for the group
$\pe(X_0),$
that is, $$\pe(X_0)=\ff_2\times \ff_2=<a_0,a_1,b_0,b_1\;\vert\;[a_i,b_j]=1, \;
i,j=0,1>.$$
We will see in Lemma \ref{fuji} below that $\pe(X)=\pe(X_0)/N$ where
$$N=<<a_i^{m_{ij}}b_j^{n_{ij}}\;\vert\; i,j=0,1>>$$ is the
minimal normal subgroup generated by these four products.
Thus, we have the following relations in $\pe(X):$
$a_0^{m_{00}}b_0^{n_{00}}=1,\quad a_1^{m_{10}}b_0^{n_{10}}=1,$
and so $a_0^{m_{00}n_{10}}=a_1^{m_{10}n_{00}}.$
Also, $\quad a_0^{m_{01}}b_1^{n_{01}}=1, \quad a_1^{m_{11}}b_1^{n_{11}}=1,$
and hence $a_0^{m_{01}n_{11}}=a_1^{m_{11}n_{01}}.$
It follows that
$$a_0^{m_{00}n_{10}m_{11}n_{01}}=a_0^{m_{01}n_{11}m_{10}n_{00}},$$
whence $a_0^{\dt T}=1,$ that is $a_0=1.$ In the same way we obtain
$a_1=b_0=b_1=1,$ and therefore $\pe(X)=\id.$
$\Box$

\bigskip

\lemma {\it $(a)$ Let $D$ be a closed hypersurface in a
complex manifold $M,$ $\dim M\ge 2.$ Then the group Ker$\,\{i_{*}\,:\,
\pe(M\setminus D) \to \pe(M)\}$ is generated by the vanishing loops
of $D$.
In particular, if $D$ is irreducible, this kernel is generated, as a
normal subgroup, by any of these loops.}

$(b)$ ({\bf Fujita} \cite[(7.18)]{Fu 2}) {\it Let $M$ be a surface,
$D_1,D_2$ be two curves in $M,$ and
$p$ be an intersection
point which is an ordinary double point of $D_1\bigcup D_2$.
Let $\sigma_p:M^\prime\longrightarrow M$ be the blow-up at $p.$
Then (the class of) a vanishing loop $\alpha_E$ of the exceptional
$(-1)-$curve $E\subset M^\prime$ of $\sigma_p$ in the group
$\pe(M\setminus (D_1 \cup D_2)) = \pe(M'\setminus (E \cup D_1' \cup D_2'))$
can be represented
as $\alpha_E=\alpha_{D_1^\prime}\alpha_{D_2^\prime}$,
where $\alpha_{D_1^\prime},\,\alpha_{D_2^\prime}$ are vanishing loops
of the proper transforms $D_i^\prime$ of $D_i,$ $i=1,2.$ Moreover, the classes
$\alpha_{D_1^\prime},\,\,\alpha_{D_2^\prime}$ commute. }

$(c)$ ({\bf Fujita} \cite[(7.18)]{Fu 2}) {\it
In the notation as in $(b)$, let $\sigma:\widehat M\longrightarrow M$
be a sequence of blow-ups over $p$ such that $D_1^{*}=mE+\ldots,$
$D_2^{*}=nE+\ldots$ for $E$ being the exceptional $(-1)-$curve of the
last blow-up. Then we have
$\alpha_E=\alpha_{D_1^\prime}^m\alpha_{D_2^\prime}^n.$}

\label{fuji}
\medskip

\proof
$(a)$ Denote $\Delta$ the unit disc in $\C$. Let $\gamma:\partial
\Delta=\sph^1\longrightarrow M\setminus D.$
Observe that $\gamma_{*}\in {\rm Ker}\, i_{*}$ iff there exists $ \tilde
\gamma: \Delta \longrightarrow M$ such that $\tilde
\gamma \;\vert\;_{\partial \Delta}=\gamma.$ After applying a small smooth
deformation
we may assume that $\tilde \gamma(\Delta)$ meets $D$ transversally at smooth
points
$p_1,\ldots,p_k \in D;$ let $q_1,\ldots, q_k \in \Delta$ be the corresponding
disc points.

Choose disjoint vanishing loops of $q_1,\ldots, q_k$ in $\Delta$, and contract
the circle
$\sph^1 = \partial \Delta$ onto their union. Being composed with
$\tilde \gamma$
this yields a desired homotopy of $\gamma$ to a product of vanishing loops of
$D.$

$(b)$ Let $D_i = \{z_i = 0\},\,i=1,2,$ in a local chart $(z_1,\,z_2)$ in $M$
centered at $p$. Representing the product $\alpha_{D_1} \alpha_{D_2}$ on the
torus
$\;\vert\;z_1\;\vert\;=\varepsilon,$
$\;\vert\;z_2\;\vert\;=\varepsilon,\,\,\varepsilon>0,$ as its diagonal section,
after blowing up at the origin
this loop becomes a vanishing one $\alpha_E.$

$(c)$ Apply induction on the number of blow-ups.
$\Box$

\bigskip

\exrs (\nexr) (after Fujita \cite{Fu 2})
 Let $X$ be a smooth acyclic (resp. contractible) surface, and let $C\subset X$
be an irreducible simply connected curve. Consider the blow-up
$\sigma_p:\widehat X\longrightarrow X$ at a smooth point $p\in C,$
and set $X^\prime=\widehat X\setminus C^\prime$ where
$C^\prime \subset \widehat X$ is the
proper transform of $C.$ Show that the surface $X^\prime$ is also acyclic
(resp. contractible).

\smallskip

(\nexr) Draw the dual graph $\Gamma_{D_T}$ where $D_T$ is as in
Example \ref{ex3} above. Deduce that in many cases the surface $X_T$
is not isomorphic to $\C^2.$

\subsection{Elements of classification: the logarithmic Kodaira dimension}

\defin
Let $L\longrightarrow V$ be a line bundle
(i.e. a one-dimensional algebraic vector bundle)
over a smooth projective variety
$V,$ and let $H^0(V,L)$ be the space of its regular sections. By the
Cartan-Serre Theorem, $h^0(V,L):=\dm H^0(V,L)<\infty.$ Suppose that
$h^0(V,L)>0$, and fix a basis $s_0,\ldots,s_n$ of the vector space
$H^0(V,L)$ where
$n=h^0(V,L)-1.$ Then $Z:=\{z\in V\;\vert\; s_0(z)=\ldots=s_n(z)=0\}$ is
a proper subvariety of $V.$ For $z\in V\setminus Z$ fix a vector space
structure in the fiber $L_z\cong \C;$ then the point
$\Phi_L(z):=(s_0(z):\ldots:s_n(z))\in \Pp^n$ is well-defined, and the
rational map $\Phi_L:V\longrightarrow \Pp^n$ is regular in $V\setminus Z.$
The line bundle $L$ is called {\it very ample} if $\Phi_L$ is an embedding
(assuming $Z = \emptyset$); {\it ample}
if $mL$ is very ample for some $m>0;$ {\it big}
if
$\ldm (V):=\overline {\lim}_{m\rightarrow \infty} \dm \Phi_{mL} =\dm_{\C} V.$
(In section 6.2 below we also use the notation
$k(V,\,L) := \ldm (V).$)
Put $\ldm (V)=-\infty$ if $h^0(mL)=0$ $\forall m.$

\bigskip

\nthm{({\bf  Serre-Siegel-Kodaira}; see e.g. \cite[Thm.  10.2]{Ii 3}).}
{\it For some
$m_0>0$ we have
$h^0(mm_0L)\sim m^l$ where $l=\ldm (V).$}

\label{sk}
\bigskip

\defin
If $L=K_V$ is the canonical line bundle (i.e.
$K_V=\Lambda^n T^{*}V$ where $n=\dm_{\C} V$), then
$k(V):=\kdm V$ is called the {\it Kodaira dimension} of $V$;
$k(V) \in \{-\infty,\,0,\,1,\dots,\dm V\}$.
If $k(V)=\dm V,$ then $V$ is said to be of {\it general type}.

\medskip

\noindent Thus, the projective variety $V$ is of general type
iff the canonical line bundle
$K$ of $V$ is big, i.e.
for some $m>0,$ $\Phi_{mK}:V\hookrightarrow \Pp^n$ is a birational
embedding.

\bigskip

\exr (\nexr) A smooth irreducible projective curve $V$ is
of general type, i.e. $k(V)=1,$ iff $g(V)\ge 2;$ $k(V)=0$ iff
$g(V)=1$, i.e. if $V = \T_{\Lambda} := \C/{\Lambda}$ is an elliptic curve,
where $\Lambda = \Z + \tau \Z \subset \C,\,\,\tau \in \C,\,\,{\rm Im}\,\tau>0,$
is a plane lattice;
$k(V)=-\infty$ iff $g(V)=0$, i.e. if $V\cong \Pp^1$ is a rational curve.

\bigskip

\defin
Let $(V,D)$ be an SNC-completion of a smooth quasi-projective variety
$X=V\setminus D.$ The {\it log-Kodaira
dimension}, or the {\it Iitaka--Kodaira dimension} of $X$ is
$\overline k(X):=\ldm (V)$ where $L=K+D$ (Iitaka \cite{Ii 1}).
$X$ is said to be {\it of log-general type} if $\overline k(X)=\dm X.$

Sometimes
$K+D$ is called {\it the log-canonical divisor}; the holomorphic sections of
$O(K + D)$ correspond to
the meromorphic forms regular in $X$ which can be written as
$$a\frac{dz_{i_1}}{z_{i_1}}\wedge \ldots \wedge \frac{dz_{i_k}}{z_{i_k}}
\wedge dz_{i_{k+1}}\wedge \ldots \wedge z_{i_n}$$
in local coordinates in $V$ where $D=\{z_{i_1}=
\ldots =z_{i_k}=0\}.$

\bigskip

\nthm{({\bf Iitaka} \cite{Ii 1}, \cite[Ch. 11]{Ii 3}).}
{\it $\overline k(X)$ is an
invariant of
$X$ which does not depend on the choice of an SNC-completion $(V,D)$
of $X.$}

\bigskip

\exr (\nexr)
Let $X$ be a smooth irreducible quasi-projective curve. Show that
$\overline k(X)=-\infty$ iff $X=\C$ or $X=\Pp^1;$
$\overline k(X)=0$ iff $X=\T_\Lambda$ or $X=\C^{*} := \C \setminus \{0\};$
$\overline k(X)=1$ otherwise.

\bigskip

\nthm{({\bf The main properties of the log-Kodaira
dimension}).}

\noindent $(a)$ (\cite{Ii 1}, \cite[Thm.  11.3]{Ii 3})
{\it $\overline k(X\times Y) =
\overline k(X)+\overline k(Y).$}

\smallskip

\noindent $(b)$ (\cite{Ii 1}, \cite[Prop. 11.5]{Ii 3}) {\it If $Y$ is a
Zariski open subset of $X,$ then
$\overline k(Y)\ge \overline k(X),$ and $\overline k(Y)= \overline k(X)$
if ${\cdm}_X(X \setminus Y) \ge 2.$}

\smallskip

\noindent $(c)$ ({\bf The Iitaka Easy Addition Theorem} \cite[Thm.  4]{Ii 1},
\cite[Thm.  11.9]{Ii 3})
{\it If $\pi :Y\longrightarrow X$ is a surjective morphism of
smooth quasi-projective varieties with a connected generic
fiber $F,$ then $\overline k(Y)\le \overline k(F)+\dm X$.}

\smallskip

\noindent $(d)$ ({\bf The Kawamata-Viehweg Addition Theorem} \cite{Kaw 1, Vie})
{\it If, in addition, $\dm F=1$ then $\overline k(Y)\ge
\overline k(F)+\overline k(X).$}

\smallskip

\noindent $(e)$ ({\bf The Logarithmic Ramification Formula} \cite{Ii 1},
\cite[Thm.
11.3]{Ii 3})
{\it Let $\dm X=\dm Y,$ and let
$f :Y\longrightarrow X$
be a dominant morphism\footnote{i.e. $f(Y)$ contains a Zariski open subset of
$X$.}. By the Hironaka
Resolution of Singularities Theorem, $f$ can be extended to a morphism
$\overline f:V_Y\longrightarrow V_X$ where
$(V_X,D_X)$ (resp. $V_Y, D_Y)$) is an appropriate SNC-completion
of $X$ (resp. of $Y$). Then there
exists an effective divisor $R_{\overline f}\subset V_Y$ (which is
called the {\rm logarithmic ramification divisor}) such that
$$K_{V_Y}+D_{V_Y}={\overline f}^{*}(K_{V_X}+D_{V_X})+R_{\overline f}.
\eqno{(R)}$$
In particular,
$$H^0(V_X,m(K_{V_X}+D_{V_X})) \hookrightarrow H^0(V_Y,
\overline f^{*}m(K_{V_X}+D_{V_X}))\subset
$$
$$H^0(V_Y,
m\overline f^{*}(K_{V_X}+D_{V_X})+mR_{\overline f})=H^0(V_Y,
m(K_{V_Y}+D_{V_Y})).$$
Therefore, $\overline k(X)\le \overline k(Y).$}

$(f)$ (\cite[Prop. 1, Thm.  3]{Ii 1}, \cite[Thms. 10.5, 11.10]{Ii 3}) {\it If,
in addition, $f$ is either a proper birational morphism, or an {\'e}tale
covering,
then we may assume
$R_{\overline f}$ being an $f-$exceptional divisor, i.e. $\cdm {\overline
f}(R_{\overline f})\ge 2,$
and we have $\overline k(Y)=\overline k(X)$.}

\label{prr}
\bigskip

\nthm{({\bf Classification of acyclic surfaces}).} {\it Let $X$ be
an acyclic surface. Then the following assertions hold.}

\smallskip

\noindent $(a)$ ({\bf Miyanishi-Sugie-Fujita} \cite{MiySug, Fu 1}) {\it
$\overline
k(X)=-\infty$ iff $X\simeq \C^2.$}

\smallskip

\noindent $(b)$ ({\bf Fujita} \cite{Fu 2}) {\it If $X$ is non-isomorphic to
$\C^2,$
then $\overline k(X)\ge 1.$}

\smallskip

\noindent $(c)$ ({\bf Iitaka-Kawamata} \cite[Thm.  5]{Ii 1}, \cite{Kaw 2})
{\it If $\overline k(X)=1,$ then there exists
a morphism $X\longrightarrow \Gamma$ onto a smooth curve
$\Gamma$ with generic
fibers isomorphic to $\C^{*}=\C\setminus \{0\}$ (called a}
$\C^{*}-$fibration)\footnote{Actually, any (not necessarily acyclic) affine
surface $X$ with $\overline k(X)=1$ possesses a $\C^{*}-$fibration.}.

\smallskip

\noindent $(d)$ ({\bf Gurjar-Miyanishi} \cite{GuMiy 1};
cf. also \cite{tDP 1, FlZa 1})
{\it There exists a complete list of acyclic surfaces with
$\overline k(X)=1.$ Any such surface can be obtained from a
tom Dieck--Petrie line
configuration in $\Pp^2$ of the first kind (see Figure 2 above)
by a composition $\sigma\,:\,V \to \Pp^2$ of blowing ups over the points
$p_0,\ldots,p_s$
to get, as $\sigma^{-1}(p_i),\,i=1,\dots,s,$ a linear chain
of rational curves with only one $(-1)-$curve.
All these curves except for the
last $(-1)-$curve of each chain are components
of the boundary divisor $D\subset V,$ as well
as the proper preimage of the
original line arrangement.
In addition, all the blow ups over $p_0$ are {\it outer}, i.e.
each of them is done at a smooth point of the total transform of $l_0,$
whereas under the {\it cutting cycle} procedure over $p_i,\,i=1,\dots,s,$
the blow ups are {\it inner}, i.e. they are done only at double points
of the preimage of the original line arrangement.

For each $i=1,\ldots,s,$ we fix a rational number $\frac{m_i}{n_i}.$
The blow-up process over $p_i$ is done
according to the data
$(m_i,n_i);$ this means that locally at $p_i$ it resolves the
point of indeterminacy of the rational function
$x^{m_i}/y^{n_i},\,\,i=1,\dots,s.$
The numbers $(m_i,\,n_i)_{i=1,\dots,s}$ must satisfy
a Diophantine equation of
unimodularity which guarantees the acyclicity
of the resulting open surface
$X = V \setminus D.$

In particular, all the
contractible surfaces with
$\overline k(X)=1$ are obtained in this way for $s=2$
and $m_1n_2+m_2n_1-m_1m_2=\pm 1,\,\,m_i > n_i, \,i=1,2.$
Their minimal dual graphs look as follows:

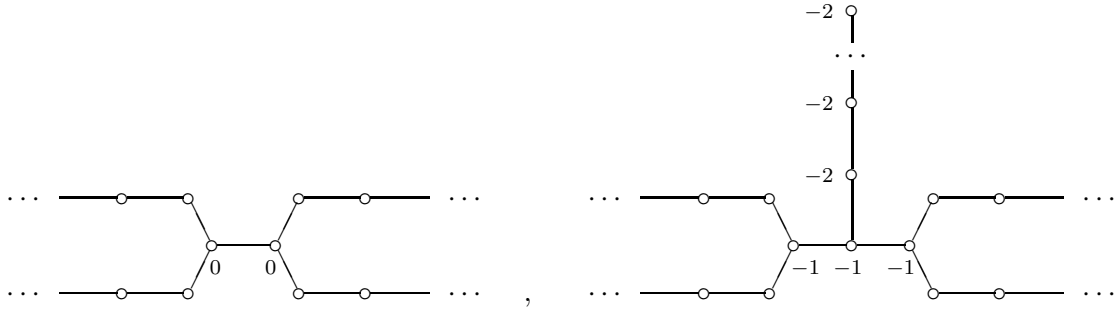
\begin{picture}(450,140)

\put(2,10){$\cdots$}

\put(25,13){\line(1,0){22}}
\put(43,10){$\circ$}

\put(51,13){\line(1,0){21}}
\put(68,10){$\circ$}

\put(75,15){\line(1,2){7}}

\put(2,46){$\cdots$}

\put(25,49){\line(1,0){22}}
\put(43,46){$\circ$}

\put(51,49){\line(1,0){21}}
\put(68,46){$\circ$}

\put(75,47){\line(1,-2){7}}

\put(77,28){$\circ$}
\put(85,31){\line(1,0){15}}

\put(79,20){{\scriptsize $0$}}

\put(169,10){$\cdots$}

\put(165,13){\line(-1,0){22}}
\put(135,10){$\circ$}

\put(139,13){\line(-1,0){21}}
\put(110,10){$\circ$}

\put(115,15){\line(-1,2){7}}

\put(169,46){$\cdots$}

\put(165,49){\line(-1,0){22}}
\put(135,46){$\circ$}

\put(139,49){\line(-1,0){21}}
\put(110,46){$\circ$}

\put(115,47){\line(-1,-2){7}}

\put(101,28){$\circ$}
\put(105,31){\line(-1,0){15}}
\put(100,20){{\scriptsize $0$}}
\put(196,10){$\,,\,$}


\put(222,10){$\cdots$}

\put(245,13){\line(1,0){22}}
\put(263,10){$\circ$}

\put(271,13){\line(1,0){21}}
\put(288,10){$\circ$}

\put(295,15){\line(1,2){7}}

\put(222,46){$\cdots$}

\put(245,49){\line(1,0){22}}
\put(263,46){$\circ$}

\put(271,49){\line(1,0){21}}
\put(288,46){$\circ$}

\put(295,47){\line(1,-2){7}}

\put(297,28){$\circ$}
\put(305,31){\line(1,0){18}}

\put(299,20){{\scriptsize $-1$}}

\put(409,10){$\cdots$}

\put(405,13){\line(-1,0){22}}
\put(375,10){$\circ$}

\put(379,13){\line(-1,0){21}}
\put(350,10){$\circ$}

\put(355,15){\line(-1,2){7}}

\put(409,46){$\cdots$}

\put(405,49){\line(-1,0){22}}
\put(375,46){$\circ$}

\put(379,49){\line(-1,0){21}}
\put(350,46){$\circ$}

\put(355,47){\line(-1,-2){7}}

\put(341,28){$\circ$}
\put(345,31){\line(-1,0){18}}
\put(335,20){{\scriptsize $-1$}}


\put(319,28){$\circ$}
\put(315,20){{\scriptsize $-1$}}

\put(325,33){\line(0,1){23}}
\put(319,55){$\circ$}
\put(304,55){{\scriptsize $-2$}}

\put(325,60){\line(0,1){23}}
\put(319,82){$\circ$}
\put(304,82){{\scriptsize $-2$}}

\put(325,87){\line(0,1){10}}
\put(315,100){$\cdots$}
\put(325,108){\line(0,1){10}}
\put(319,117){$\circ$}
\put(304,117){{\scriptsize $-2$}}
\end{picture}

\noindent The choice of centers of the
blow-ups over $p_0$ (besides the first one), and the positions of the points
$p_3,\dots,p_s$ on $l_{s+1}$ (once the first three intersection points on
$l_{s+1}$ have been fixed) give rise to the parameters of a versal deformation
of the surface $X$ or, more accurately, of its minmal SNC-completion
$(V_{\rm min},\,D_{\rm min}),$ which is defined in a unique way} \cite{FlZa 1}.

\label{clthm}
\bigskip

\nthm{({\bf Simply connected curves on acyclic surfaces}
\cite{AM, LiZa, GuMiy 2, Suz 1, Za 1}).}
{\it Let $X$ be a smooth acyclic surface, and let $\Gamma$
be an irreducible simply connected curve in  $X$. Then either

\smallskip

\noindent $*$ $(X,\,\Gamma) \simeq (\C^2,\,\Gamma_{k,\,l})$,
where $\Gamma_{k,\,l}:=\{x^k - y^l = 0\} \subset \C^2,\,\,\,k\ge l \ge
1,\,\,(k,\,l)=1$, or

\smallskip

\noindent $*$ $\ok(X) = 1$ and $\Gamma = E \setminus D \simeq \C$ where
$E\subset V$ is the last $(-1)-$curve over the point $p_0$ in the
reconstruction process as in Theorem 3.3(d) above, and $D$ is the boundary
divisor of the corresponding SNC-completion of $X$.}

\label{cusuthm}
\bigskip

In particular, this theorem shows that, up to automorphisms of the affine
plane, there is only a sequence of
irreducible simply connected curves in $\C^2$ (namely, $\{\Gamma_{k,\,l}\}$).
Each smooth acyclic surface
of logarithmic Kodaira dimension $1$ contains exactly one such curve, and
this curve is smooth. At last, there is no simply connected curves at all
on acyclic surfaces of log-general type.
See also \cite{GuMiy 2, GuPa} for some further information.

\bigskip

\ex {\bf Tom Dieck-Petrie surfaces} \cite{tDP 2}. \label{ex3.1} The surface
$X_{k,\,l}\subset \C^3$ given
by the equation $$\frac{(xz+1)^k-(yz+1)^l}{z}=1$$
where $k>l\ge 2,\,\,(k,\,l)=1$ is a smooth contractible one with
$\overline k(X_{k,\,l})=1$ (see Examples \ref{exDP} and \ref{ex6.1} below).
The only
simply connected curve in $X_{k,\,l}$ is given by the equation $z = 0.$

In a similar way,
any smooth contractible surface with $\overline k=1$
can be properly embedded into
$\C^3$ (Kaliman, Makar-Limanov \cite{KaML 2}).

\bigskip

\ex
It can be shown (see \cite{Za 3, Za 4}) that $\overline k(X_T)=1$ for
a surface $X_T$ as in Example \ref{ex3} iff
$m_{ij}=n_{ij}=1$  for a pair of diagonal points from the square vertices
$(z_{ij} = (i,\,j))_{i,\,j = 0,\,1}$ (see Figure 1 above). If so,
then the only simply connected curve in $X_T$ is the proper transform
of the
corresponding diagonal line. Otherwise, $\overline k(X_T)=2$, i.e.
$X_T$ is of log-general type.

\bigskip

\rem There is a number of examples of acyclic or even contractible
surfaces of log-general type (see e.g. \cite{tD 2, FlZa 1, GuMiy 1, Sug}), but
no classification is known. While acyclic surfaces of log-Kodaira dimension $1$
admit deformations (see the Classification Theorem 3.3(d)),
those of log-general type are rigid in all known examples \cite{FlZa 1, FlZa
2}. So, the problem arises \cite{OPOV, FlZa 1}:

\smallskip

\noindent {\it Is it true that any smooth ($\Q-$)acyclic
surface of log-general type is rigid?}

\section{Exotic product structures}

We begin this section by recalling

\medskip

{\bf The Zariski Cancellation Problem.}
{\it Given an isomorphism $X \times \C^k \stackrel{\Phi}
\simeq \C^{n + k},$ does it follow
that $X\simeq \C^n$?}

\medskip

\noindent Take $\C^n$ generic in $\C^{n+k}$, and combine $\Phi$ with the
first projection.
This yields a surjective morphism $\C^n\longrightarrow X.$
Thus, by Theorem \ref{prr}.(e), $\ok(X)=-\infty.$ Clearly, $X$ is
homotopically trivial; in particular, for $n = 2$ $X$ is an acyclic
surface $\ok(X)=-\infty.$ By the Miyanishi-Sugie-Fujita Theorem 2.3.$(a)$,
$X\simeq \C^2$.
This provides the positive answer to the Zariski Cancellation Problem
for $n=1,2.$ For $n\ge 3$ the problem is open.
In this respect the following fact could be useful.

\bigskip

\nthm{({\bf The Iitaka-Fujita Strong Cancellation Theorem} \cite{IiFu}).}
{\it Let $X,Y$ be smooth quasi-projective varieties of the
same dimension, and let $\Phi:Y\times \C^n\longrightarrow X\times \C^n$
be an isomorphism. Assume that $\ok(X)\ge 0.$ Then
there is a commutative diagram

\begin{picture}(200,95)
\unitlength0.2em
\put(62,25){$Y \times \C^k$}
\put(105,25){$X \times \C^k$}
\put(84,27){$\vector(1,0){15}$}
\put(66,5){$Y$}
\put(84,6){$\vector(1,0){15}$}
\put(109,5){$X\,\,$}
\put(69,22){$\vector(0,-1){11}$}
\put(72,16){{\rm pr}}
\put(111,22){$\vector(0,-1){11}$}
\put(114,16){{\rm pr}}
\put(91,10){$\varphi$}
\put(91,30){$\Phi$}
\end{picture}

\noindent where $\varphi$ is an isomorphism.}

\label{thif}
\bigskip

We use below the following well known corollary of the Smale
h-cobordism Theorem.

\bigskip

\nprop{(see \cite[\S 9]{Mil 1}).} {\it Let $D^n$ be a smooth
simply connected manifold of (real) dimension $n\ge 5$
with a simply connected boundary. Then the following conditions
are equivalent:

$1)$ $D^n$ is diffeomorphic to the closed unit $n-$ball $\overline B^n.$

$2)$ $D^n$ is homeomorphic to $\overline B^n.$

$3)$ $D^n$ is contractible.

$4)$ $D^n$ is acyclic.}

\label{prsm}
\bigskip

\nthm{({\bf Dimca-Ramanujam} \cite{Di 1, Ram}).} {\it Let $X$ be a contractible
smooth affine algebraic variety. If dim$_{\C}\,X =n \ge 3$ then $X$ is
diffeomorphic to $\R^{2n}.$}

\label{thdr}
\medskip

\proof\footnote{Cf. the proof of the Lefschetz Hyperplane Section Theorem in
\cite{Mil 2}.} Fix a closed embedding $X \hookrightarrow \C^N$ such that the
smooth
function $\varphi:= ||z||^2\,\vert\,X$ on $X$ is a {\it Morse function}, i.e.
it has only non-degenerate critical points (see \cite[Thm. 6.6]{Mil 2}).
Since the smooth mapping $\varphi\,:\,X \to \R$
is proper and has only finite number of critical
values, for $R > 0$ large enough the domain
$X_R := \{\varphi < R\}$ in $X$ is diffeomorphic to the whole manifold
$X.$ Denote $S_R = \partial \overline{X_R};$ that is,
$\overline{X_R}$ is a smooth manifold with the boundary $S_R.$ By the Morse
Theory applied to the Morse function $\psi:= R - \varphi$ on $X,$ the manifold
$\overline{X_R}$ can be obtained, starting
with the boundary $S_R$, by  successively gluing handles
of indices equal to those of the critical points of $\psi$ on $X_R$.

If $p \in X_R$ is a critical point of $\psi$, then ind$_p \psi = 2n -
$ind$_p\varphi$. But ind$_p\varphi \le n$ \cite[the proof of Thm.  7.2]{Mil 2}.
Hence, ind$_p \psi \ge n \ge 3$. Therefore, $\overline{X_R}$ is obtained
from $S_R = \partial \overline{X_R}$ by attaching handles of indices at least
$3$. Consequently, $\overline{X_R}$ is homotopically equivalent to a cell
complex obtained from $S_R$ by successively attaching cells of dimension at
least $3$. It follows that the first two relative homotopy groups
$\pi_i(\overline{X_R}, \,S_R),\,\,i=1,\,2,$ are trivial. Since
$\overline{X_R}$
is contractible, applying the exact homotopy sequence of a pair
$${\bf 1} = \pi_2(\overline{X_R}, \,S_R)
\stackrel{{\partial_{*}}}\longrightarrow \pe(S_R)
\stackrel{{i_{*}}}\longrightarrow \pe(\overline{X_R})={\bf 1}$$
we conclude that $\pe(S_R)={\bf 1}$. Now the theorem follows from Proposition
\ref{prsm}.
$\Box$

\bigskip

\rem In section 2 above some examples have been given
(see e.g. Examples \ref{exR}, \ref{ex3})
of  smooth contractible affine surfaces $S$
with non-simply connected attached boundaries $\partial S$ (in other words,
$S$ is not simply connected at infinity: $\pi_1^{\infty} (S) \neq {\bf 1}$).
Therefore, these contractible surfaces are not homeomorphic to $\R^4$.
This shows that the restriction $n \ge 3$ in the above theorem is crucial.

\bigskip

\cor {\it Let $X$ be a smooth contractible surface. Then $X\times \C$
is diffeomorphic to $\C^3\simeq \R^6,$ and so $X\times \C^k$
is diffeomorphic to $\R^{2k+4},$ $k\ge 1.$}

\label{lmo}
\medskip

\proof We indicate, following Ramanujam \cite{Ram}, an alternative direct
proof of this corollary.
According to Proposition \ref{prsm},
it suffices to show that $X\times \C$ is diffeomorphic to
the interior of a smooth compact manifold $D$ with a simply connected
boundary $\partial D.$ There are two natural ways to
compactify $X\times \C.$ First,
consider any smooth affine variety $Z\hookrightarrow \C^N.$ Then
the restriction $\varphi$ of the real polynomial
$\vert \vert z\vert \vert^2$ to $Z$ has only a finite number
of critical values, and hence, $\varphi^{-1}[R,\infty[$ for $R$ large enough
is diffeomorphic to $[R,\infty[\times T$ where $T:= \varphi^{-1}(R).$
Thus, $Z$ is diffeomorphic to $Z_0:=\varphi^{-1}[0,R[,$
the interior of the manifold with boundary
$\overline Z_0:=\varphi^{-1}[0,R],$ $\partial \overline Z_0=T.$
Represent in this way $X\hookrightarrow \C^n$ attaching the boundary
$\partial X,$
and $Y:=X\times \C\hookrightarrow \C^{n+1}$ attaching the boundary
$\partial_1 Y = \varphi^{-1}(R_1).$
Since $Y$ is diffeomorphic to $X\times \Delta$ where
$\Delta=\{\vert z\vert <1\},$ $Y$ can also be compactified
by attaching the non-smooth boundary
$\partial_2 Y:=(\partial X\times {\overline \Delta}) \bigcup ({\overline X}
\times \sph^1).$  In fact, $\partial_2 Y= \psi^{-1}(R_2)$ where $\psi({\bar
x},\,z) :=
$max$\,\{\vert \vert {\bar x}\vert \vert^2,\,\vert  z \vert^2\},$ and $R_2>0$
is large enough. By the Van Kampen Theorem, $\partial_2 Y$ is simply connected.

We may assume that sufficiently large $R'_1,\, R''_1,\,R'_2,\,R''_2$ are chosen
in such a way that
$\partial_1 Y = \varphi^{-1}(R_1)\subset \psi^{-1}([R'_2,\,R''_2])\subset
\varphi^{-1}([R'_1,\,R''_1]),$ and that $\varphi^{-1}([R'_1,\,R''_1])\approx
              \partial_1 Y \times [R'_1,\,R''_1],\,\,\psi^{-1}([R'_2,\,R''_2])
\approx
\partial_2 Y \times [R'_2,\,R''_2]$. Thus, the composition of embeddings
$\partial_1 Y \hookrightarrow \partial_2 Y \times [R'_2,\,R''_2]
\hookrightarrow \partial_1 Y \times [R'_1,\,R''_1]$ provides a homotopical
equivalence. Respectively, the induced isomorphism $\pi_1(\partial_1 Y)
\stackrel{\simeq}\longrightarrow \pi_1(\partial_1 Y \times [R'_1,\,R''_1])$
factors through the trivial one
$\pi_1(\partial_1 Y) \to \pi_1(\partial_2 Y \times [R'_2,\,R''_2]) \simeq
\pi_1(\partial_2 Y) = {\bf 1}.$
This proves simply connectedness of the boundary
$\partial_1 Y,$ and the assertion follows.
\qed

\bigskip

\defin By an {\it exotic $\C^n$} we mean a smooth
affine variety diffeomorphic to $\R^{2n}$ but non-isomorphic to $\C^n.$

\medskip

\nthm{({\bf Exotic product structures}).}

\noindent {\it $(a)$ Let $S$ be a smooth contractible surface non-isomorphic to
$\C^2.$ Then the product $S\times \C^{n-2}\,\,\,(n > 2)$ is an exotic $\C^n.$

\smallskip

\noindent $(b)$ Furthermore, if two smooth contractible surfaces
$S_1,\,\,S_2$ are not isomorphic, then $S_1\times \C^{n-2},\,\,S_2\times
\C^{n-2} \,\,\,(n > 2)$ are two non-isomorphic exotic $\C^n-$s.}

\label{contrsur}
\medskip

\proof $(a)$ By Lemma \ref{lmo}, $S\times \C^{n-2}$ is diffeomorphic to
$\R^{2n}$
for $n\ge 3.$ By the
Miyanishi-Sugie-Fujita Theorem 2.3.$(a)$, $\ok(S)\ne -\infty$
(otherwise $S \simeq \C^2$), whence $\ok(S)\ge 0.$
But if $S\times \C^{n-2}$ were isomorphic to $\C^n,$ then by Theorem
\ref{prr}$(e)$ we would have
$\ok(S)= -\infty,$ a contradiction.

\smallskip

\noindent $(b)$ We have to show that
the classification of
exotic product structures on $\C^n$ of the type $S\times \C^{n-2}$ where
$S$ is a surface as above is reduced to the classification
of surfaces $S$ themselves. Indeed,
$S_1\times \C^{n-2}\simeq S_2\times \C^{n-2}$ and $\ok(S_1)\ge 0$
would imply that $S_1\simeq S_2.$ Since $S_1\not\simeq S_2,$ and both surfaces
are acyclic, by the Miyanishi-Sugie-Fujita Theorem 3.3$(a)$, $\ok(S_i) \ge 0$
for
at least one value of $i$, say, for $i=1$, and so, the assertion
follows from the Iitaka-Fujita Strong Cancellation Theorem \ref{thif}.
$\Box$

\bigskip

\rems \nrem For instance, a sequence of pairwise non-isomorphic
surfaces $X_T$ of log-general type (see Example \ref{ex3} above)
yield a sequence of pairwise non-isomorphic exotic $\C^n$-s \cite{Za 3}.
Since contractible surfaces $S$ with
$\ok(S)=1$ admit deformations (see the Classification Theorem 3.3(d)), the
corresponding
exotic $\C^{n}$-s of product type $S\times \C^{n-2}$ admit deformations, too
\cite{FlZa 1}.

\bigskip

\nrem Let $X=\prod_{i=1}^nS_i$ be a product of
$n\ge 2$ smooth contractible surfaces. Then $X$ is diffeomorphic to the
interior of a compact contractible variety with boundary. By the Van Kampen
Theorem, the boundary $\partial X$
is simply connected. Therefore, by the h-cobordism Theorem, $X$ is
diffeomorphic to
$\C^{2n}.$ Also, $\ok(X)=\sum_{i=1}^n\ok(S_i).$ Hence,  $X$ is of log-general
type iff $S_i$ are so for all $i=1,\ldots,n;$ $\ok(X)=-\infty$
if $\ok(S_i)=-\infty$ for at least one value of $i.$ This shows, in particular,
that for any $n \ge 2$ there exist exotic $\C^{2n}$-s of log-general
type.

\bigskip

\nrem If $\ok(S)=2,$ then the exotic $\C^3$ $\,\,X=S\times \C$ contains no
copy of $\C^2,$ i.e. there is no embedding $\C^2\hookrightarrow S\times \C$
\cite{Za 3}.
(This is based on the fact that the surface $S$ contains no simply connected
curve; see
\cite{Za 1} and Theorem \ref{cusuthm} above.)
In the next section we present examples of exotic $\C^3$ with many copies of
$\C^2$ (see Example \ref{exkl}).

\bigskip

\nrem Due to the Ramanujam Theorem \ref{ram}$(b)$, there is no exotic
$\C^2.$

\bigskip

\nrem Actually, the Zariski Cancellation Problem
can be reformulated as follows:

\smallskip

\noindent {\it Given an exotic $\C^n,$ denote it $X,$ should also the product
$X \times \C^k$ be an exotic $\C^m\,\,(m=n+k)$?}

\bigskip

\exr (\nexr) Verify that a smooth irreducible quadric hypersurface
$X$ in $\C^{n+1}$ is contractible if and only if it is isomorphic to $\C^n,$
and if so, then the embedding $X \hookrightarrow \C^n$ is rectifiable.

\section{Contractible affine modifications}
\subsection{The Kaliman modification}

\defin Consider a triple $(M,D,C)$,
where  $M\supset D\supset C,$ $M$ and $C$ are smooth
quasi-projective (or, more generally, complex)
varieties, $D$ is an irreducible hypersurface in $M,$ and $C$ is
a proper subvariety
in $D$ contained in the smooth part $\reg D := D \setminus \sing D$ of $D,$
and so $\cdm_MD=1$ and $\cdm_MC\ge 2.$ Let
$\sigma_C: \hM\longrightarrow M$
be the blow up of $M$ with center $C$ and with the
exceptional divisor $E=\sigma^{-1}(C).$ Then
$\sigma_C \vert \,E:E\longrightarrow C$ is a fiber bundle
with the fiber $\Pp^k,$ $k=\dm E-\dm C;$ $E$ and the
proper transform $D^\prime$ of $D$ meet transversally, and
$\sigma_C: E\bigcap D^\prime\longrightarrow C$ is a fiber
bundle with the fiber $\Pp^{k-1}.$ The variety
$M^\prime:=\hM\setminus D^\prime$ is called the {\it Kaliman modification}
of the triple  $(M,D,C)$
along $D^{\prime}$ with center $C$ and with the exceptional divisor
$E^\prime=E\setminus D^\prime$ \cite{Ka 2}. Clearly, the restriction
$\sigma_C \vert \,{E^\prime}: E^\prime\longrightarrow C$  is a fiber
bundle with the fiber $\C^k.$

\medskip

One can show that the affine modification $X'$ of an affine variety $X$ is
again an affine variety \cite[Lemma 3.3]{Ka 2}.

\bigskip

In the proof of the next lemma we use the following notation.

\medskip

{\bf Notation.} Let $G$ be a group. For a subset $S$ in $G$ denote by $<<S>>$
the subgroup of $G$ generated by all the conjugacy classes of the elements of
$S,$
that is, $<<S>>$ is the minimal normal subgroup of the group $G$ which contains
$S.$
We also say that the subgroup $<<S>>$ is {\it normally generated} by $S.$

\medskip

\nlemma{({\bf Kaliman} \cite[Lemma 3.4]{Ka 2}).}
{\it The induced homomorphism
$(\si_C)_*\,:\,\pe(M^\prime)\to \pe(M)$ is an isomorphism.}

\label{lmm1}
\medskip

\proof
The restriction $\sigma_C \vert \,({M^\prime\setminus E^\prime})\,:\,
M^\prime\setminus E^\prime\longrightarrow M\setminus D$ is an isomorphism.
Thus, we may consider the following commutative diagram (left) and the induced
commutative triangle (right):

\begin{picture}(400,90)
\put(101,53){\vector(1,-1){30}}

\put(126,41){$\sigma_C$}

\put(13,10){$M\setminus D$}

\put(51,23){\vector(1,1){30}}

\put(45,41){{\small $i$}}
\put(83,59){$M^\prime$}
\put(63,12){\vector(1,0){60}}

\put(88,18){{\small $j$}}
\put(131,10){$M\,\,,$}

\put(313,53){\vector(1,-1){30}}
\put(313,53){\vector(1,-1){25}}

\put(333,41){$(\sigma_C)_{*}$}

\put(205,10){$\pe(M\setminus D)$}

\put(251,23){\vector(1,1){30}}
\put(251,23){\vector(1,1){25}}

\put(246,41){{\small $i_{*}$}}
\put(278,59){$\pe(M^\prime)$}
\put(267,12){\vector(1,0){58}}
\put(267,12){\vector(1,0){52}}

\put(288,18){{\small $j_{*}$}}
\put(331,10){$\pe(M)\,\,.$}

\end{picture}

\noindent It is easily seen that both $i_{*}$ and $j_{*}$ are surjections
(since
a complex hypersurface has real codimension 2). Thus, $(\sigma_C)_{*}$ is also
surjective. Denote by $\alpha_D$ a vanishing loop of $D.$
By Lemma \ref{fuji}.$(a)$, $\Ker j_{*}=\,<<\alpha_D>>\,.$ We choose
$\alpha_D$ in such a way that near $D$ it is a boundary circle of a small
transversal disc $\omega$ centered at a point $c_0\in C.$ Then the proper
transform
$\omega^\prime$ of $\omega$ in $M^\prime$ is a disc centered at a point of
$E^\prime=
E\setminus D^\prime.$ Thus, $i_{*}(\alpha_D)=1\in \pe(M^\prime),$ i.e.
$\alpha_D\in \Ker i_{*}.$ This implies that
$\Ker j_{*}\subset \Ker i_{*}.$ But since $j_{*}=(\sigma_C)_{*}\circ i_{*},$
$\Ker i_{*}=\Ker j_{*},$ and so $(\sigma_C)_{*}:\pe(M^\prime)\longrightarrow
\pe(M)$ is an isomorphism.
\qed

\medskip

\rem The proof can be word-by-word applied in a more
general case when the cener $C$ of the blow up is not necessarily smooth and
contained in the regular part $\reg D$ of $D;$
it is only important that it should not be contained in the singular
locus $\sing D$ of the hypersurface $D.$

\label{rempi}
\bigskip

\lemma (cf. Kaliman \cite[Proof of Thm. 3.5]{Ka 2}) {\it
Suppose that $(i)$ $D$ is a topological manifold, and $(ii)$ the varieties
$D$ and $C$ are acyclic. Then the variety $M'$ is acyclic iff $M$ is.}

\label{lmm2}
\medskip

\proof
As follows from Lemma \ref{lmm1}, $(\sigma_C)_*\,:\,H_1(M') \to H_1(M)$ is an
isomorphism (hereafter all the homology groups are with coefficients in $\Z$).
Note that

\smallskip

\noindent $\bullet$ $\sigma_C\,:\,E^\prime\longrightarrow C$
is a smooth fibration with a contractible fiber, and so, it yields a
homotopy equivalence between $E'$ and $C.$ Therefore,
$(\sigma_C)_*\,:\,H_*(E') \to H_*(C)$ is an isomorphism.
Hence, the exceptional divisor $E^\prime$ is also acyclic.

\smallskip

\noindent $\bullet$ Let $\dX$ be the one-point compactification of a manifold
$X.$
Then we have
$$\tH^i(\dX)\simeq H^i(\dX,*)\simeq H^i_c(X)\cong H_{m-i}(X)\,$$
where $\cong$ stands for the Lefschetz--Poincar{\'e} duality, and
$m=\dm_{\R}X.$ Thus, under our assumptions $\dD$ and $\dE^\prime$ are homology
spheres, and so is $M$ resp. $M'$ iff it is acyclic.

\smallskip

Assume first that $M$ is acyclic. Then by the exact cohomology sequence of
a pair

\begin{picture}(400,80)
\put(242,26){\vector(-1,1){30}}
\put(242,41){{\small $r^{*}$}}
\put(109,10){$H^{*}(\dD)$}
\put(176,56){\vector(-1,-1){30}}
\put(129,41){{\small $i^{*}$}}
\put(174,60){$H^{*}(\dM)$}
\put(165,12){\vector(1,0){60}}
\put(188,18){{\small $\delta^{*}$}}
\put(228,10){$H^{*}(\dM,\dD)\,\,,$}
\end{picture}

\noindent where deg $r^*=$ deg $i^*=0$, deg $\delta^{*}=1,$ we have
$H^{2n-1}(\dM,\dD)\simeq H^{2n}(\dM,\dD)\simeq \Z,$
where $n=\dm_{\C}M$,
and $\tH^{j}(\dM,\dD)=0$ for $j \le 2n-2$.
Since $\dM\setminus \dD=M\setminus D\approx M^\prime\setminus E^\prime,$
we have the homeomorphisms\footnote{More generally,
one can show that if $D$ is a
non-compact connected closed subspace of a smooth connected manifold $M$,
then the identity mapping of the complement $M \setminus D$ extends to
a homeomorphism of Hausdorff compact spaces
$\dM/\dD \stackrel{{\approx}}\longrightarrow (M\setminus D)\dot{\,}.$}
$$\dM/\dD\approx (M\setminus D)\dot{\,}\approx (M^\prime\setminus
E^\prime)\dot{\,} \approx \dM^\prime/\dE^\prime\,.$$
Hence,
$$H^{*}(\dM^\prime,\dE^\prime)\simeq \tH^{*}(\dM^\prime/\dE^\prime)
\simeq \tH^{*}(\dM/\dD)\simeq H^{*}(\dM,\dD).$$
Thus,  $H^{2n-1}(\dM^\prime,\dE^\prime)\simeq
H^{2n}(\dM^\prime,\dE^\prime)\simeq \Z,$
and the other groups are zero. From the exact cohomology sequence
of the pair $(\dM^\prime,\dE^\prime)$ we obtain
$H^j(\dM^\prime)\simeq H^j(\dE)=0,$ $1\le j\le 2n-3,$
and
$$0=H^{2n-2}(\dM^\prime,\dE^\prime)
\longrightarrow H^{2n-2}(\dM^\prime)\longrightarrow
H^{2n-2}(\dE^\prime)\simeq \Z\stackrel{\partial^{*}}
\longrightarrow $$
$$\longrightarrow H^{2n-1}(\dM^\prime,\dE^\prime)\simeq \Z\longrightarrow
H^{2n-1}(\dM^\prime)\longrightarrow
H^{2n-1}(\dE^\prime)=0. \eqno{(*)}$$

\noindent By the Poincar{\'e} duality, we have
$$H_{2n-j}(M^\prime)=\tH^j(\dM^\prime).$$
Hence, $H_i(M^\prime)=0$ for $i\ge 3,$ and $H_1(M^\prime)\simeq H_1(M)=0.$
Thus, by the Poincar{\'e} duality, $H^{2n-1}(\dM^\prime)=0$ in $(*),$ whence
$\partial^{*}:H^{2n-2}(\dE^\prime)\simeq \Z\longrightarrow H^{2n-1}
(\dM^\prime,\dE^\prime)\simeq \Z$ is onto, and so, it is an
isomorphism. This implies that $H^{2n-2}(\dM^\prime)=0,$
and also, by the Poincar{\'e} duality,  $H_2(M^\prime)=0$. Finally, we have
that $\tH_{*}(M^\prime)=0,$ which means that the variety $M^\prime$ is acyclic.

Vice versa, assuming that $M^\prime$ is an acyclic manifold,
one can prove that so is $M$
repeating word-in-word the above arguments, but exchanging the roles of the
pairs $(M,\,D)$ and $(M',\,E').$ This completes the proof.
$\Box$

\bigskip

\nthm{({\bf Kaliman} \cite[Thm. 3.5]{Ka 2}).} {\it Suppose that
$(i)$ $D$ is a topological manifold, and $(ii)$ $D$ and $C$ are acyclic. Then
$M^\prime$ is contractible iff $M$ is.}

\label{thkl}
\medskip

\proof By the Theorems of Hurewicz and Whitehead, $M$ resp. $M'$
is contractible iff it is acyclic and simply connected. Thus, the statement
follows immediately from Lemmas \ref{lmm1} and \ref{lmm2}.
\qed

\bigskip

\nlemma{({\bf Kaliman} \cite{Ka 2}).} {\it $\ok(M^\prime)\ge \ok(M).$}

\label{lmp2}
\medskip

\proof
Indeed, $M^\prime=\hM\setminus D^\prime$ implies $\ok(M^\prime)\ge
\ok(\hM)$ (see Theorem \ref{prr}.$(a)$). Since $\sigma_C:\hM\longrightarrow M$
is a proper birational morphism, by Theorem \ref{prr}.(e),
we get $\ok(M^\prime)\ge \ok(M),$ as claimed.
\qed

\bigskip

\ex (Kaliman \cite{Ka 2}) Let $X=S\times \C$ be an exotic $\C^3$ where
$S$ is a contractible surface of log-general type.
Choose a finite sequence of points $\{(s_i,z_i)\}_1^n\subset X$
(where the complex numbers $(z_i)_{1}^n$ are pairwise distinct) as the center
$C$ of the Kaliman modification $\si_C\,:\,X'\to X$ along the union of the
fibers
$H_i=S\times \{z_i\}$ of the second projection
$X\longrightarrow \C,\,\,i=1,\dots,n.$
Then $E^\prime_i\simeq \C^2,$ and one can show that
$E_i'\,, \,\,i=1,\dots,n,$ are the only copies of $\C^2$ in $X^\prime.$
Thus, by Kaliman's Theorem \ref{thkl}, $X'$ is an exotic $\C^3.$ It contains
precisely $n$ copies  $E_i'\,, \,\,i=1,\dots,n,$ of $\C^2,$
and their positions in $X'$
or, what is the same,  the positions of the points
$\{(s_i,z_i)\}_1^n\subset X$ up to automorphisms of $X$
provide deformation parameters.
\label{exkl}

\subsection{ Affine modifications}

Here we follow the recent paper \cite{KaZa}.

\medskip

\defin More generally, consider a triple $(A,\,I,\,f)$
where $A$ is an affine domain, i.e. a finitely generated integral domain over
$\C,$ $I$ is an ideal of $A,$ and $f \in I$ is a nonzero element.
By the {\it affine modification} of the domain $A$ along
the principal divisor $(f)$ with center $I$
we mean the affine domain $$A' = \Sigma_{I,\,f}(A) := A[It]/(1-ft)\,$$
where
$$A[It] := A \oplus \bigoplus_{n=1}^{\infty} (It)^n \simeq
A \oplus I\oplus I^2 \oplus \dots = Bl_I(A) $$
is the {\it blow up algebra}, or
the {\it Rees algebra} of the ideal $I$ (see e.g. \cite[\S 5.2]{Ei}).

The affine algebra $A$ resp. $A'$ coincides with the algebra
of regular functions $\C[X]$ (resp. $\C[X']$) on its spectrum
$X = {\rm spec}\,A$ resp. $X' = {\rm spec}\,A'$ which is an irreducible
affine variety. The variety $X' = \Sigma_{I,\,f}(X)$ is also called
the {\it affine modification} of the variety $X$
along the principal divisor $D_f$ with center $I.$

\medskip

\rem
It is easily seen that $X' = \hX \setminus D'_f$ where $\hX := Bl_I X$
is the blow up of the variety $X$ with center
$I,$ and $D'_f$ is the proper transform
in $\hX$ of the divisor $D_f,$ defined in an appropriate way
\cite[Proposition 1.1.a]{KaZa}. Thus, the Kaliman modification is
a particular case of the affine modification.
Moreover, we have the following theorem.

\medskip

\nthm{{\cite[Theorem 1.1]{KaZa}}.}
{\it Any birational morphism $X' \to X$ of
reduced irreducible affine varieties is an affine modification.}

\bigskip

\rem A choice of a system of generators $a_1,\dots,a_r$
of the algebra $A$ defines the proper embedding
$X \hookrightarrow \C^r_{\bx},\,\,x_i = a_i(x), \,\,i=1,\dots,r.$
If also a system of generators
$b_0=f,\,b_1,\dots,b_s$ of the ideal $I$ is given, then the formulas
$$x_i = a_i, \,\,i=1,\dots,r,\,\,\,y_j = b_jt,\,\,j=1,\dots,s,$$
yield a proper embedding $X' \hookrightarrow \C^{r+s}_{(\bx,\,\by)}.$
The blowup morphism $\si_I\,:\,X' \to X$ coincides with
the restriction to $X'$ of the natural projection
$\C^{r+s}_{(\bx,\,\by)} \to \C^r_{\bx}.$

\medskip

Recall such a notion.

\medskip

\defin A system of generators $b_0=f,\,b_1,\dots,b_s$ of the ideal $I$
is called {\it regular} if for any $j=0,\dots,s-1$
the image of the element $b_{j+1}$ in the quotient
algebra $A/(b_0,\dots,b_j)$ is not a zero divisor.

\medskip

\nprop{({\bf Davis} \cite{Dav}).} {\it Let $b_0=f,\,b_1,\dots,b_s$
be a regular system of generators of the ideal $I.$ Then the image
of the variety
$X'$ under the embedding $X' \hookrightarrow X \times \C^{s}_{\by}$
coincides with the subvariety given by the equations
$f(x)y_j = b_j(x),\,\,j=1,\dots,s.$}

\label{prDa}
\medskip

In Examples \ref{exCI} - \ref{exMS}
below we wright down explicit equations of certain affine modifications
by making use of Proposition \ref{prDa}.

\medskip

Let $E \subset \hX$
be the exceptional divisor of the blow up $\si_I\,:\,\hX \to X.$ The divisor
$E' := E \setminus D_f' = E \cap X'$ on $X'$ is also called the
{\it exceptional divisor} of the affine modification
$\si_I\,:\,X' \to X.$ Its image $\si_I(E')$ is contained in the subvariety
$C:=V(I) \subset {\rm supp}\,D_f.$ Denote
$\tau\,=\si_I\,|\,E'\,:\,E' \to C \hookrightarrow {\rm supp}\,D_f.$

\medskip

The next proposition is a generalization of Kaliman's Lemmas
\ref{lmm1} and \ref{lmm2} above. It provides a control
on preservation of the topology under affine modifications.

\medskip

\nprop{{(see \cite[Prop. 3.1 and Thm. 3.1]{KaZa}).}}
{\it Suppose that the following conditions $(i)-(iii)$ are fulfilled:

\smallskip

\noindent $(i)$ the affine varieties $X$ and $X'=\Sigma_{I,\,f}(X)$ are smooth;

\smallskip

\noindent $(ii)$ the divisors $E'$ and ${\rm supp}\,D_f$ are irreducible and
$E'=\si_I^*({\rm supp}\,D_f);$

\smallskip

\noindent $(iii)$ these divisors $E'$ and ${\rm supp}\,D_f$ are topological
manifolds.

\smallskip

\noindent Then the following statements hold:

\smallskip

\noindent $(a)$ the induced homomorphism
$(\si_I)_*\,:\,\pi_1(X') \longrightarrow \pi_1(X)$ is an isomorphism, and

\smallskip

\noindent $(b)$
the homomorphism $(\si_I)_*\,:\,H_*(X';\,\Z) \longrightarrow H_*(X;\,\Z)$
is an isomorphism iff the homomorphism
$\tau_*\,:\,H_*(E';\,\Z) \longrightarrow H_*({\rm supp}\,D_f;\,\Z)$ is.}

\label{ptop}
\bigskip

\rem Actually, under the condition $(ii)$ we have that
$\si_I(E') \cap \reg D_f \neq \emptyset.$
This allows us to apply Remark \ref{rempi} above to prove $(a).$

\bigskip

The next statement follows from  Proposition \ref{ptop} in the same way as
Theorem \ref{thkl} follows from Kaliman's Lemmas \ref{lmm1}, \ref{lmm2}.

\medskip

\nthm{{\cite[Corollary 3.1]{KaZa}.}}
{\it Under the conditions $(i)-(iii)$ of
Proposition \ref{ptop} the variety $X'$ is
contractible (resp. acyclic) iff the variety $X$ is.}

\label{thmST}
\bigskip

We give below several examples of application of this theorem.

\medskip

\ex {\bf An affine modification of the affine space along a divisor
with center at a codimension two complete intersection.} \label{exCI}
Let $A = \C[X]=\C^{[r]}$ be a polynomial algebra, i.e. $X = \C^r$
is an affine space, and set $I =(f,\,g)$
where $f,\,g \in A$ are two non-constant relatively prime polynomials.
Then $\{f,\,g\}$ is a regular system of generators of the ideal $I.$
In virtue of Proposition \ref{prDa}, the affine modification
$X'= \Sigma_{I,\,f} (X)$ is a hypersurface in $\C^{r+1}$
given by the equation
$f({\overline x})y - g({\overline x}) = 0$ where ${\overline x} =
(x_1,\dots,x_r).$
The blowup morphism $\sigma_I\,:\,X' \to X$ coincides with the restriction to
$X$ of the projection
$\C^{r+1} \to \C^{r},\,\,\,({\overline x},\,y) \to {\overline x}.$
The exceptional divisor
$E' \subset X'$ is given in $\C^{r+1}$ by the equations
$f({\overline x})=g({\overline x}) = 0;$ thus,
$E' \simeq C \times \C$ where $C = V(I)$ is the center of the blow up.

\medskip

\ex {\bf The Russell cubic threefold as affine modification.} \label{exKR}
In particular, set $A = \C[x,\,z,\,t]$ (i.e. $X = \C^3$),
$f = -x^2$ and $I =(f,\,g)$ where $g=x + z^2 + t^3.$
Consider the affine modification
$X' = \Sigma_{I,\,f} (X)$ along the divisor $D_f = 2D_x$ with center at the
ideal $I = (-x^2, \,x + z^2 + t^3) \subset \C^{[3]}$
supported by the affine plane curve
$C =V(I) = \Gamma_{2, \,3} :=\{x=z^{2}+t^{3}=0\}\subset \C^2.$
Then $X'$ is the smooth 3-fold $x + x^2y + z^2 + t^3=0$ in $\C^4,$ which has
been called in the Introduction the {\it Russell cubic}
(see also \cite{Ru 1} and Examples \ref{ex8.1} and \ref{ex8.3} below).
It birationally dominates the affine space $X=\C^3$ under the blowup morphism
$\si_I \,:\,X' \to X = \C^3,\,\,\,\si_I \,:\, (x,\,y,\,z,\,t)
\longmapsto (x,\,z,\,t).$
The exceptional divisor $E'$ coincides with the `book-surface'
$B:=\{x=0\} \subset X',\,\,\,B \simeq \C\times \Gamma_{2, \,3}.$
It is easily seen that the conditions $(i)-(iii)$ of Proposition \ref{ptop}
are fulfilled. Therefore, by Theorem \ref{thmST}, the Russell cubic
$X'$ is contractible. Moreover, by the Dimca--Ramanujam Theorem \ref{thdr},
it is diffeomorphic to $\R^6.$ However, as we will see in $\S 7$ below,
the Russell cubic is not isomorphic to the affine space $\C^3,$ and whence
it is an exotic $\C^3.$

\medskip

\ex {\bf The tom Dieck--Petrie surfaces as affine modifications} \cite{tDP 1,
tDP 2}. \label{exDP}
Recall (see Example \ref{ex3.1} above)
that these are smooth surfaces $X_{k,\,l}$ in $\C^3$
defined by the polynomials
$$p_{k,l}= {(xz + 1)^k - (yz + 1)^l - z \over z} \in \C[x,\,y,\,z]\,$$
where $k,\,l \ge 2,\,\,\,(k,\,l)=1.$
Actually, one can see that
$X_{k,\,l} = \Sigma_{\rho,\,\Gamma_{k,\,l}} (\C^2)$
where $\Gamma_{k,\,l} := \{x^k - y^l = 0\} \subset \C^2$ and
$\rho = (1,\,1) \in \Gamma_{k,\,l}$ (cf. \cite[Example 2.2]{KaZa}).
Moreover, the conditions $(i)-(iii)$ of
Proposition \ref{ptop} hold. In view of Theorem \ref{thmST},
the surface $X_{k,l}$ is contractible.

\medskip

\ex \label{exMS}
Let $M = f^*(0),\,\,\,f \in \C^{[n]},$ be
a smooth reduced irreducible hypersurface in $\C^n,$
and let
$f_1,\dots,f_k \in \C^{[n]}$ be polynomials without common zeros on $M.$
Consider the smooth affine varieties
$$D := M \times \C^k \subset \C^{n+k}_{(\bx,\,\bu)}=\C^n_{\bx}
\times \C^k_{\bu}\,\,\,\,\,\,\,\,{\rm and}\,\,\,\,\,\,\,
C := \{f(\bx) = 0 = g(\bx,\,\bu)\} \subset D \subset\C^{n+k}_{(\bx,\,\bu)}\,$$
where $g(\bx,\,\bu) := \sum_{i=1}^k u_if_i \in \C^{[n+k]}.$
It is easily seen that the natural embeddings
$M \times {\overline 0} \hookrightarrow C \hookrightarrow D$ provide hoomotopy
equivalences.
Therefore, by Theorem \ref{thmST}, the affine modification
$$X' := \Sigma_{C,\,D} (\C^{n+k}) = \{f(\bx)v - g(\bx,\,\bu) = 0\} \subset
\C^{n+k+1}_{(\bx,\,\bu,\,v)}$$ of the affine space $\C^{n+k}$
along the hypersurface $D$ with center $C$
is a smooth contractible hypersurface in $\C^{n+k+1}.$

Since the hypersurface $M$ is assumed being smooth one may take e.g.
$k=n$ and $f_i = {\partial f \over \partial x_i}, \,\,i=1,\dots,n;$
then we have $C \simeq TM$ and $D \simeq T\C^n\,|\,M.$

\subsection{The hyperbolic modification}

Here we follow, up to minor changes, tom Dieck \cite{tD 1} (cf. another
treatment in Petrie \cite{Pe}). We restrict the consideration to
the simplest possible case.

\bigskip

\defin Let $h\in \C[x_1,\ldots,x_n]$ be an irreducible polynomial such that
$h(\overline 0)=0.$
Suppose that $\grad_{\overline 0}\,h\ne \overline 0,$ and so
the hypersurface $X=\{h=0\}\subset \C^n$ is smooth at the origin.
Define  the {\it hyperbolic modification} $q$ of $h$ as follows:
$$q(\ox,u)=\frac{h(u\ox)}{u}\in \C[x_1,\ldots,x_n,u].$$
Since $h(u\ox)=uq(\ox,u),$ we have the equalities
$$u\frac{\partial q}{\partial u}(\ox,u)+q(\ox,u)=\sum_1^nx_i\frac{\partial
h(u\ox)}
{\partial x_i}\,,$$
$$\frac{\partial q}{\partial x_i}(\ox,u)=\frac{\partial h(u\ox)}
{\partial x_i}, \qquad i=1,\ldots,n.$$
It follows that, once $(\ox_0,u_0)$ is a critical point of $q,$ i.e.
$\grad_{(\ox_0,u_0)}\,q={\overline 0},$ then also
$\grad_{\ox_0u_0}h={\overline 0}=h(u_0\ox_0),$
that is, $u_0\ox_0\in X$ is a singular point of the hypersurface $X,$ and
$(\ox_0,u_0)\in Y_0:=\{q=0\}\subset \C^{n+1}.$ Thus, all
the fibers $Y_c=\{q=c\},$ $c\in \C^{*}=\C\setminus \{0\}$,
of the polynomial $q$ are smooth hypersurfaces, and the
fiber $Y_0=\{q=0\}$ is smooth iff so is $X,$ which will be assumed in the
sequel.

\bigskip

\lemma {\it The restriction $q\,\vert \,{(\C^{n+1}\setminus Y)}
:\C^{n+1}\setminus Y\longrightarrow \C^{*}$ is a trivial algebraic fiber
bundle with the fiber $Y_1:=\{q=1\}.$}

\label{lmh1}
\medskip

\proof
Consider the commutative triangle

\begin{picture}(400,90)
\put(201,53){\vector(1,-1){30}}

\put(226,41){{\small $q$}}

\put(115,10){$Y_1\times \C^{*}$}

\put(150,23){\vector(1,1){30}}

\put(145,41){{\small $\Phi$}}
\put(164,60){$\C^{n+1}\setminus Y_0$}
\put(163,12){\vector(1,0){60}}

\put(188,18){{\small pr}}
\put(231,10){$\C^{*}\,\,$}
\end{picture}

\noindent where the mapping $\Phi$ is defined as follows:
$$(\oy,\lambda) := ((\ox,u),\,\lambda) \stackrel{\Phi}\longmapsto
(\lambda\ox,\lambda^{-1}u):=\oy_\lambda\in Y_\lambda.$$
It is easy to check that $\Phi$ is a fibrewise (biregular) isomorphism,
so we are done.
$\Box$

\bigskip

Define a $\C^{*}-$action on $\C^{n+1}: (\lambda, (\ox,u))\stackrel{G_\lambda}
\longrightarrow (\lambda\ox,\lambda^{-1}u),$ $\lambda\in \C^{*}.$
Then $$q(G_\lambda (\ox,u))=\frac{h(u\ox)}{\lambda^{-1}u}=\lambda q(\ox,u).$$
That means that $q$ is a quasi-invariant of weight $1$ of the $\C^{*}-$action
$G.$
In particular, the hypersurface $Y_0$ is invariant with respect to $G,$ and
$G_\lambda(Y_c)=Y_{\lambda c}.$
In the above diagram to the action $G$ there corresponds the canonical
$\C^{*}-$action
on the direct product, whence $\Phi$ is an equivariant morphism.

The monomials $ux_1,\ldots,ux_n\in \C[x_1,\ldots,x_n,u]$ are $G-$invariants.
Moreover, they generate the algebra of $G-$invariants
$\C[x_1,\ldots,x_n,u]^G=\C[ux_1,\ldots,ux_n].$ Hence,
the {\it algebraic quotient} of $\C^{n+1}$
by this $\C^{*}-$action is isomorphic to
the affine space $\C^n:$
$$\C^{n+1}//G\simeq \C^n=\spec \C[x_1,\ldots,x_n,u]^G.$$
The $\C^{*}-$action $G$ on $\C^{n+1}$ is
{\it hyperbolic}\footnote{This is why the modification was called {\it
hyperbolic.}},
that is, it has only one fixed point
(the origin ${\overline 0} \in \C^{n+1}$), and the weights $(1,\ldots,1,-1)$
of the action $G$ at the origin are of
different signs. The origin belongs to the closure of each
$G-$orbit which is contained either in the hyperplane $\{u=0\}$
or in the axis $OU:= \{{\bar 0}\} \times \C;$ all the other
$G-$orbits are closed.

Denote by $M$ the complement of the axis $OU$ in $\C^{n+1}.$
Then the $\C^{*}-$action $G$ restricts to $M$ with closed orbits only.
Let $\pi:M\longrightarrow T$ be the canonical morphism onto the orbit space,
or the {\it geometric quotient}, $T=M/G.$ Also, consider the morphism
$\tau:M\longrightarrow \C^n,$ $(\ox,u)\longmapsto u\ox.$ Since $\tau$
is constant on any orbit, it factors as $\tau=\sigma\circ \pi$:

\begin{picture}(400,90)
\put(221,53){\vector(1,-1){30}}

\put(246,41){{\small $\sigma$}}

\put(147,10){$M$}

\put(171,23){\vector(1,1){30}}

\put(165,41){{\small $\pi$}}
\put(205,59){$T$}
\put(183,12){\vector(1,0){60}}

\put(208,18){{\small $\tau$}}
\put(251,10){$\C^n\,$}
\end{picture}

\noindent The restriction of the morphism $\pi$ to the
hypersurface $M \bigcap \,\{u=0\}:= \hE\simeq \C^n\setminus \{\overline 0\}$
coincides with the standard projection
$\C^n\setminus \{0\}\longrightarrow \Pp^{n-1},$ $\ox \longmapsto \{\lambda
\ox\}_
{\lambda \in \C^{*}}.$ Set $\pi(\hE)=E\subset T.$ Then
$\sigma(E)=\{\overline 0\},$ i.e. $E$ is the exceptional divisor
of $\sigma;$ it is straightforward that $\sigma \vert\,(T\setminus E)\,
:\,T\setminus E\longrightarrow \C^n\setminus \{0\}$ is an isomorphism.
It is easily seen that $\sigma:T\longrightarrow \C^n$ is the blow up of the
origin.

Furthermore, $\pi(Y_0\bigcap M):=X^\prime$ is the proper transform of $X$
in $T.$ Indeed, $Y_0$ is saturated by the orbits, whence $\pi(Y_0\bigcap M)$
is an irreducible closed hypersurface in $T$ which contains the proper
transform
$\sigma^\prime(X).$
Therefore, $T\setminus X^\prime$ is the Kaliman transform of the affine space
$\C^n$ along
$X$ with center at the origin ${\overline 0} \in \C^n.$

\bigskip

\lemma {\it There is an isomorphism $Y_1\simeq T\setminus X^\prime.$}

\label{lmh2}
\medskip

\proof
Fix a point $y=(\ox,u)\in Y_1.$ Since $q$ is a $G$--quasi-invariant of weight
$1,$ $G_\lambda (y)\in Y_\lambda,$ and hence, the orbit $Gy$ of $y$
meets the hypersurface $Y_1$ at the point $y$ only. This means that the
morphism
$\pi \vert\,{Y_1}:Y_1\longrightarrow T\setminus X^\prime$ is injective.
On the other hand, any $G-$orbit outside the hypersurface $Y_0$ meets the
$q-$fibre
$Y_1;$ thus, this morphism is also surjective. Finally, a bijective morphism of
smooth varieties is an isomorphism.
\qed

\bigskip

\cor {\it The hypersurface $Y_1 \subset \C^{n+1}$ is isomorphic to the Kaliman
modification of  the affine space $\C^n$ along the
hypersurface $X\subset \C^{n}$ with center at the origin.}

\label{corh1}
\bigskip

\lemma {\it  The hypersurface $Y_0 \subset \C^{n+1}$ is isomorphic to the
Kaliman modification of the product $X\times \C$ along the
hypersurface $X\times\{0\}$ with center at the point $(\overline 0,0)\in
X\times \C.$}

\label{lmh3}
\medskip

\proof
The morphism $$Z^\prime:= \C^{n+1}\stackrel{\sigma}
\longrightarrow \C^{n+1}=: Z, \quad (y_1,\ldots,y_n,u)\longmapsto
(uy_1,\ldots,uy_n,u)\,,$$ is nothing but the affine modification
of the variety $Z$ along the
hyperplane $H_0:=\{u=0\}$ with center at the origin, and with
the exceptional divisor $E^\prime=\{u=0\}\subset Z^\prime.$
Consider the natural embedding
$i\,:\,X\times \C\hookrightarrow Z.$ Set $\hh(\ox,u)=h(\ox);$ then the
image of $i$ is the hypersurface $\hh=0$ in $Z\simeq \C^{n+1}.$

We have $$\hh\circ \sigma(\oy,u)=h(u\oy)=uq(\oy,u).$$
Hence, $\hh\circ \sigma(\oy,u)=0$ for any point $(\oy,u)\in Y_0$
(i.e. such that $q(\oy,u)=0$).
Thus, $\sigma(\oy,u)\in X\times \C,$
and so $\sigma(Y_0)\subset X\times \C.$ Furthermore, the total
preimage of the product $X\times \C$ in $Z^\prime$ is the union of
the hypersurface $Y_0$ and of the exceptional divisor $E'=\{u=0\}.$
Therefore, $Y_0$ is the proper
transform of the variety $X\times \C$ in $Z^\prime,$
and the assertion follows.
\qed

\bigskip

\rem The $\C^{*}-$action $\lambda(y_1,\ldots,y_n,u)=
(\lambda y_1,\ldots,\lambda y_n,\lambda^{-1}u)$ on $Z^\prime$
provides the $\C^{*}-$action
$\lambda(x_1,\ldots,x_n,u)=
(x_1,\ldots,x_n,\,\lambda^{-1}u)$ on the affine space $Z$
and on the product $X\times \C.$

\bigskip

\exr
(\nexr) \label{exr6.1} Show that under the embedding
$\C^n\hookrightarrow \C^{n+1}$ given as $\ox\longmapsto (\ox,1)$ the variety
$\,X$ is naturally isomorphic to the hyperplane section
$Y_0\bigcap H_1$ where $H_1:=\{u=1\} \subset \C^{n+1}.$
Furthermore, show that the exceptional divisor $E' \subset Y_0$ of the Kaliman
transform
$\sigma:Y_0\longrightarrow X\times \C$ coincides with the linear subspace
$Y_0\bigcap \,\{u=0\}$ in $\C^{n+1}.$ Let $\sigma^\prime\,:\,Y_0\longrightarrow
X$
be composed of the contraction $\sigma$ and the first projection. Verify that
$$\sigma^\prime:(\ox,u)\longmapsto (\frac{\ox}{u},1)\in Y_0\bigcap H_1\simeq
X$$
outside the exceptional divisor $E',$ and
$\sigma^\prime(\ox,0)=(\overline 0,1)$ on $E'.$ Deduce that the hypersurface
$Y_0$ is the closure in $\C^{n+1}$ of the $\C^{*}-$orbit of the subvariety
$Y_0\bigcap H_1\simeq X.$

\bigskip

\nthm{({\bf tom Dieck} \cite{tD 1}).} {\it Let $X\subset \C^n$ be a smooth
contractible
hypersurface given by an irreducible polynomial $h\in \C[x_1,\ldots,x_n]$
where $h(\overline 0)=0.$ Then any fiber $Y_c=q^{-1}$(c)$,$ $c\in\C,$
of the hyperbolic modification $q(\ox,u)=\frac{h(u\ox)}{u}\in
\C[x_1,\ldots,x_n,u]$ of the polynomial $h$ is a smooth contractible
hypersurface in
$\C^{n+1}.$
Thus, $q:\C^{n+1}\longrightarrow \C$ yields a foliation of $\C^{n+1}$
by smooth contractible hypersurfaces.}

\medskip

\proof
Indeed, by Lemma \ref{lmh1}, $Y_c\simeq Y_1$
for any $c\ne 0.$ By Corollary \ref{corh1} and Lemma \ref{lmh3},
the hypersurfaces $Y_0$ and $Y_1$ are both
Kaliman modifications
of triples of smooth contractible varieties.
By Kaliman's Theorem \ref{thkl}, these hypersurfaces
$Y_0$ and $Y_1$ are contractible.
\qed

\bigskip

\rems \nrem The inequality $\ok(M^\prime)\ge \ok(M)$ of
Lemma \ref{lmp2} does not help to answer to the question whether the
contractible affine varieties
$Y_0$ and $Y_1$ are exotic $\C^n-$s. Indeed, we have
$\ok(\C^n)=\ok(X\times \C)=-\infty$ (see Theorem \ref{prr}$(a)$).
However, in certain cases
the {\it intermediate Eisenman-Kobayashi intrinsic measures}
serve as appropriate analytic invariants
(see Kaliman \cite{Ka 2}).

\medskip

\nrem The Kaliman Theorem \ref{thkl} is still applied if
$X$ is only assumed being a contractible topological manifold smooth at the
origin. In that case we still have that all the
hypersurfaces $Y_c,\,c \neq 0,$ are smooth and contractible, but the central
fiber $Y_0$ can be singular, as it is in the following example.

\bigskip

\ex {\bf The foliations of the affine spaces with contractible leaves
arising from the tom Dieck-Petrie surfaces} \cite{tDP 2}
(see Examples \ref{ex3.1}, \ref{exDP} above). \label{ex6.1}
Recall (see \cite{LiZa} and Theorem \ref{cusuthm} above)
that up to automorphisms of $\C^2$ the only irreducible simply connected
singular
affine plane curves are the curves
$\,\,\Gamma_{k,l}:=\{x^k-y^l=0\} \subset \C^2,$
$(k,l)=1,$ $k>l\ge 2.$
Starting with $\Gamma_{k,l}$ perform the hyperbolic modification at
the smooth point
$(1,1)\in \Gamma_{k,l}.$ We obtain a foliation
$p_{k,l}:\C^3\longrightarrow \C$ of $\C^3$ by the
fibers of the polynomial
$$p_{k,l}:=\frac{(xz+1)^k-(yz+1)^l}{z}\in \C[x,y,z].$$
All of them are irreducible contractible surfaces;
all but the central one $p^{-1}_{k,l}(0)$ are smooth. One can see
that $\ok(X_{k,l})=1$ where $X_{k,l}:=p^{-1}_{k,l}(1)$
(see Exercise \ref{exr6.2} below).
Next,
starting with $X_{k,l},$ by means of hyperbolic modifications one can construct
non-trivial foliations of $\C^4,$ $\C^5,$ etc. by smooth contractible
hypersurfaces. Moreover, the corresponding polynomials are quasi-invariants
of hyperbolic $\C^{*}-$actions on $\C^n.$ In particular, for $n \ge 4$
the zero fiber of such a polynomial
is a smooth contractible hypersurface in $\C^n$ endowed with a hyperbolic
$\C^{*}-$action. Furthermore, new contractible affine hypersurfaces
can be obtained by passing
to cyclic $\C^{*}-$coverings over these ones (see the next section).

\bigskip

\exr \label{exr6.2} (\nexr)
Verify that $\ok(X_{k,l})=1$ when $(k,l)=1,$ $k>l\ge 2$.

\medskip

{\it Hint.} One can proceed, for instance, as follows.
Lifting the meromorphic function $x^k/y^l$ on $\C^2$ to the function
$f:= (xz+1)^k/(yz+1)^l\,\vert \,X_{k,l}$ on $X_{k,l}$
(see Exercise \ref{exDP}) we obtain a
$\C^*-$fibration $f\,:\,X_{k,l} \to \Pp^1$. Hence, by Iitaka's Easy Addition
Theorem \ref{prr}$(c)$, $\ok(X_{k,l})\le 1.$ Since $X_{k,l}$
is acyclic, by the Classification Theorem \ref{clthm}$(b)$, $\ok(X_{k,l})=1$ as
soon as
$X_{k,l} \not\simeq \C^2.$ Recall that the surface $X_{k,l}$ is the
Kaliman modification of $\C^2$ along the curve $\Gamma_{k,l} \subset \C^2$
with center at the point $(1,1)\in \Gamma_{k,l}$.
Resolving singularities of the plane projective curve ${\overline \Gamma_{k,l}}
\cup l_{\infty} \subset \Pp^2$ and blowing up
at the point $(1,1)\in {\overline \Gamma_{k,l}}$,
we obtain an SNC-completion $(V_{k,l},\,D_{k,l})$ of $X_{k,l}.$
Contracting, if necessary, the $(-1)-$ components of the boundary
divisor $D_{k,l}$ of valence at most
two of the dual graph $\Gamma_{D_{k,l}}$
we come to a minimal completion $(V_{k,l}^{\rm min},\,D_{k,l}^{\rm min})$ of
the open surface $X_{k,l}.$
The dual graph of its boundary divisor $D_{k,l}^{\rm min}$ is
non-linear
(what is this graph?). Therefore, by the Ramanujam Theorem \ref{ram},
$X_{k,l} \not\simeq \C^2.$

\section{Cyclic $\C^{*}-$coverings}

\ndefin{(cf. \cite[Prop. 2.11]{KoRu 2}).}
Let $X$ be an affine variety, and let $q\in\C[X]$ be
a regular function on $X.$ For an integer $s>1$ set
$Y_s=\{(x,u)\in X\times \C \;\vert\; q(x)=u^s\}.$ The projection
$\varphi_s\,:\, Y_s\longrightarrow X,$ $(x,u)\stackrel{{\varphi_s}}\longmapsto
x,$
yields a {\it cyclic covering of $X$ branched to order $s$ along the
principal divisor $F_0=q^{*}(0)$ }.
We suppose that $X$ is a smooth affine variety and $F_0$ is a smooth reduced
divisor on $X;$ then the variety $Y_s$ is also smooth
(indeed, $\grad_{(x,u)}(q(x)-u^s)=(\grad_xq,\,-su^{s-1})$), as well as
the hypersurface $F_{s,\,0} := \varphi_s^{-1}(F_0)$ in $Y_s.$

Let $X$ be endowed with a regular $\C^{*}-$action $t:\C^{*}\times
X\longrightarrow X.$
 Suppose that the regular function $q$ is a quasi-invariant of $t$
of weight $d,$ i.e.
$$q(t_\lambda x)=\lambda^d q(x)$$
where
$\,d\in \Z.$ Then the $\C^{*}-$action
$\lambda(x,u):=(\lambda^s(x),\,\lambda^d u)$
on $X\times \C$ restricts to $Y_s$ making the
following commutative diagram equivariant

\begin{picture}(200,95)
\unitlength0.2em
\put(66,25){$Y_s$}
\put(99,25){$X$}
\put(74,27){$\vector(1,0){21}$}
\put(66,5){$\C$}
\put(74,6){$\vector(1,0){21}$}
\put(99,5){$\C\,\,$}
\put(69,22){$\vector(0,-1){11}$}
\put(60,16){${\rm pr_2}$}
\put(101,22){$\vector(0,-1){11}$}
\put(104,16){{\rm q}}
\put(78,9){{\small $u\mapsto u^s$}}
\put(73,30){{\small $(x,u)\mapsto x$}}
\end{picture}

\noindent where the original $\C^{*}-$action $G$ on $X$ is replaced by its
's-th power'
$(\lambda,x)\longmapsto \lambda^s (x):=t(\lambda^s,\,x).$ Indeed, for
$(x,u)\in Y_s$ we have:
$$q(\lambda^s(x))=\lambda^{sd}q(x)=\lambda^{sd}u^s=(\lambda^du)^s,$$
whence $(\lambda^s(x),\,\lambda^du)\in Y_s,$ which shows that the above diagram
is equivariant.

If, in addition, $(d,s)=1,$ then the monodromy of the cyclic covering
$\varphi_s\,:\, Y_s \to X$ is represented via
the action on $Y_s$ of the subgroup $\omega_s\subset \C^{*}$ of the
$s-$th roots of unity. Indeed, since $(s,\,d)=1$
the $\omega_s-$orbit of a point $(x,u)$ in $Y_s$
is
$$\omega_s(x,u)=\{(x,\lambda^du)\;\vert\; \lambda^s=1\}=\varphi_s^{-1}(x)\,.$$
The fixed point set
$Y_s^{\omega_s} = \{(x,u)\in Y_s\;\vert\; u=0\}=F_{s,\, 0} \subset Y_s$
of the monodromy action on the variety
$Y_s$ can be
identified with the hypersurface $F_0 \subset X.$ Thus,
we get $X=Y_s/\omega_s$ with the quotient action of $\C^{*}/\omega_s\simeq
\C^{*}$  on $X.$

The equivariant covering $Y_s\longrightarrow X$ as above is called a
{\it cyclic $\C^{*}-$covering.}

\medskip

\rems \nrem \label{vrem} The action of the monodromy group
$\omega_s\simeq \Z/s\Z$ on $Y_s$ is homologically trivial. Indeed,
this is so for the action  on $Y_s$ of the connected group
$\C^{*}\supset \omega_s.$

\medskip

\nrem \label{vvrem} The above observations are equally applied in the more
general
setting when the regular $\C^{*}-$action is only given on the Zariski open
subset $X^* := X \setminus F_0$ of $X.$ In particular, if $(d,s)=1,$ then the
monodromy group $\omega_s$ of the cyclic covering
$\varphi_s\,:\, Y^*_s \to X^*$ where
$ Y^*_s := Y_s \setminus \vi_s^{-1}(F_0),$ acts
trivially in the homology $H_*(Y^*_s;\,\Z).$

\bigskip

The following result provides a generalization of
Theorem A in\footnote{Exposing this result in
\cite[Thm. 6.9]{Za 5},
the condition ($\sharp_1$) below has been missed.
In the proof of Theorem \ref{thm7.3}
below it replaces Lemma 6.8 in \cite{Za 5}
which is wrong.} \cite{Ka 1}.

\bigskip

\nthm{({\bf Kaliman}).} {\it Let $X$ be a smooth contractible
affine variety, and let the principal divisor $F_0=q^{*}(0)\subset X$ where
$q\in \C[X]$ be smooth,
reduced and irreducible. Denote $G = \pi_1(X \setminus F_0),$ and fix a
vanishing loop $\alpha = \alpha_{F_0}\in G$ of the divisor $F_0.$
Suppose that the following conditions are fulfilled.

\smallskip

\noindent ($\sharp$) The regular function $q$ is a quasi-invariant of weight
$d\ne 0$ of a regular $\C^{*}-$action defined on $X \setminus F_0.$

\smallskip

\noindent ($\sharp_1$) For some integer $c\ne 0,$ $\alpha^c$ is an element of
the center $Z(G)$ of the group $G.$

\smallskip

\noindent ($\sharp_2$) For an integer $s > 0$ such that $(s,c)=(s,d)=1,$
the hypersurface $F_0$
is $\Z_p-$acyclic\footnote{hereafter $\Z_p= \Z/p\Z$.}
for each prime divisor $p$ of $s.$

\smallskip

\noindent Consider the cyclic
covering $\varphi_s\,:\,Y_s \to X$ branched to order $s$ along $F_0$.
Then $Y_s$ is a smooth contractible affine variety.}
\label{kalim}
\medskip

Due to the Theorems of Hurewicz and Whitehead, it is enough to show that
$Y_s$ is acyclic and simply connected. This is done, respectively,
in Theorems \ref{thm7.2} and \ref{thm7.3} below. Notice that the conditions
$\sharp$ and $\sharp_2$ guarantee acyclicity of $Y_s$ whereas the condition
$\sharp_1$ provides its simply connectedness.

\subsection{Acyclicity of cyclic $\C^*-$coverings}

\nthm{({\bf Kaliman} \cite{Ka 1}; {\bf tom Dieck} \cite{tD 2}).} {\it Let $X$
be an
acyclic smooth affine variety, and $F_0=q^{*}(0)$ where $q\in \C[X],$ be a
smooth reduced irreducible principal divisor on $X.$ Consider a cyclic
covering $\varphi_s\,:\,Y_s \to X$ branched to order $s$ along $F_0$  where
$(s,d)=1.$
Suppose that the following conditions hold.

\smallskip

\noindent ($\sharp$) The regular function $q$ is a quasi-invariant of weight
$d\ne 0$ of a regular $\C^{*}-$action defined on $X \setminus F_0.$

\smallskip

\noindent ($\sharp_2$) The hypersurface $F_0$
is $\Z_p-$acyclic for each prime divisor $p$ of $s.$

\smallskip

\noindent Then $Y_s$ is acyclic, too.}
\label{thm7.2}
\medskip

Before proving Theorem \ref{thm7.2} we recall the
Smith theory (see \cite[Ch.III]{Bre}).

\medskip

{\bf Elements of Smith's Theory.}
Consider a finite simplicial polyhedron
$Y$ endowed with a simplicial action of a finite group $\om.$
Usually, passing to the second barycentric subdivision one obtains
certain additional regularity properties of the action, which
are always to be assumed (see \cite[III.1]{Bre}).
Let $k$ be a field, and let
$\Z[\om],\, \Z_p[\om], \,k[\om]$ be the group rings of $\om$
(e.g. $\Z[\om]=\{\sum_{g\in \om}n_gg\;\vert\; n_g\in \Z\}$
with natural ring operations). The simplicial chain complexes
$C(Y),\, C(Y)\otimes \Z_p,\,C(Y)\otimes k$ are, respectively,
$\Z[\om]-,\, \Z_p[\om]-, \,k[\om]-$modules (indeed, given a simplex
$\delta$ of $Y$
we have
$$(\sum n_gg)(\delta)=\sum n_gg(\delta)\in C(Y)).$$
In the sequel $\om$ is assumed to be a finite cyclic group
$\Z_s=\Z/s\Z$ acting on $Y$ in such a way that the fixed point
set $Y^{\om}$ of the $\om-$action on $Y$ coincides with the individual
fixed point set $Y^g$ for every $g\in \om,$ $g\ne e.$ In particular,
the $\om-$action on the complement $Y\setminus Y^{\om}$ is free.
We denote by $X=Y/\om$ the orbit space, by
$\pi:Y\longrightarrow X$ the natural projection, and we identify
$Y^{\om}$ with its image in $X.$ Consider the following three
homomorphisms of chain complexes:
$$\pi_{*}:C(Y)\longrightarrow C(X),$$
$$\sigma: C(Y)\longrightarrow C(Y),\qquad \sigma=\sum_{g\in \om}g\in
\Z[\om]\,,$$
$$\mu=\pi^{*}:C(X)\longrightarrow C(Y), \qquad
\mu(\delta)=\pi^{-1}(\delta)\,\,\,\,{\rm if}\,\,\,\,\delta \cap Y^{\om} =
\emptyset; \qquad \mu(\delta)=\sigma(c)\,\,\,\,{\rm if}\,\,\,\,\,\pi_{*}(c) =
\delta\,.$$
Note that the induced homomorphism $\pi_{*}$ is surjective.
Then we have \cite[III.2]{Bre}
$\Ker \pi_{*}=\Ker \sigma,$ and there is an isomorphism
$$\sigma C(Y)\simeq C(Y)/\Ker \sigma\simeq C(Y)/\Ker \pi_{*}=C(X),$$
whence $\mu\pi_{*}=\sigma.$ But $\pi_{*}\mu(c)=\vert \om \vert {\pi_{*}(c)}.$
On the homology level, this leads to the following assertions.

\bigskip

\lemma \cite[III(2.2), (2.3)]{Bre} {\it
$$\pi_{*}\mu_{*}=\vert \om\vert: H_{*}(X)\longrightarrow H_{*}(X),$$
$$\mu_{*}\pi_{*}=\sigma_{*}=\sum_{g\in \om}g_{*}:H_{*}(Y)\longrightarrow
H_{*}(Y).$$
Here $\mu_{*}$ is called a {\rm transfer}. On the invariant part of homology
we have
$$\mu_{*}\pi_{*}\vert {H_{*}(Y)^{\om}}=\vert
\om\vert:H_{*}(Y)^{\om}\longrightarrow
H_{*}(Y)^{\om}.$$}
This implies such a corollary.

\bigskip

\ncor{\cite[III(2.4)]{Bre}.} {\it If $k$ is a field of characteristic
$\chr k=q$ with $(q,\vert \om\vert)=1,$ then the restriction
$$\pi_{*}\vert {H_{*}(Y;k)^{\om}}:H_{*}(Y;k)^{\om}\longrightarrow H_{*}(X;k)$$
is an isomorphism, and its inverse is the transfer $\mu_{*}.$ Moreover,
$$H_{*}(Y;k)=\mu_{*}H_{*}(X;\,k)\oplus \Ker \pi_{*},$$
where $\Ker \pi_{*}=\Ker \sigma_{*}.$}

\label{corb1}
\bigskip

\cor {\it Suppose that ${\om}$ acts trivially in the homology, i.e.
$${\om}_{*}\vert {H_{*}(Y)}={\rm id}\,.$$ Then for any field $k$ with
$(\chr k,\,\vert {\om}\vert)=1$ we have the isomorphism
of transfer
$$\pi_{*}=\mu^{-1}_{*}:H_{*}(Y;k)\stackrel{\simeq}\longrightarrow H_{*}(X;k).$$
In particular, if ${\om}\simeq \Z_q$ where $q$ is a prime number, then
the elements of the kernel and of the cokernel
of the homomorphism
$\pi_{*}:H_{*}(Y;\Z)\longrightarrow H_{*}(X;\Z)$ are torsions of
order $q.$}

\label{corb2}
\medskip

The last assertion follows by the Universal Coefficient Formula:
$$\tH_j(Y;\Z_q)=\tH_j(Y;\Z)\otimes\Z_q\oplus \Tor (H_{j-1}(Y;\Z);\Z_q)\qquad
\forall j.$$

\bigskip

\ndefin{\cite[III.3]{Bre}.}
In what follows ${\om}\simeq \Z_p$ is a multiplicative cyclic group of
prime order $p$ with a generator $t\in {\om},$ so that
$\sigma=1+t+\ldots+t^{p-1}\in \Z_p[{\om}].$ We set $\tau=1-t\in \Z_p[{\om}].$
We have $t^p=1,\,\,\sigma \tau=\tau \sigma=0,$ and $\sigma=\tau^{p-1}$
(indeed, $(-1)^i {p-1 \choose i} \equiv 1 \,\,\bmod{p}$).
For an element $\rho=\rho_i:=\tau^i\in \Z_p$ set $\overline \rho=\tau^{p-i};$
then $\overline \sigma=\tau$ and $\overline \tau=\sigma.$ Given
$\rho=\rho_i=\tau^i,$
consider the chain complex $\rho C(Y;\Z_p).$ Its graded homology group
$H^\rho_{*}(Y;\Z_p):= H_{*}(\rho C(Y;\Z_p))$ is called the
{\it special Smith's homology group}.

\bigskip

There are the exact sequences of chain complexes with coefficients in $\Z_p$
\cite[III(3.1),(3.8)]{Bre}:
$$0\longrightarrow \overline \rho C(Y)\oplus C(Y^{\om})
\stackrel{i}\longrightarrow
C(Y)\stackrel{\rho}\longrightarrow\rho C(Y)\longrightarrow 0\,,$$
$$0 \longrightarrow \sigma C(Y) \stackrel{i}\longrightarrow \tau^j C(Y)
\stackrel{\tau}\longrightarrow \tau^{j+1} C(Y) \longrightarrow 0,\,\,\,j =
1,\dots,p-1\,.$$
Besides, the kernel of the homomorphism
$\sigma\, : \,C(Y;\,\Z_p) \to C(Y;\,\Z_p)$
and those of the composition
$C(Y;\,\Z_p) \to C(Y,\, Y^{\om};\,\Z_p) \to  C(X,\, Y^{\om};\,\Z_p)\,$
where $Y^{\om}$ is indentified with its image in $X,$
are the same \cite[p. 124]{Bre}. These observations lead to the following

\bigskip

\nprop{\cite[III(3.3),(3.4),(3.8)]{Bre}.}
{\it For the homology groups with $\Z_p$ coefficients, one has

\smallskip

\noindent $(a)$ an isomorphism $H^\sigma_{*}(Y)\simeq H_{*}(X;Y^{\om})$,

\smallskip

\noindent and the following two {\rm Smith's exact homology sequences}:

\smallskip

\noindent $(b)$

\begin{picture}(400,80)
\put(270,26){\vector(-2,1){60}}
\put(254,41){$i_{*}$}

\put(105,10){$H_{*}^\rho (Y)$}
\put(176,56){\vector(-1,-1){30}}

\put(140,41){$\rho_{*}$}

\put(174,60){$H_{*}(Y)$}
\put(165,12){\vector(1,0){60}}
\put(188,18){{\small $\delta_{*}$}}

\put(230,10){$H_{*}^{\overline \rho}(Y)\oplus H_{*}(Y^{\om})$\,\,}
\end{picture}

\smallskip

\noindent $(c)$

\begin{picture}(400,80)
\put(247,26){\vector(-1,1){30}}
\put(244,41){$i_{*}$}

\put(105,10){$H_{*}^{\rho_{j+1}} (Y)$}
\put(176,56){\vector(-1,-1){30}}

\put(140,41){$\tau_{*}$}

\put(174,60){$H_{*}^{\rho_j}(Y)$}
\put(165,12){\vector(1,0){60}}
\put(188,18){{\small $\delta_{*}$}}

\put(230,10){$H_{*}^{\sigma} (Y)\,\,$}
\end{picture}

\noindent where $\dg \rho_{*}=\dg \tau_{*} =\dg i_{*}=0,\dg \delta_{*}=-1.$}

\label{prop3}
\bigskip

\prop {\it Suppose that\footnote{The numeration of the conditions that
we use here agrees with those in the next Corollary \ref{cor3}
and Exercise \ref{exr7.1}.}

\smallskip

\noindent $(ii)$ the fixed point set $Y^{\om}$ is non-empty and
$\Z_p$--acyclic:
$\tH_{*}(Y^{\om};\,\Z_p) = 0\,,$ and

\smallskip

\noindent $(iii)$ the quotient $X = Y/{\om}$ is $\Z_p$--acyclic:
$\tH_{*}(X;\,\Z_p) = 0\,.$

\smallskip

\noindent Then also $Y$ is $\Z_p$--acyclic:
$\tH_{*}(Y;\,\Z_p) = 0\,.$}

\label{prop4}
\medskip

\proof In view of the vanishing
$$\tH_{*}(X;\Z_p)=\tH_{*}(Y^{\om};\Z_p) = 0\,,$$ from the
exact homology sequence of a pair
$$\dots\stackrel{i_{*}}\longrightarrow H_{j}(X;\Z_p) \stackrel{r_{*}}
\longrightarrow H_{j}(X,\,Y^{\om};\Z_p)\stackrel{\delta_{*}}\longrightarrow
H_{j-1}(Y^{\om};\Z_p)\stackrel{i_{*}}\longrightarrow\dots$$
it follows that $H_{*}(X,\,Y^{\om};\Z_p)=0$, and thus,
by Proposition \ref{prop3}$(a)$, also $H^\sigma_{*}(Y;\Z_p) = 0.$
Therefore, by Smith's exact sequence $(c)$,
$H^\rho_{*}(Y;\Z_p)=0\,\forall \rho = \rho_j,\,j=1,\dots,p-1\,.$
Now, by Smith's exact sequence $(b)$,
$\tH_{*}(Y;\Z_p)\simeq H^\rho_{*}(Y;\Z_p) = 0.$
$\Box$

\bigskip

\cor {\it Suppose that the following conditions hold:

\smallskip

\noindent $(i)$ ${\om} \simeq \Z_p$ acts trivially in homology:
${\om}_{*}\vert {H_{*}(Y)}={\rm id}\,,$

\smallskip

\noindent $(ii)$ the fixed point set $Y^{\om}$ is non-empty and
$\Z_p$--acyclic:
$\tH_{*}(Y^{\om};\,\Z_p) = 0\,,$ and

\smallskip

\noindent $(iii)$ the quotient $X = Y/{\om}$ is acyclic:
$\tH_{*}(X;\,\Z) = 0\,.$

\smallskip

\noindent Then also $Y$ is acyclic:
$\tH_{*}(Y;\,\Z) = 0\,.$}
\label{cor3}

\medskip

\proof
By Corollary \ref{corb2}, $H_{*}(Y;\,\Z_q) \simeq H_{*}(X;\,\Z_q) = 0\,$
for any prime $q \neq p$. By Proposition \ref{prop4}, also $H_{*}(Y;\,\Z_p)=
0\,$.
Thus, by the Universal Coefficient Formula, $\tH_{*}(Y;\,\Z)\otimes\Z_q = 0$
for all prime $q$. Then $\tH_{*}(Y;\,\Z)= 0$.
\qed

\bigskip

\exr (\nexr) \label{exr7.1}
Assume that ${\om}\simeq \Z_s$ acts on $Y$ in such a way that

\noindent $(0)$ $Y^g=Y^{\om}\ne \emptyset$ for every $g\in {\om},$ $g\ne e;$

\noindent $(i)$ the action is homologically trivial, i.e.
${\om}_{*}\vert {H_{*}(Y;\Z)}={\rm id};$

\noindent $(ii)$ the fixed point set $Y^{\om}$ is $\Z_p-$acyclic for
any prime divisor $p$ of $s;$

\noindent $(iii)$  the quotient $X=Y/{\om}$ is acyclic: $\tH_{*}(X;\,\Z) =
0\,.$

\noindent Show that $Y$ is acyclic, too:
$\tH_{*}(Y;\,\Z)= 0$.

\bigskip

\rem Assume that the $\C^*-$action in Theorem
\ref{thm7.2} is regular on the whole variety $X.$ Then the monodromy group
action on
the covering variety $Y_s$
is homologically trivial (see Remark \ref{vrem} above),
that is,
the above condition $(i)$ is fulfilled. This provides a proof of
Theorem \ref{thm7.2} in that particular case. Notice that the
assumptions of smoothness of $X$ and $F_0$ are not used in
this proof.

In general case,
following tom Dieck \cite{tD 2}, we need to consider
branched coverings over smooth varieties and to make use of the Thom classes.
We recall below their definition and some properties
(see e.g. \cite[VIII.11]{Do}, \cite[\S  9, \S 10]{MilSta}).

\bigskip

{\bf Thom's classes and Thom's isomorphisms.}
Consider an oriented connected smooth real manifold
$X$ and a codimension $2$ closed oriented submanifold  $F_0$ of $X.$
Let $N \to F_0$ be the (oriented) normal
bundle of $F_0$ in $X$ with the zero section $Z_0 \simeq F_0.$ Fix
a tubular neighborhood
$U \subset X$ be of the submanifold $F_0$ in $X$ such that the pair
$(U,\,F_0)$ is diffeomorphic to the pair $(N,\,Z_0).$ Denote
$U^* := U \setminus F_0$ and
$N^* := N \setminus Z_0.$ By excision, we have the
isomorphisms $\tH_{*}(X, \,X^*;\,\Z) \simeq \tH_{*}(U, \,U^*;\,Z) \simeq
\tH_{*}(N, \,N^*;\,\Z),$ and similarly for the cohomology groups. The {\it
Thom class} $t(F_0) \in H^2(X, \,X^*;\,\Z)\simeq H^2(N, \,N^*;\,\Z)\,$ is a
unique cohomology class which takes the value $1$ on any
oriented relative two-cycle $(F,\,F^*) \in H_2(N, \,N^*;\,\Z)$ defined by
a fiber $F$ of the normal bundle $N.$

The cap-product
with the Thom class $t(F_0) \in H^2(X, \,X^*;\,\Z_q)$
yields the {\it Thom isomorphism}\footnote{the homology groups with negative
indices are considered
being zero.}
$$\tH_i(X, \,X^*;\,\Z_q) \simeq H_{i-2}(F_0;\,\Z_q),\,\,i=0,\,1,\dots\,.$$

Let $\vi_s\,:\,Y_s \to X$ be a smooth cyclic ramified covering of $X$
branched to order $s$ along $F_0,$ i.e. $Y_s$ is an
oriented manifold equipped with an action of a group ${\om} \simeq \Z_s$ of
orientation preserving diffeomorphisms;
the fixed point set $Y_s^{\om}\subset Y_s$ is a codimension $2$ closed oriented
submanifold; ${\om}$ acts
freely in the complement $Y_s \setminus Y_s^{\om},$ and
$\vi_s\,:\,Y_s \to X$ is the orbit map which
provides a natural identification of $Y_s^{\om}$ with the submanifold
$F_0 \subset X$.

\bigskip

Note that under the
assumption ($\sharp$) of Theorem \ref{thm7.2} the monodromy group
${\om} \simeq \Z_s$
acts trivially in the homology $H_*(Y_s \setminus F_0;\,\Z).$ Thus,
the next proposition yields Theorem
\ref{thm7.2} (cf. Corollary \ref{cor3}).

\bigskip

\nprop{(see \cite[Thm. 2.9]{tD 2}).}
{\it Let in the notation as above $\vi_s\,:\,Y_s \to X$ be a smooth cyclic
ramified covering of a smooth manifold $X$ branched to order $s$ along
a codimension $2$ submanifold $F_0 \subset X.$ Suppose that

\smallskip

\noindent $(i)$ the covering group ${\om} \simeq \Z_s$ acts trivially in the
homology of the complement $Y_s \setminus Y_s^{\om}\,:$
$\,\,{\om}_{*} \vert {H_{*} (Y_s \setminus Y_s^{\om};\Z)} = {\rm id};$

\smallskip

\noindent $(ii)$ the fixed point set $Y_s^{\om}$ is $\Z_p-$acyclic for
any prime divisor $p$ of $s;$

\smallskip

\noindent $(iii)$ the quotient $X=Y_s/{\om}$ is acyclic: $\tH_{*}(X;\,\Z) =
0\,.$

\smallskip

\noindent Then the manifold $Y_s$ is acyclic, too: $\tH_{*}(Y_s;\,\Z)= 0$.}

\label{p5}
\medskip

\proof
Assume, for simplicity, that $s = p$ is a prime number; like in Exercise
\ref{exr7.1} above
the general case can be reduced to this one.
By Proposition \ref{prop4}, we have $\tH_*(Y_p;\,\Z_p)=0$.
By the Universal Coefficient Formula, it suffices to prove that
$\tH_*(Y_p;\,\Z_q)=0$ (i.e. $Y_p$ is $\Z_q-$acyclic) for any prime
$q \neq p$.

Denote $Y_p^* = Y_p \setminus Y_p^{\om}$ and
$X^* = X \setminus F_{0}.$ The restriction $\pi\,:\,Y_p^*\to X^*$ is a
non-ramified cyclic covering of order $p$.
By Corollary \ref{corb1},
$(\vi_s)_*\,:\,H_*(Y_p^*;\,\Z_q) \to H_*(X^*;\,\Z_q)$ is an isomorphism for any
prime $q \neq p.$

We have the Thom isomorphisms
$$\tH_i(X, \,X^*;\,\Z_q) \simeq H_{i-2}(F_0;\,\Z_q),\,\,\,\,\,\,\,\,
\tH_i(Y_p, \,Y_p^*;\,\Z_q) \simeq H_{i-2}(Y_p^{\om};\,\Z_q),$$
given by cap-products
with the Thom classes $t(F_0) \in H^2(X, \,X^*;\,\Z_q)\,$ resp.
$\,\,t(Y_p^{\om}) \in H^2(Y_p, \,Y_p^*;\,\Z_q).$ It is easily seen that
$(\vi_s)^*(t(F_0)) = p\cdot t(Y_p^{\om}).$ Since the multiplication by $p$
is an invertible operation in $\Z_q-$(co)homology for $q \neq p,$
it follows that
$(\vi_s)_*\,:\,H_*(Y_p, \,Y_p^*;\,\Z_q)\longrightarrow H_*(X, \,X^*;\,\Z_q)$
is an isomorphism.

Consider the following commutative diagram where the horizontal lines are
exact homology sequences of pairs with $\Z_q-$coefficients:

\begin{picture}(200,95)

\put(0,65){$\dots\longrightarrow H_{j+1}(Y_p,\,Y_p^*)
\longrightarrow H_{j}(Y_p^*)\longrightarrow H_{j}(Y_p)
\longrightarrow H_{j}(Y_p,\,Y_p^*) \longrightarrow H_{j-1}(Y_p^*)
\longrightarrow \dots$}

\put(77,50){$\vector(0,-1){25}$}
\put(82,35){$\simeq$}
\put(152,50){$\vector(0,-1){25}$}
\put(157,35){$\simeq$}
\put(216,50){$\vector(0,-1){25}$}
\put(287,50){$\vector(0,-1){25}$}
\put(292,35){$\simeq$}
\put(357,50){$\vector(0,-1){25}$}
\put(362,35){$\simeq$}

\put(3,5){$\dots\longrightarrow H_{j+1}(X,\,X^*)
\longrightarrow H_{j}(X^*)\longrightarrow H_{j}(X)
\longrightarrow H_{j}(X,\,X^*) \longrightarrow H_{j-1}(X^*)
\longrightarrow \dots$}

\end{picture}

\noindent By the above observations, we may conclude that the four indicated
vertical arrows are isomorphisms induced by
the projection $\vi_s$. By the 5-lemma, the middle vertical arrow is
an isomorphism, too.
Hence, since $X$ is acyclic, $\tH_*(Y_p;\,\Z_q)\simeq \tH_*(X;\,\Z_q) = 0$
for any prime $q.$ This yields the assertion.
\qed

\medskip

Thus, the proof of Theorem \ref{thm7.2} is completed.

\bigskip

\ex {\bf The acyclic surfaces $\bf Y_{k,\,l,\,s}$ in $\C^3$} \label{ex7.1}.
Let $X$ be a smooth acyclic surface, $F_0 = q^*(0)$
be a smooth reduced irreducible simply connected curve in\footnote{see Theorem
\ref{cusuthm} above for a description of such pairs
$(X,\,F_0).$} $X$ where $q\in \C[X]$
is a quasi-invariant of weight $d \neq 0$ of a regular $\C^*-$action on $X
\setminus F_0$. Then by Theorem \ref{thm7.2},
$Y_s := \{z^s = q(x)\}\subset X \times \C,$ where
$(d,\,s) = 1,$ is a smooth acyclic surface, too.

\smallskip

For instance, for $k, \,l,\,s$ pairwise relatively prime
the surface $Y_{k,\,l,\,s}\subset \C^3$ given
by the equation $$\frac{(xz^s+1)^k-(yz^s+1)^l}{z^s}=1$$
is a smooth acyclic one, and $\ok(Y_{k,\,l,\,s}) = 1.$
Indeed, there is a cyclic
$\C^*-$covering $Y_{k,\,l,\,s} \to X_{k,\,l}$ over the tom Dieck-Petrie surface
$X_{k,\,l}$ (see Examples \ref{exDP}, \ref{ex6.1} and Exercise \ref{exr6.2}
above)
branched to order $s$ along the
curve $L_{k,\,l} := X_{k,\,l} \cap \{z=0\} \simeq \C$ in $X_{k,\,l}.$
The $\C^*-$action in $X_{k,\,l} \setminus L_{k,\,l}$ is induced via the
isomorphism
$X_{k,\,l} \setminus L_{k,\,l} \simeq \C^2 \setminus {\Gamma}_{k,\,l}$
by the linear $\C^*-$action $(\lambda,\,(x,\,y)) \longmapsto
(\lambda^l x,\,\lambda^k y)$ on $\C^2.$ Thus, we may apply Theorem
\ref{thm7.2} to show that  the surface $Y_{k,\,l,\,s}$ is acyclic.

\subsection{Simply connectedness of cyclic $\C^{*}-$coverings}

\nthm{({\bf Kaliman}).} {\it Let $X$ be a simply connected smooth affine
variety,
and let $F_0=q^{*}(0),$ where $q\in \C[X],$ be a smooth
reduced irreducible principal divisor in $X.$ Fix a vanishing loop
$\alpha=\alpha_{F_0} \in G:= \pi_1(X \setminus F_0)$ of the divisor
$F_0.$ Consider a cyclic
covering $\varphi_s\,:\,Y_s \to X$ branched to order $s$ along $F_0.$
Assume that

\smallskip

\noindent ($\sharp_1$) For an integer $c\ne 0$ such that $(s,\,c) = 1,$
$\alpha^c$ is an element of the center $Z(G)$ of the group $G.$

\smallskip

\noindent Then $Y_s$ is simply connected, too.}
\label{thm7.3}
\bigskip

\rem In \cite[Lemmas 7 and 8]{Ka 1} conditions on a polynomial
$q \in \C[x_1,\dots,x_n]$ are given which ensure that
$\pi_1(\C^n \setminus F_0) \simeq \Z.$ In particular, repeating word-in-word
the proof of Lemma 8 in \cite{Ka 1} (based on the Seifert-van Kampen Theorem)
one can easily see that $\pi_1(X \setminus F_0) \simeq \Z$ if $F_0$ is a
generic fibre of a regular function $q$ on a simply connected smooth affine
variety $X,$ that is, the restriction of $q$ onto a preimage
$q^{-1}(\Delta_{\epsilon})$ of a small disc $\Delta_{\epsilon} \subset \C$
centered at the origin yields a  (trivial) smooth fibre bundle over
$\Delta_{\epsilon}.$
Thus, in this case also the assumption ($\sharp_1$) of Theorem \ref{thm7.3}
holds.

\bigskip

We need the following definition.

\medskip

\defin
We say that a subgroup $H$ of a group $G$ is {\it normally
generated by elements} $a_1,\dots, a_n \in H$ if it is generated by the set
of all elements conjugate with $a_1,\dots, a_n,$ i.e. if $H$ is the minimal
normal subgroup of $G$ which contains $a_1,\dots, a_n.$ We denote it by
$\,<<a_1,\dots, a_n>>.$ $G$ is said to be {\it normally one--generated} if
$G = \,<<a>>\,$ for some element $a \in G.$

\bigskip

\lemma {\it Let $X$ be a smooth irreducible affine variety,
and let $F_0=q^{*}(0),$ where $q\in\C[X],$ be a reduced irreducible principal
divisor in $X.$  Fix a vanishing loop
$\alpha =\alpha_{F_0}\in G:= \pi_1(X \setminus F_0)$ of the divisor $F_0.$
Then the following statements hold.

\smallskip

\noindent $(a)$  $\pi_1 (X) = {\bf 1}$ iff $G = <<\alpha>>\,.$

\smallskip

\noindent $(b)$ Let $\varphi_s\,:\,Y_s \to X$ be a cyclic
covering  branched to order $s$ along $F_0.$  Set
${\widehat G}_s = <<\alpha^s>>\,.$
Assume that $F_0$ is a smooth divisor. Then $\pi_1(Y_s) = {\bf 1}$
iff $G/{\widehat G}_s \simeq \Z/s\Z.$}

\label{lmk2}
\bigskip

\proof $(a)$ By Lemma \ref{fuji}$(a)$, we have that Ker$\,(i_{*}\,:\,\pi_1 (X
\setminus F_0) \to \pi_1 (X)) = <<\alpha>>$, and the assertion follows.

\smallskip

\noindent $(b)$ Set $F_{s,0} = \varphi_s^{-1}(F_0) \subset Y_s,$
$X^* = X \setminus F_0$ and $Y^*_s = Y_s \setminus F_{s,0}.$
Then $\varphi_s\,:\,Y^*_s \to X^*$
is a non--ramified cyclic covering of order $s$. The induced homomorphism
$$(\varphi_s)_*\,:\,\pi_1(Y^*_s) \to \pi_1(X^*) =:G$$
is an injection onto a normal subgroup $G_s$ of $G$ of index $s,$ and
$G/G_s \simeq \Z/s\Z.$ Observe that by $(a)$,
$G = <<\alpha>>,$ and that
$\alpha^s \in G_s$ is the image of a vanishing loop
$\beta \in \pe(Y^*_s)$ of the smooth irreducible divisor
$F_{s,0}\subset Y_s\,,$ i.e. $\,(\varphi_s)_* (\beta) = \alpha^s.$
Therefore, $\,{\widehat G}_s:=<<\alpha^s>> \,\subset G_s,$  and
$\,{\widehat G}_s = G_s$ iff $G/{\widehat G}_s\simeq G / G_s \simeq \Z/s\Z.$

Denote also ${\widehat{\widehat G}}_s  = \,<<\alpha^s>>_{G_s}$
the subgroup of $G_s$ normally generated (in $G_s$)
by the element $\alpha^s\in G_s.$

\medskip

 {\bf Claim.} ${\widehat{\widehat G}}_s = {\widehat G}_s.$

\medskip

{\bf Proof of the claim.} Clearly,
${\widehat{\widehat G}}_s \subset {\widehat G}_s \subset G_s.$ Since the
quotient $G/G_s \simeq \Z/s\Z$ is Abelian we have that
$K := [G,\,G] \subset G_s.$ Since $G = \,<<\alpha>>\,$ the abelianization
$G_{\rm ab} := G/K$ is a cyclic group generated
by the class $K\alpha$ of the vanishing loop $\alpha.$
Hence, any element $g \in G$ can be written as $g = g'\alpha^t$ where
$g' \in K\subset G_s$ and $t \in \Z.$ Thus, we have
$g\alpha^s g^{-1} = g'\alpha^s g'^{-1} \in {\widehat{\widehat G}}_s$
for any $g \in G.$ Therefore,
${\widehat G}_s \subset {\widehat{\widehat G}}_s,$
and the claim follows.
$\Box$

\medskip

By $(a)$, $\pi_1(Y_s)=\id$ iff $\pi_1(Y^*_s) = \,<<\beta>>\,,$ or, what is the
same, iff ${\widehat{\widehat G}}_s = G_s.$ Due to the above Claim, the latter
holds iff ${\widehat G}_s = G_s,$ or, equivalently, iff
$G/{\widehat G}_s \simeq \Z/s\Z.$ This proves $(b)$.
$\Box$

\bigskip

{\bf Proof of Theorem \ref{thm7.3}.} Since $G = \,<<\alpha>>\,,$ any element
$g \in G$ can be written as $g = \prod_{i=1}^n g_i\alpha^{r_i}g_i^{-1}\,,$
where $g_i \in G$ and $r_i \in \Z,\,\,i=1,\dots,n.$ Let $\rho\,:\,G \to G/K
\simeq H_1(X \setminus F_0;\,\Z) \simeq \Z$ be the canonical surjection.
Then, clearly, $\rho(\alpha) = 1,$ and so,
$\rho(g) = \sum_{i=1}^n r_i \in \Z.$

Since $K \subset G_s := (\vi_s)_*(\pi_1(Y_s^*))$ and
$G/G_s \simeq \Z/s\Z,$ we have that $\rho(G_s) = s\Z.$ That is,
$g = \prod_{i=1}^n g_i\alpha^{r_i}g_i^{-1} \in G_s$ iff
$\rho(g) = \sum_{i=1}^n r_i \equiv 0 (\mod s).$

Using the assumption $(s,\,c) = 1$ write
$r_i = k_is + l_ic$ where $k_i,\,l_i \in \Z,\,\,i=1,\dots,n.$
By our assumption ($\sharp_1$), $\alpha^c \in Z(G).$ Hence,
$g_i\alpha^{r_i}g_i^{-1} =
g_i\alpha^{sk_i}g_i^{-1}\alpha^{cl_i},\,\,i=1,\dots,n,$ and furthermore,
$$g = \prod_{i=1}^n g_i\alpha^{r_i}g_i^{-1} = \left(\prod_{i=1}^n
g_i\alpha^{sk_i}g_i^{-1}\right)\alpha^{mc}\,$$
where $m = \sum_{i=1}^n l_i.$ For an element $g \in G_s$ it follows that
$\rho(g) = s\sum_{i=1}^n k_i + mc \equiv 0 \,(\mod s),$ or, equivalently,
$m \equiv 0\, (\mod s).$ Set $m = ls,\,l \in \Z.$ Whence, we have
$g = \left(\prod_{i=1}^n g_i\alpha^{sk_i}g_i^{-1}\right)(\alpha^s)^{lc}
\in {\widehat G}_s.$  Therefore,
${\widehat G}_s \subset G_s \subset {\widehat G}_s,$ and so,
${\widehat G}_s = G_s,$ as required (see Lemma \ref{lmk2}$(b)$).
$\Box$

\medskip

Now the proof of Kaliman's Theorem \ref{kalim} is completed.
In Exercises \ref{exr7.2} - \ref{exr7.7} below we expose some additional
properties
of the fundamental group $G = \pi_1(X \setminus F_0)$ in the situation where
the variety $X \setminus F_0$ is equipped with a $\C^*-$action. After that, in
Example\footnote{These exercises and example were elaborated in
\cite{KaZa}.}\ref{ex7.2} we show that
without the assumption ($\sharp_1$) (or, perhaps, a weaker one which has to
be precised) the fundamental group of a cyclic $\C^*-$covering $Y_s$ of a
contractible smooth affine variety (even surface) $X$ can be quite big.

\bigskip

\exrs
(\nexr) \label{exr7.2}\footnote{See \cite[Appendix]{Za 5}. A shorter proof
was suggested by H. Flenner.} Let $G = \,<<\alpha>>\,$
be a normally one-generated group. Denote by $K = [G,\,G]$ the commutator
subgroup of the group $G.$ Show that
$\alpha^c \in Z(G)$ iff \footnote{For a group $G$ and two subsets
$A, \,B \subset G$ we denote by $[A, \,B]$ the
subgroup generated by all the commutators $[a,\,b] = aba^{-1}b^{-1}$ where
$a \in A,\,b \in B.$}$[a^c,\,K] = {\bf 1},$ and that under this condition
$K \subset {\widehat G}_s := \,<<\alpha^s>>\,$ for any $s\in \Z$ prime to $c.$

\smallskip

(\nexr)  \label{exr7.3} Let $X$ and $q \in \C[X]$ be as in Theorem
\ref{thm7.3};
in particular, $\pi_1 (X) = {\bf 1}.$
Assume that $(\sharp')$ the restriction
$q\,|\,(X \setminus F_0)\,:\,X \setminus F_0 \to \C^*$
is a smooth fiber bundle with a connected fiber
$F_1 := q^{-1}(1)$. Show that $q_*(\alpha) = 1 \in \Z,$
$i_*\pi_1(F_1) = K :=[G,\,G],$
$G_{\rm ab} := G/K \simeq H_1 (X \setminus F_0;\,\Z) \simeq \Z,$
and $q_* = \rho\,:\,G \to G/K = \Z$ is the canonical surjection.
Deduce that $G = \pi_1(X \setminus F_0) \simeq \Z$ if and only if
$\pi_1(F_1) = {\bf 1},$ which in turn implies the condition ($\sharp_1$).

\smallskip

(\nexr) \label{exr7.4} Show, furthermore, that under the condition
$(\sharp)$ of Theorem \ref{kalim} the above assumption $(\sharp')$ is
fulfilled,
and, furthermore, the group $G$
contains a normal subgroup $G_d$ of index $d$ with the cyclic quotient
$G/G_d \simeq \Z/d\Z$ such that $G_d \simeq K \times \Z \simeq \pi_1(F_1)
\times \Z .$ Let an element $\gamma \in G_d $
correspond to a generator of the second factor $\Z$
of this decomposition.
Verify that\footnote{In the proof of Lemma 6.8 in \cite{Za 5} it was taken
$\gamma = \alpha^{d}.$ In general, this is not true;
see the next exercise.} $\gamma\alpha^{-d} \in K,$ and that
the centralizer subgroup $C_{\gamma}$ of $\gamma$ in $G$ contains $G_d.$

\smallskip

\noindent {\it Hint.} Put $G_d = (\vi_d)_*(\pi_1(Y_d^*))$ where
$\vi_d\,:\,Y_d \to X$ is the $d-$fold branched cyclic covering, $X^* = X
\setminus F_0$ and $Y_d^* = Y_d \setminus \vi_d^{-1}(F_0).$ The induced
$\C^*-$action on $Y_d^*$ yields an equivariant isomorphism
$F_1 \times \C^* \to Y_d^*,\,\,\,(x,\,\mu) \longmapsto (t_{\mu}(x),\,\mu),$
which provides, in turn, the desired decomposition of the subgroup $G_d$ of
$G.$ The element $\gamma$ corresponds to the image of a generator of the group
$\pi_1(\C^*) \simeq \Z$ under the homomorphism induced by the mapping
$\C^* \to {\cal O}_{x_0}\subset X^*$ onto the $\C^*-$orbit ${\cal O}_{x_0}$
of a base point $x_0 \in F_1.$

\smallskip

(\nexr)  \label{exr7.5} Set
$X = \C^2$ and $q(x,\,y) = x^2 - y^3 \in \C[x,\,y].$ Show that the group
$G = \pi_1(X \setminus F_0)$ can be identified with the 3-braid group
$$B_3 := <\sigma_1,\,
\sigma_2\,|\,\sigma_1\sigma_2\sigma_1=\sigma_2\sigma_1\sigma_2>\,,$$
the generators $\sigma_1,\, \sigma_2 \in G$ being vanishing loops of
the divisor $F_0 = \Gamma_{2,\,3} := q^*(0) \subset X.$
Describe the subgroups $G_d$ and $K=[G,\,G]$ in this example
(see the preceding exercise).
Verify that one can take for $\gamma$ the
element $(\sigma_1\sigma_2)^3$ which generates the center
$Z(G) \simeq Z(B_3)\simeq \Z$ of the braid group
$B_3.$ Putting $\alpha = \sigma_1$ check that $G = <<\alpha>>,$ and
$\alpha^c \notin Z(G)$ whatever $c \in \Z \setminus \{0\}$ is.

\smallskip

\noindent {\it Hint.} One can use the presentation of $-q$
as the discriminant of the cubic polynomial $t^3 -(y/{\sqrt[3]{4}})t +
x/{\sqrt{27}} \in \C[t].$ Consider, further, the
Vieta covering $\C^2 \to \C^2$
which is a branched Galois covering ramified
over $\Gamma_{2,\,3}$ with the symmetric group $S_3$ as the Galois group.

\smallskip

(\nexr) \label{exr7.6} Let $X$ be an irreducible quasi-projective variety,
$D \subset X$ be an irreducible hypersurface which contains the singularity
locus sing$X$ of $X,$ and let $C \subset D$ be a non-empty smooth subvariety
such that $C \cap ($sing$X \cup $sing$D) = \emptyset.$ Consider the Kaliman
modification $\sigma_C\,:\,X' \to X$ of the variety $X$ along $D$ with center
$C.$
Show that $\pi_1(X') \simeq \pi_1(X \setminus $sing$X).$

\smallskip

\noindent {\it Hint.} Replace the triple $(X,\,D,\,C)$ by the
triple $(X\setminus $sing$\,X,\,D\setminus $sing$\,X,\,C),$ and then
apply Lemma \ref{lmm1}.

\smallskip

(\nexr) \label{exr7.7} For a polynomial
$q \in \C[x_1,\,\dots,x_n]$ such that $q({\overline 0}) = 0$ denote
$h_q = q\circ \sigma_n (x_1,\,\dots,x_n) / x_n$ where
$$\sigma_n\,:\,\C^n \to \C^n,\,\,\,\sigma_n (x_1,\,\dots,x_n) :=
(x_1x_n,\,\dots,x_{n-1}x_n,\,x_n)\,$$
is the affine modification of the affine space
$\C^n$ along the hyperplane $\{x_n = 0\}$
with center at the origin (cf. the proof of Lemma \ref{lmh3}).
Put $X = q^{-1}(0) \subset \C^n$ and $X' = h_q^{-1}(0) \subset \C^n.$
Suppose that ${\overline 0} \in X$ is a smooth point. Verify that
$X'$ is the  affine modification of the variety $X$ along
the divisor $D:= X \cap \{x_n = 0\}$ with center
at the origin. Wright down an explicit equation of the hypersurface
$X' := \Sigma_{{\overline 0},\,D}(X) \subset \C^n$ in the general case
(see \cite[Example 2.1]{KaZa}).

\bigskip

\ex\label{ex7.2}\footnote{We are thankful to V. Sergiescu for useful
discussions related to this example.} Consider the smooth surface
$X_{k,l,s,m} := p_{k,l,s,m}^{-1}(0) \subset \C^3$
defined by the polynomial
$$p_{k,l,s,m}:= {(xz^m + 1)^k - (yz^m + 1)^l - z^s \over z^m} \in
\C[x,\,y,\,z]\,$$ where $0 \le m \le s.$ For $m>0$
this surface is smooth, and for $m=0$ it is smooth
outside of the point $P_0 = (1,\,1,\,0) \in X_{k,l,s,0}.$
Suppose that $(k,\,l) = (k,\,s) = (l,\,s) = 1.$ Then the surface
$X_{k,l,s,s} = Y_{k,l,s}$ is acyclic
(see Example \ref{ex7.1})\footnote{Actually, by
Lemma \ref{lmm2}, all the surfaces
$X_{k,l,s,m},\,\,m=1,\dots,s,$ are acyclic,
too.}.
Moreover, it can be presented as a cyclic
$\C^*-$covering over the
smooth contractible tom Dieck-Petrie
surface
$X_{k,l} = X_{k,l,1,1}\subset \C^3$ (see Examples \ref{ex3.1},
\ref{exDP} and \ref{ex6.1}
above)
branched to order $s$ along the curve
$L_{k,l} := X_{k,l} \cap \{z=0\} \simeq \C$ in $X_{k,l}.$ However,
the smooth acyclic surface $Y_{k,l,s}$ is not contractible
and possesses quite a big
fundamental group which we describe below\footnote{More generally, see
\cite{GuMiy 1} for a description of the fundamental groups of acyclic surfaces
with log-Kodaira dimension $1.$}.

Indeed, Exercise \ref{exr7.7} above shows that $\sigma_3\,|\,X_{k,l,s,m}\,:\,
X_{k,l,s,m} \to X_{k,l,s,m-1}$ is the Kaliman transform of the surface
$X_{k,l,s,m-1}$
along the curve $D:= X_{k,l,s,m-1} \cap \{z = 0\}$
with center at the origin.
The repeated application of Lemma \ref{lmm1} and Exercise \ref{exr7.6}
yields the isomorphisms
$$\pi_1(Y_{k,l,s}) = \pi_1(X_{k,l,s,s}) \simeq \pi_1(X_{k,l,s,s-1})
\simeq \ldots \simeq \pi_1(X_{k,l,s,1}) \simeq \pi_1(X_{k,l,s,0} \setminus
\{P_0\})\,$$ where $X_{k,l,s,0} \simeq X_{k,l,s}:= \{x^k - y^l - z^s = 0\}
\subset \C^3.$
Whence, the surface $X_{k,l,s,0} \simeq X_{k,l,s}$ is homotopically equivalent
to the cone over the {\it Pham--Brieskorn 3-manifold} $M_{k,l,s} := X_{k,l,s}
\cap S^5$ (the {\it link} of the surface singularity of $X_{k,l,s}$
in the sphere $S^5$). In turn, $X_{k,l,s} \setminus\{{\overline 0}\}$ is
homotopically equivalent to the link $M_{k,l,s},$ and thus
$\pi_1(Y_{k,l,s}) \simeq \pi_1(M_{k,l,s}).$
We denote the latter group as $G'_{k,l,s}.$

The structure of these groups is well known (see \cite{Mil 3}).
The groups $G'_{k,l,s}$ are finite iff $1/k + 1/l + 1/s > 1,$
infinite nilpotent iff $1/k + 1/l + 1/s = 1.$ If  $1/k + 1/l + 1/s < 1,$ then
$G'_{k,l,s} = [G_{k,l,s},\,G_{k,l,s}],$  where
$$G_{k,l,s}:=<\gamma_1,\,\gamma_2,\,\gamma_3\,|
\,\gamma_1^k=\gamma_2^l=\gamma_3^s=\gamma_1\gamma_2\gamma_3>$$ is a central
extension of the {\it Schwarz triangular group}
$$T_{k,l,s} := <b_1,\,b_2,\,b_3\,|\,b_1^2=b_2^2=b_3^2=1,\,\,(b_1b_2)^k =
(b_2b_3)^l = (b_3b_1)^s =1>\,,$$ which is a discrete group of isometries of
the non-Euclidean plane generated by reflections in the sides of an appropriate
triangle.

Note that for $1/k + 1/l + 1/s < 1$ the triangular group $T_{k,l,s}$
contains a free subgroup with two generators. Therefore,
the group $G'_{k,l,s}$ also contains such a subgroup;
in particular, it is not solvable. Observe that this group is {\it perfect},
i.e. it coincides with its commutator subgroup; indeed, its abelianization
$H_1(Y_{k,l,s};\,\Z)$ is trivial. On the other hand, it is known that
for $(k,\,l) = (k,\,s) = (l,\,s) = 1$ and only under this condition the
Pham--Brieskorn manifold
$M_{k,l,s}$ is a homology 3-sphere; see \cite[Appendix I.8]{Bri, HNK}.

Recall that
$X_{k,l} \setminus L_{k,l} \simeq \C^2 \setminus \Gamma_{k,l}$ where
$\Gamma_{k,l}:=\{x^k - y^l=0\}\subset \C^2$ (see Example \ref{ex6.1}).
The group
$B_{k,l}:=\pi_1(X_{k,l} \setminus L_{k,l}) \simeq \pi_1(\C^2 \setminus
\Gamma_{k,l})$ has the presentation $B_{k,l}= <a,\,b\,|\,a^k = b^l>$
(see e.g. \cite{Di 2}).
In turn, the group $\pi_1(Y_{k,l,s} \setminus \{z=0\})$ is isomorphic to
an index $s$ subgroup of the group $B_{k,l}$ with a cyclic quotient.
By Lemma \ref{fuji}$(a)$, Ker$\,\left(i_{*}\,:\,
\pi_1 (Y_{k,l,s} \setminus \{z=0\}) \to
\pi_1 (Y_{k,l,s})\right) = \,<<\alpha^s>>\,$ where $\alpha\in B_{k,l}$
represents a vanishing loop of the line $L_{k,l} \subset X_{k,l}.$
Let $p,\,q \in \Z$ be such that $kp + lq =1.$ Then one may take
$\alpha = a^qb^p \in B_{k,l}.$

Therefore,
for $1/k + 1/l + 1/s < 1$ and $(k,\,l) = (k,\,s) = (l,\,s) = 1$ the group
$G'_{k,l,s} \simeq \pi_1(Y_{k,l,s})$
is isomorphic to an index $s$ subgroup of the quotient
$$B_{k,l,s}:= B_{k,l}/\,<<\alpha^s>>\, =
 <a,\,b\,|\,a^k = b^l,\,\,(a^qb^p)^s = 1>\,.$$
In particular, for $k = 2,\,l=3,$ and $s \ge 7$ we have that $B_{2,3} = B_3$
is the 3-braid group with generators
$\sigma_1,\,\sigma_2\in B_3$ being vanishing loops of $L_{2,3}$ in $X_{2,3}$
(see Exercise \ref{exr7.5} above), $a =  \sigma_1\sigma_2\sigma_1,\,\,\,b =
\sigma_1\sigma_2,$ and $G'_{2,3,s}$ is isomorphic to an index $s$ subgroup of
the group $$B_{2,3,s} = B_3/<<\sigma_1^s>> =
<\sigma_1,\,\sigma_2\,|\,\sigma_1\sigma_2\sigma_1=\sigma_2\sigma_1\sigma_2,
\,\,\sigma_1^s = \sigma_2^s = 1>\,$$
which consists of the words in the generators $\sigma_1,\,\sigma_2$ with the
sum of exponents divisible by $s.$

\section{Multicyclic $\C^{*}-$coverings}

\subsection{Contractibility of multicyclic $\C^{*}-$coverings. Examples}

To clarify the very idea of the construction of contractible multicyclic
$\C^{*}-$coverings
due to Koras and Russell \cite{KoRu 1, KoRu 2} let us start with simple
examples. We exhibit two different approaches (see
Examples \ref{ex8.1} and \ref{ex8.3} below). Recall first the following
definitions.

\medskip

\defin The {\it suspension} over a topological space $X$ is
the cylinder
$X \times [0,\,1]$ with the bases
$X \times \{0\},\,X \times \{1\}$ being
contracted each one to a point. The {\it join}
$X\star Y$ of two topological spaces $X,\,Y$ is the cylinder
$(X \times Y)\times [0,\,1]$ with the base $(X \times Y) \times \{0\}$
resp. $(X \times Y) \times \{1\}$ being contracted to $X \times \{0\}$
resp. to $Y \times \{1\}.$ Clearly, the join $X\star \Z/2\Z$
is nothing but the suspension over $X.$

\bigskip

We also use the following fact.

\medskip

\nthm{({\bf N\'emethi--Sebastiani--Thom}
\cite{Ne, Di 2}).}
{\it A generic fibre of a polynomial $p(x) + q(y),\,\,p \in
\C[x_1,\dots,x_k],\,\,\,q \in \C[y_1,\dots,y_l],$
is homotopically equivalent
to the join of generic fibres of the polynomials $p$ and $q.$}

\label{NSTthm}
\bigskip

\nex{({\bf The Russell cubic threefold}; see also Examples \ref{exKR},
\ref{ex8.3}).}
The polynomial $q_0 = x(xy+1) \in \C[x,\,y]$ is the hyperbolic modification of
the polynomial
$h = x + x^2 \in \C[x]$. Thus, it is a quasi-invariant of weight $1$ of the
$\C^{*}-$action
$(\lambda,\, (x,y)) \longmapsto (\lambda x,\lambda^{-1}y)$ on $\C^2.$ The
zero fiber $\Gamma_0=q_0^{-1}(0)$ is a disjoint union of two
affine curves isomorphic to $\C$ and to $\C^{*}$, respectively.
Consider the two-fold
$\C^*-$covering $F_0 \to \C^2$ branched along $\Gamma_0$, given as
the surface
$F_0 =\{x+x^2y+z^2=0\}\subset \C^3$ with the projection
$\vi_2\,:\,F_0\longrightarrow \C^2,$ $(x,y,z)\longmapsto (x,y).$
Thus, $\vi_2$ is $\C^{*}-$equivariant with respect to the actions
$(\lambda,\,(x,y,z))\longmapsto (\lambda^2x,\lambda^{-2}y,\lambda z)$
on $F_0$ and $(\lambda,\,(x,y)) \longmapsto
(\lambda^2 x,\lambda^{-2}y)$ on $\C^2.$
The restriction of the above $\C^{*}-$action on $F_0$ to the subgroup
$\omega_2=\{\lambda^2=1\}\simeq \Z/2\Z$ of $\C^{*}$ yields the
monodromy of the covering $F_0 \to \C^2$.
Since this monodromy acts trivially in the homology of $F_0$,
by Corollary \ref{corb2}, $F_0$ is $\Z_3-$acyclic.

By the N\'emethi--Sebastiani--Thom Theorem \ref{NSTthm},
the generic fibre
$F_1 = q^{-1}(1)$ of the polynomial $q=(x+x^2y)+z^2 \in \C[x,y,z]$ is
homotopically equivalent to the join $\Gamma_1 \star \Z/2\Z$ where
$\Gamma_1 := q_0^{-1}(1) \subset \C^2,$
i.e. to the suspension over $\Gamma_1.$ Since the curve $\Gamma_1
\simeq \C^*$ is connected, the fibre $F_1$ is simply connected, and hence,
$G := \pi_1(\C^3 \setminus F_0) \simeq \Z$ (see Exercise \ref{exr7.3}).

The Russell cubic threefold $$X=\{x+x^2y+z^2+t^3=0\}\subset \C^4$$
can be regarded as a three-sheeted cyclic $\C^{*}-$covering over $\C^3$
branched along the surface $F_0$ under the projection
$\vi_3\,:\,(x,y,z,t)\longmapsto (x,y,z)$ onto $\C^3.$ Since the polynomial
$q$ is a quasi-invariant of weight $2$ of the
above $\C^{*}-$action on $\C^3,$ we are under the assumptions of
Kaliman's Theorem \ref{kalim}. Due to this theorem, the Russell cubic
$X$ is a contractible smooth affine variety.

\label{ex8.1}
\bigskip

\exrs (\nexr) Verify that the smooth cubic threefold
$$X'=\{x+x^2y+z^2+t^2=0\}\subset \C^4$$
is simply connected, but not acyclic (what are the homology groups of $X'?$).

\medskip

(\nexr) Show that a generic fiber $p^{-1}$(c)$,\,\,c \in \C^*,$
of the Russell polynomial $p = x+x^2y+z^2+t^3$ is not contractible.

\smallskip

\noindent {\it Hint.} Apply the N\'emethi--Sebastiani--Thom Theorem
\ref{NSTthm}
to find the homotopical type of this fiber.

\bigskip

\nex{(see tom Dieck \cite[Thm. B]{tD 2}).}
More generally, let  $X$ be a smooth contractible affine variety equipped with
a regular $\C^*-$action $t,$ and let $q \in \C[X]$ be a quasi-invariant of $t$
of weight $d \neq 0$ such that $F_0 := q^{*}(0)$ is a smooth reduced (not
necessarily irreducible) principal divisor in $X.$ Fix $s_1,\,s_2 \in \N$ such
that
$d,\,s_1,\,s_2$ are pairwise relatively prime. Consider the smooth affine
hypersurface $$Y_{s_1,\,s_2} := \{q(x) + z^{s_1} + t^{s_2} = 0\} \subset
X \times \C^2\,.$$ We assert that this hypersurface is contractible.

Indeed, consider first the cyclic $\C^*-$covering
$Y_{s_1}\to X,\,\,\,Y_{s_1} = \{q(x)+z^{s_1}= 0\} \subset X \times \C,$
branched to order $s_1$ along $F_0.$
Then $Y_{s_1} \subset X \times \C$ is a smooth reduced irreducible divisor
defined by the quasi-invariant
$q_1(x,\,z) := q(x)+z^{s_1}\in \C[X \times \C]$ of weight
$ds_1$ of the $\C^*-$action
$(\lambda,\,(x,\,z)) \longmapsto (t(\lambda^{s_1},\, x),\,\lambda^d z)$
on $X \times \C.$ Since the monodromy group $G \simeq \Z_{s_1}$ of the
covering acts trivially in the homology of  $Y_{s_1},$ by Corollary
\ref{corb2},
$Y_{s_1}$ is $\Z_p-$acyclic for any prime $p$ which is prime to $s_1,$ and
hence, for any prime divisor $p$ of $s_2.$

Besides, the fibre $F_1 = q^{-1}(1)$ of the regular function $q\in \C[X]$
is connected (see Exercise 7.3). As above, it follows
from the N\'emethi--Sebastiani--Thom Theorem \ref{NSTthm} that the fibre
$q_1^{-1}(1) = \{q(x) + z^{s_1} = 1\}$ of the function
$q_1 \in \C[X \times \C]$ is simply connected. Hence, we have
$\pi_1((X \times \C) \setminus Y_{s_1}) \simeq \Z.$

Therefore, by Kaliman's Theorem 7.1, the total space
of the cyclic $\C^*-$covering $Y_{s_1,\,s_2}\to X \times \C$ branched to order
$s_2$ over $Y_{s_1}$ is a smooth contractible affine variety.

\label{ex8.2}
\bigskip

Applying
Kaliman's Theorem \ref{kalim} successively in the same way as above, one can
derive the following result (cf. Koras and Russell \cite[(7.14)]{KoRu 2}).

\bigskip

\thm {\it Let $X$ be a smooth contractible affine variety
equipped with an effective $\C^{*}-$action. Let $q_i\in \C[X],$
$i=1,\ldots,k,$ be a sequence of quasi-invariants of positive weights
$d_1,\ldots,d_k,$ respectively, and let $s_1,\ldots, s_k$
be a sequence of positive integers. Suppose that the following conditions are
fulfilled:

\smallskip

\noindent $(i)$ For each $i=1,\ldots,k,$
$F_i:=q_i^{*}(0)$ is a smooth reduced irreducible divisor,
the union $\bigcup_1^kF_i$ is a divisor with normal
crossings, and the group $\pi_1(X \setminus \bigcup_1^kF_i)$ is Abelian;

\noindent $(ii)$ $(d_i,s_i)=(s_i,s_j)=1$ for all $i,j=1,\ldots,k;$

\noindent $(iii)$ $F_i$ is $\Z_p-$acyclic for each prime divisor
$p\vert s_i,$ $i=1,\ldots,k.$

\smallskip

\noindent Let $Y\longrightarrow X$ be a multicyclic covering branched
to order $s_i$ along $F_i,$ $i=1,\ldots,k,$ i.e.
$Y=Y_{s_1\ldots s_k}$ is the last one in the tower of cyclic
$\C^{*}-$coverings
$$Y_{s_1\ldots s_k}\longrightarrow
Y_{s_1\ldots s_{k-1}}\longrightarrow \ldots \longrightarrow
Y_{s_1s_2}\longrightarrow Y_{s_1}\longrightarrow X,$$
where $Y_{s_1\ldots s_i}\longrightarrow
Y_{s_1\ldots s_{i-1}}$ is a $\C^{*}-$covering branched
to order $s_i$ over the preimage
of the divisor\footnote{In other words, if $A=\C[X],$ then
$Y=\spec A[\sqrt[{s_1}]{q_1},\ldots,\sqrt[{s_k}]{q_k}].$} $F_i$ in
$Y_{s_1\ldots s_{i-1}}$.

Then $Y$ is a smooth contractible affine variety given in
$X\times \C^k$ by the equations $z_i^{s_i}=q_i(x),$ $i=1,\ldots,k.$}

\label{t81}
\bigskip

\rems \nrem In the case when $X\subset \C^n,$ and $q_j$ is
a variable, $q_j=x_j$ say, the equations
of the cyclic covering $Y_{s_j}\longrightarrow X$ can be obtained
from the equations $P_i(x_1,\ldots,x_n)=0,$ $i=1,\ldots,m,$
which define $X,$ by the substitution $x_j\longmapsto x_j^{s_j},$
i.e.
$Y_{s_j}=\{P_i(x_1,\ldots,x_j^{s_j},\ldots,x_n)=0,\,\,i=1,\ldots,m\} \subset
\C^n.$ In particular,
if $X$ is a hypersurface in $\C^n,$ so is $Y_{s_j}.$

\medskip

\nrem To construct contractible
$\C^{*}-$invariant hypersurfaces
one can use the hyperbolic modification in the same way as in Example
\ref{ex8.1} above.
Recall that if $X=\{h=0\}\subset \C^n$
is a smooth contractible hypersurface,
then the hypersurface $Y_0:=q^{-1}(0)\subset \C^{n+1},$
where $q(\ox,y)=\frac{h(u\ox)}{u},$
is isomorphic to the Kaliman modification of the product $X\times \C$
along the divisor $X$ with center at the origin (Lemma \ref{lmh3}),
and so, by Theorem \ref{thkl}, $Y_0$
is a smooth contractible hypersurface, too.
Moreover, $q$ is a quasi-invariant of weight $1$ of the regular
$\C^{*}-$action
$(\lambda,\, (\ox,u)) \longmapsto (\lambda \ox,\lambda^{-1}u)$ on $\C^{n+1}$
with the only fixed point at the origin (of hyperbolic type with the weights
$(1,\ldots,1,-1)$).

Let smooth hypersurfaces $H_i=\{h_i(\ox)=0\}\bigcap X,$ $i=1,\ldots,k,$ in $X$
satisfy the condition $(i)$ of Theorem \ref{t81}. Put
$q_i(\ox,u)=
\frac{h_i(u\ox)}{u}.$ Then the $q_i$ are $\C^{*}-$invariants
of weight $1$  of the above $\C^{*}-$action on $Y_0.$ The hypersurfaces
$F_i:=q_i^{-1}(0)\bigcap Y,\,\,i=1,\ldots,k,$ in $Y$
also satisfy the condition $(i)$ (note that $F_i$ is the closure of the
$\C^{*}-$orbit of the subvariety $H_i\subset X\simeq Y_0\bigcap H\subset Y_0,
\,\,i=1,\ldots,k,$
where $H := \{u=1\}$; see Exercise \ref{exr6.1}).
This construction can be illustrated by the following simple example.

\bigskip

\ex {\bf The Russell cubic threefold once again} (see Koras-Russell
\cite{KoRu 2}; cf. Examples \ref{exKR}, \ref{ex8.1}).
Starting with $X=\C^2,$ fix two smooth curves
$(f)$ and $(g)$ where $f,\,g \in \C[x,\,z]$, isomorphic to $\C$ and such that
$(f),(g)$ meet transversally
at the origin and in $k$ other points, $k\ge 1.$ For instance,
take $f=z,$ $g=z+x+x^2.$ Then the Kaliman modification $Y$ of
$X\times \C$ along $X$ with center at the origin is nothing but $\C^3.$
The plane curves $(f)$ and $(g)$ give rise, respectively, to the surfaces
$(F)$ and $(G)$ in $\C^3$ where $F=\frac{f(yx,yz)}{y},$ and
$G=\frac{g(yx,yz)}{y}$.
Observe that $(F)$ and $(G)$ are isomorphic to $\C^2,$ meet transversally,
and $\pi_1(\C^3 \setminus ((F) \cup (Q))) \simeq \Z.$

In our particular example $F=z$ and $G=z+x+x^2y.$ The polynomials $F$ and $G$
are $\C^{*}-$quasi-invariants
of weight $1$ with respect to the action
$(\lambda,\,(x,y,z)) \longmapsto (\lambda x,\lambda^{-1}y,\lambda z)$
on $\C^3$. We may also take the plane $H_0=\{y=0\}$ for the
third surface transversal to the first two $(F)$ and $(G)$.

Fix two relatively prime positive integers $s_1,s_2.$ Passing to the
bicyclic $\C^{*}-$covering of $Y\simeq \C^3$ branched to order $s_1$
along $(F)$ and to order $s_2$ along $(G)$ we obtain
a hypersurface $Y_{s_1,s_2}\subset \C^4$ given by the equation
$$x+x^2y+z^{s_1}+t^{s_2}=0$$
which is a smooth contractible threefold. For $s_1=2, s_2=3$ we get
once again the Russell cubic. More generally, passing to the tricyclic covering
 of
$\C^3$ branched
to order $s_0$ resp. $s_1,\,s_2$ along the surface $H_0$ resp. $(F),\,(G)$
where
$(s_i,s_j)=1,\,\,i\neq j$ yields the smooth contractible hypersurface
$\{x+x^2y^{s_0}+z^{s_1}+t^{s_2}=0\}$ in $\C^4.$

\label{ex8.3}
\medskip

\nrem
A theorem due to Koras and Russell \cite[Thm.  4.1]{KoRu 2}
says that any smooth
contractible affine threefold with a `good' hyperbolic  $\C^{*}-$action
appears in the same way as in the above example.

\subsection{The logarithmic Kodaira dimension of multicyclic coverings}

\lemma\footnote{See section 2.2 for the terminology.}{\it Let
$V$ be a smooth projective variety, and let
$L$ be a line bundle on $V$.}

\smallskip

\noindent $(a)$ ({\bf Mori} \cite[Prop. 1.9]{Mo}) {\it The line bundle
$L$ is big (i.e.
$k(V,L)=\dm_{\C}V$) iff for some $k \in \N$ the multiple $kL$
can be written as $kL = A + E$ where $A$ is an ample line bundle on $V,$
and $E$ is an effective one (that is,
$E$ admits a non-zero holomorphic section).}

\smallskip

\noindent $(b)$ ({\bf Kleiman-Kodaira}; see e.g. \cite[(2.3)]{Wil},
\cite[Lemma 0-3-3]{KMM}) {\it If the line bundle $L$ is ample (resp. big), then
for any line bundle $L'$ on $V$ and for any $k \in \N$ large enough the line
bundle $kL - L'$ is ample (resp. big), too.}

\label{lmko}
\medskip

{\bf Proof of $(b)$.} Denote by $NE_1(V)$ the cone of numerically effective
$1-$cycles modulo numerical equivalence on a projective variety $V.$
Recall the {\it Kleiman criterion of
ampleness} \cite{Kl}: {\it a line bundle $L$ on $V$ is ample iff it is positive
on the cone $NE_1(V)$ with the origin being deleted}.

This finite dimensional cone is closed, and
hence, it has a compact intersection with the unit sphere.
Thus, the openness of
ampleness follows.

Let $L$ be a big line bundle, and let $k_0L = A + E$
be a decomposition as in $(a)$.
Then for $n_0 \in \N$ large enough we have
$n_0k_0L - L' = (n_0A - L') + n_0E$ where $n_0A - L'$ is ample. Therefore,
by $(a)$, the line bundle $n_0k_0L - L' $ is big.
It follows that for any $k \ge n_0k_0$ the line bundle
$kL - L' = (n_0k_0L - L') + (k - n_0k_0)L$ is also big.
\qed

\bigskip

\nprop{({\bf Kaliman} \cite[Lemma 11]{Ka 1}).}
{\it Let $X$ be a quasi-projective
variety,
$(V,D)$ be an SNC-completion of $X,$ and $Z=\sum_1^kZ_i$
be an SNC-divisor on $V$ such that $D\bigcup Z$ is also an SNC-divisor,
and $D$ and $Z$ have no irreducible component in common.
Let $Y = Y_{\overline s}\longrightarrow X,\,\,{\overline s}:=
(s_1,\dots,s_k),$ be a ramified covering branched to order $s_i$
over $Z_i\bigcap X,\,\,i=1,\dots,k.$ Then
$$\ok(Y)=k\Bigl(V,\,K_V+D+\sum_{i=1}^k(1-\frac{1}{s_i})Z_i\Bigr).$$}
\label{p81}
\proof
One can compactify $Y$ by an SNC-divisor
$D^\prime$ (i.e. $Y=V^\prime\setminus D^\prime$) to obtain a commutative
diagram of morphisms

\begin{picture}(200,70)
\unitlength0.2em
\put(75,25){$Y$}
\put(105,25){$V^\prime$}
\put(90,25){$\hookrightarrow$}
\put(75,5){$X$}
\put(90,5){$\hookrightarrow$}
\put(105,5){$V\,\,.$}
\put(77,22){$\vector(0,-1){11}$}
\put(107,22){$\vector(0,-1){11}$}
\put(110,16){{\rm $\varphi$}}
\end{picture}

\noindent Then $\varphi^{*}(Z_i)=s_iZ_i^\prime+E_i$ where
each divisor $E_i$ is $\varphi-$exceptional, i.e. $\cdm_V\varphi(E_i)\ge 2.$
The restriction
$\varphi\,|\,(V' \setminus (D' \cup Z'))\,:\,V' \setminus (D' \cup Z') \to
V \setminus (D \cup Z)$ is an {\'e}tale covering. Therefore,
by the Logarithmic Ramification Formula (R),
we have
$$K_{V^\prime}+D^\prime+Z^\prime=\varphi^{*}(K_V+D+Z)+R\,$$
where $R$ is an effective $\varphi-$exceptional divisor in $V'$
(see Theorem \ref{prr}.$(e),\,(f)$).
Hence,
$$K_{V^\prime}+D^\prime=\varphi^{*}(K_V+D+Z)+R-Z^\prime
=\varphi^{*}(K_V+D+Z)+$$
$$+\varphi^{*}\Bigl(\sum_{i=1}^k(-\frac{1}{s_i})Z_i\Bigr)
+R+\sum_{i=1}^k\frac{1}{s_i}E_i
= \varphi^{*}\Bigl(K_V+D+\sum_{i=1}^k(1-\frac{1}{s_i})Z_i\Bigr)+E\,$$
where $E:=\sum_{i=1}^k\frac{1}{s_i}E_i+R$
is a $\varphi-$exceptional $\Q-$divisor.
By \cite[Lemma 1]{Ii 1} or \cite[Thm. 10.5]{Ii 3},
$\ok(V',\,\varphi^{*}(D_1) +
E)=\ok(V,\,D_1)$ for any $\Q-$Cartier divisor $D_1$ on $V$ where $E$ is a
$\varphi-$exceptional divisor in $V'$
(indeed, a meromorphic section  of the associated line bundle
$[\varphi^{*}(D_1)]$ with poles at most along $E$ has no pole).  Thus, the
assertion follows.
\qed

\bigskip

\ncor{(\cite{Ka 1}, \cite[Cor.6.2]{KoRu 2}).} {\it If
${\overline s'} = (s'_1,\dots,s'_k)$ and $s'_i \ge s_i,\,i=1,\dots,k,$ then
$\ok(Y_{\overline s'}) \ge \ok(Y_{\overline s}).$}

\label{c81}
\bigskip

\cor {\it If $X\setminus Z$ is a variety of
log-general type, that is, $\ok(X\setminus  Z)=\dm_{\C}X,$ then for
$s_i,\,\,i=1,\dots,k,$ large enough $Y_{\overline s}$ is a variety of
log-general type, too.}

\label{c82}
\medskip

\proof
Indeed, by Lemma \ref{lmko}$(b)$, for
$s_i>>1,\,\,i=1,\dots,k,$ we have
$$\ok(Y_{\overline s})=k\Bigl(V,\,K_V+(D+Z)-\sum_{i=1}^k\frac{1}{s_i}Z_i\Bigr)
=k(V,\,K_V+D+Z)=\ok(X\setminus Z) = {\rm dim}_{\C}\,X\,.$$
\qed

\bigskip

\nprop{(see \cite {Ka 1}, \cite[Prop. 6.5]{KoRu 2}).} {\it
Consider the Koras-Russell threefolds
$Y_{\overline s}\subset \C^4$ where ${\overline s} = (s_1,\, s_2,\, s_3)$
and $(s_i, \,s_j) = 1, \,\,i \neq j$ given as\footnote{This is a particular
kind of the Koras-Russell
threefolds; see Example \ref{ex8.3} above.}
$$Y_{\overline s} = \{x+x^2y^{s_1}+z^{s_2}+t^{s_3}=0\}\,.$$
If $s_1,\,s_2,\,s_3 >> 1,$ then $Y_{\overline s}$ is an exotic $\C^3$, and
$\ok(Y_{\overline s})=2$.}

\label{p82}
\medskip

{\bf Proof.} Set
$$X = \{x + x^2u_1 + u_2 + u_3 = 0\} \subset \C^4\,\,\,\,\,\,{\rm
and}\,\,\,\,\,\,Z_i = \{u_i = 0\} \subset X,\,\,i=1,2,3.$$
Evidently, $X \simeq \C^3,\,\,Z_i \simeq \C^2,\,i=1,2,3,$ and $Z:= Z_1\cup Z_2
\cup Z_3$ is an SNC-divisor in $X$.
The threefold $Y_{\overline s}$ is a tricyclic covering of $X$
branched to order $s_i$ along $Z_i,\,i=1,2,3,$ with the covering morphism
$\varphi_{\overline s}\,:\,(x,\,y,\,z,\,t) \longmapsto (x,\,u_1,\,u_2,\,u_3) :=
(x,\,y^{s_1},\,z^{s_2},\,t^{s_3})$. By Theorem \ref{t81}, it follows that
$Y_{\overline s} \subset \C^4$ is a smooth contractible affine hypersurface.
Due to the Dimca-Ramanujam Theorem \ref{thdr}, the variety $Y_{\overline s}$ is
diffeomorphic to $\R^6$. It remains to show that $\ok(Y_{\overline s})=2$
when $s_1,\,s_2,\,s_3$ are large enough.

Due to Corollary \ref{c81}, $\ok(Y_{\overline s})\ge 2$ for sufficiently large
$s_1,\,s_2,\,s_3$ if it is so
for a particular choice of ${\overline s} = (s_1,\, s_2,\, s_3)$
(even without the assumption of relative primeness, which guarantees
the contractibility).

Note that the hypersurface $Y_{\overline s}\subset \C^4$ is invariant under the
hyperbolic
linear $\C^*-$action on $\C^4$
$$G\,:\,(\lambda,\,(x,\,y,\,z,\,t)) \longmapsto
(\lambda^a x,\, \lambda^{-b} y,\,\lambda^c z,\,\lambda^d t)\,$$ where
$$a =  s_1s_2s_3, \,b = s_2s_3, \,c = s_1s_3, \,d =s_1s_2\,.$$
The morphism $\varphi_{\overline s}\,:\,Y_{\overline s} \to X$ is
a $\C^*-$covering with respect to the $\C^*-$action $G$ on $Y_{\overline s}$
and the $\C^*-$action
$${\overline G}\,:\,(\lambda,\,(x,\,u_1,\,u_2,\,u_3)) \longmapsto
(\lambda^a x,\, \lambda^{-a} u_1,\,\lambda^a u_2,\,\lambda^a u_3)$$ on $X$.
We have: spec$(\C[X])^{{\overline G}}
:= X//{\overline G} \simeq S$ where $S := \{u + u^2 + v + w = 0 \} \subset
\C^3;$ clearly, $S \simeq \C^2.$
Indeed, $(\C[X])^{{\overline G}} = \C[u,\,v,\,w]$ where
$u:=u_1 x,\,v:=u_1u_2,\,w:=u_1u_3 \in (\C[X])^{{\overline G}}$ are the basic
${\overline G}-$invariants.
This yields the following commutative diagram of
morphisms:

\begin{picture}(200,95)
\unitlength0.2em
\put(75,25){$Y_{\overline s}$}
\put(105,25){$Y_{\overline s}//G = S_{\overline s}$}
\put(84,27){$\vector(1,0){15}$}
\put(92,30){{\rm $\rho_{\overline s}$}}
\put(75,5){$X$}
\put(84,6){$\vector(1,0){15}$}
\put(92,9){{\rm $\rho$}}
\put(105,5){$X//{\overline G} = S\,\,$}
\put(77,22){$\vector(0,-1){11}$}
\put(70,16){{\rm $\varphi_{\overline s}$}}
\put(114,22){$\vector(0,-1){11}$}
\put(118,16){{\rm ${\overline \varphi}_{\overline s}$}}
\end{picture}

\noindent where $S_{\overline s} := Y_{\overline s}//G=$spec$\,(\C[Y_{\overline
s}])^G$ is a normal surface. A generic fiber of the quotient morphism
$\rho_{\overline s}\,:\,Y_{\overline s} \to S_{\overline s}$ (i.e. a generic
orbit) is isomorphic to $\C^*$. Since $\ok(\C^*) = 0$, from the Addition
Theorems \ref{prr}$(c),(d)$\footnote{They are still available, although the
quotient surface $S_{\overline s}$ might be singular.} we obtain
$$2 = \dm S_{\overline s} \,\ge \,\ok(Y_{\overline s})\, \ge \,\ok(S_{\overline
s})\,.$$ Thus, it remains to find a particular triple ${\overline s} = (s_1,\,
s_2,\, s_3)$ such that $\ok(S_{\overline s}) = 2$.

Note that the threefold $Y_{\overline s}$ is the closure of the $G-$orbit of
the surface $T_{\overline s}:= Y_{\overline s} \cap H$ where
$H := \{y = 1\} \subset \C^4$ (see Exercise \ref{exr6.1}). The surface
$T_{\overline s}$ is invariant under the induced action of the cyclic subgroup
$\omega_b \subset \C^*$ on $Y_{\overline s}$, and
$S_{\overline s} = Y_{\overline s}//G \simeq T_{\overline s}/\omega_b$.

Take $s_1=pq,\,s_2=p,\,s_3=q$ where $p,\,q \in \N$ are prime
and distinct.
Then we have $Y_{\overline s} = \{x + x^2y^{pq} + z^p + t^q =0\}$, and
$$G(\lambda,\,(x,\,y,\,z,\,t)) =
(\lambda^{p^2q^2} x,\, \lambda^{-pq} y,\,\lambda^{pq^2} z,\,\lambda^{p^2q}
t)\,.$$ Therefore, the subgroup
$\omega_b = \omega_{pq} \subset \C^*$ which coincides with the
non-effectiveness kernel
of the $\C^*-$action $G$ on $Y_{\overline s}$ acts trivially
on
$T_{\overline s} = Y_{\overline s} \cap H
= \{x + x^2 + z^p + t^q =0\}\subset \C^3$.
Hence,
$S_{\overline s} = Y_{\overline s}//G \simeq T_{\overline s} \subset
Y_{\overline s}$. The projection $$\rho\circ\varphi_{\overline
s}\,|\,T_{\overline s}\,:\,T_{\overline s} \to S,\,\,\,\,\,\,(x,\,z,\,t)
\longmapsto (u,\,v,\,w) =
(x,\,z^p,\,t^q)\,,$$ is a bicyclic covering branched to order $p$ resp. $q$
over the curve $C_1:=\{v=0\} \subset S$ resp. $C_2:=\{w=0\} \subset S$.
By Corollary \ref{c82} above, we have
$\ok(S_{\overline s}) = \ok(T_{\overline s})
= \ok(S \setminus (C_1 \cup C_2)),$
if $p$ and $q$ are sufficiently large. Thus, the proof is completed by the
following simple exercises.
$\Box$

\bigskip

\exrs Show that

(\nexr) $(S,\,C_1 \cup C_2) \simeq (\C^2, \, D_1 \cup D_2),$
where $D_1 := \{y=1\},\,D_2 := \{y = x^2\} \subset \C^2;$
and that\footnote{Cf. \cite[Lemma 16]{Ka 1}, \cite[Lemma 6.3]{KoRu 2}.}

(\nexr) $\ok(\C^2 \setminus (D_1 \cup D_2))=2.$

\bigskip

\rem (see \cite[Prop.  7.8.]{KoRu 2}) However, for $s_1=1$ the
threefold $Y_{\overline s} \subset \C^4$ is dominated by
$\C^3;$ in particular, it has the log-Kodaira dimension $\ok=-\infty$.
Indeed, if $s_1=1,$ then for any $x\ne 0,$ $y$ is expressed
in terms of $z$ and $t,$ whence the part $\{x\ne0\}$ of the threefold
$Y_{\overline s}$ is isomorphic to the cylinder $\C^2\times \C^{*}.$
The `book-surface' $B:=\{x=0\} \subset Y$ is the product
$\C\times \Gamma_{s_2, s_3}$ where
$\Gamma_{s_2, s_3} :=\{z^{s_2}+t^{s_3}=0\}\subset \C^2.$ Fix a smooth point
$\rho\in \Gamma_{s_2, s_3},$
and perform the Kaliman modification
$\sigma:Y_{\overline s}^\prime\longrightarrow Y_{\overline s}$ of
the 3-fold $Y_{\overline s}$ along the divisor $B$ with the center
$C:=\C\times\{\rho\}.$
In this way, we replace the singular book-surface $B$ by a smooth surface
$E'\simeq \C^2,$
and replace the function $x$ by a function
$f \,:\,Y_{\overline s}^\prime\longrightarrow \C$ such that all the fibers of
$f$
are smooth reduced surfaces isomorphic to $\C^2.$ By the Miyanishi Theorem
\ref{thmMI} below we have: $Y_{\overline s}^\prime\simeq \C^3.$ So,
$\sigma:\C^3\simeq Y_{\overline s}^\prime \longrightarrow Y_{\overline s}$
is a birational (whence, dominant) morphism.

In the case of Russell's cubic $X_0=\{x+x^2y+z^2+t^3=0\}\subset \C^4,$
a dominant morphism $\C^3\longrightarrow X_0$ can be given explicitly as
$(u,\,v,\,w) \longmapsto (x,\,y,\,z,\,t)$ where
$$(x,\,y,\,z,\,t) = \left( -u,\,\,{u-(u^2v + 1)^2 -
(u^2w + u/3 - 1)^3 \over u^2},\,\,u^2v + 1,\,\,u^2w + u/3 - 1\right)\,.$$
We conclude this section by the following characterization of the affine
3-space $\C^3$ due to Miyanishi \cite{Miy 1}.

\medskip

\nthm{({\bf Miyanishi}).} {\it Let $X$ be a smooth affine threefold.
Then $X$ is isomorphic to the affine 3-space $\C^3$ iff the following
conditions hold:

\smallskip

\noindent $(i)$ the Euler characteristic $e(X)$ of $X$ is equal to $1;$

\smallskip

\noindent $(ii)$ the algebra $\C[X]$ of regular functions on $X$ is UFD, and
all its invertible elements are constants;

\smallskip

\noindent $(iii)$ there exists a non-empty Zariski open subset
$\Omega \subset X$ isomorphic to a cylinder $\Gamma \times \C^2$ where
$\Gamma$ is an affine curve;

\smallskip

\noindent $(iv)$ the algebra of regular functions on each irreducible component
of the divisor $X \setminus \Omega$ is UFD.}

\label{thmMI}
\medskip

\rem Observe that for the Russell cubic threefold the
conditions $(i)-(iii)$ are fulfilled, and only the last
condition $(iv)$ does not hold.

\section{The Makar-Limanov invariant of the Russell cubic threefold}

Let $X$ be an affine variety. We assume in what follows that $X$
is irreducible, so that the algebra
$A=\C[X]$ of regular functions on $X$ is an integral domain. Makar-Limanov
\cite{ML 2}
(see also \cite{KaML 3})
introduced a subring $\ML(A)$ of a ring $A$ such that $\ML(A)$ is
invariant under ring isomorphisms; that is,
if $B\simeq A,$ then $\ML(B)\simeq \ML(A).$ He proved the following theorem.

\bigskip

\nthm{({\bf Makar-Limanov} \cite{ML 2, ML 3}).} {\it Set $A_0=\C[X_0]$ where
$X_0=\{x+x^2y+z^2+t^3=0\}\subset \C^4$ is the Russell cubic threefold. Then
$\ML(A_0)$ is not isomorphic to $\C=\ML(\C[x,y,z]).$
Thus, $X_0$ is not isomorphic to $\C^3,$ and hence, $X_0$ is an exotic $\C^3.$}

\label{thmLML}
\medskip

Later on, Kaliman and Makar-Limanov \cite{KaML 3}
extended this result to
all the Koras-Russell threefolds. This was one of the crucial steps
in the recent proof of the Linearization Conjecture for $n=3$
(Koras and Russell \cite{KoRu 2, KoRu 3}, Kaliman, Koras, Makar-Limanov,
Russell \cite{KaKoMLRu}).

\bigskip

\nthm{\cite{KoRu 2, KoRu 3, KaKoMLRu, KrPo, Po}.} {\it Any regular
$\C^{*}-$action on $\C^3$
is linearizable (i.e. it is conjugate with a linear $\C^{*}-$action on $\C^3$).
Moreover, any regular action of a connected reductive group on $\C^3$
is linearizable.}

\bigskip

Here we give an exposition of Makar-Limanov's result following
a simplified approach due to Derksen \cite{De}.

\subsection{$\C_{+}-$actions and locally nilpotent derivations}

{\bf Notation and definition.} As before, $\C_{+}$ stands for the additive
group
of the complex number field. Further, $\C^{[n]}=
\C[x_1,\ldots,x_n]$ denotes the polynomial algebra in $n$ variables.
Let $A$ be an algebra over $\C.$ A {\it locally nilpotent
derivation} of $A$ is a $\C-$linear homomorphism
$A\stackrel{\partial}\longrightarrow A$
satisfying the Leibniz rule and such that
for any $a\in A,$ $\partial^n a=0$ for some $n=n(\partial,\,a)\in \N.$
Denote by $\LND(A)$ the set of all locally nilpotent
derivations of the algebra $A.$

\smallskip

Let, as above, $X$ be an affine variety and $A = \C[X]$ be the algebra of
regular
functions on $X.$
Any regular $\C_{+}-$action
$\lambda:\C_{+}\times X\longrightarrow X$ induces
an algebra homomorphism $$A\stackrel{\varphi}\longrightarrow A[t]
\simeq \C[\C_{+}\times X], \quad p \in A\stackrel{\varphi}\longmapsto
p(\lambda(t,x))\in A[t].$$
Set $\partial p=\frac{d}{dt} \vert_{t=0}\, (p\circ \lambda).$
Then $\partial \in \LND(A).$ Vice versa, any locally nilpotent
derivation $\partial \in\LND(A)$ corresponds to the algebra homomorphism
$\varphi_\partial\,:\,A\longrightarrow A[t]$ given by
$$\varphi_\partial(a)=\exp (t\partial)(a)=\sum_{i=0}^\infty \frac{t^i\partial
^ia}{i!}, \quad a\in A,$$ and thus,
to a regular $\C_{+}-$action on the variety $X.$

\bigskip

\exrs
(\nexr) \label{exr9.1}
Prove the equality $A^{\varphi_\partial}=A^\partial$ where
$A^\partial:=\Ker \partial,$ and $A^{\varphi_\partial}$ is the subalgebra
of invariants of the $\C_{+}-$action $\varphi_\partial$ on $X.$
Deduce that the action $\varphi_\partial$ is trivial iff $\partial=0.$
Verify that
the subalgebra $A^\partial$ is algebraically closed in the algebra $A.$

\smallskip

(\nexr) \label{exr9.2}
Let $\partial \in \LND(A)\setminus \{0\}.$ Verify that
the transcendence degree of the algebra extension $[A\,:\,A^\partial]$
is 1. More precisely, let $r_0$ be any element of $A$
such that $\partial r_0\in A^\partial$ and $r_0\notin A^\partial.$
Show that the subalgebra $A^\partial [r_0]\subset A$ is a free $A^\partial-$
module, and for any element $a\in A$ there exists an  element $b\in
A^\partial\setminus\{0\}$
such that $ba\in A^\partial [r_0].$

\smallskip

(\nexr) \label{exr9.3}
Given a linear representation $\varphi\,:\, \C_{+}\longrightarrow
\GL_n(\C),\,\,\,
t \stackrel{\varphi}\longmapsto e^{tB}$ where $B \in L_n(\C),$ verify that it
provides a regular
$\C_{+}-$action on $\C^n$ iff it is unipotent, i.e. iff
$B$ is a nilpotent matrix. Or, equivalently, iff the associated derivation
$\partial_{\varphi} (p)=<Bx, \grad p>$ of the polynomial algebra
$\C^{[n]}=\C[x_1,\ldots,x_n]$ is locally nilpotent.

\smallskip

(\nexr) \label{exr9.4} Given a locally nilpotent derivation $\partial \in
\LND(A)$ show that $a\partial \in \LND(A)$ for any {\it
$\partial-$ constant} $a \in A^\partial.$ Thus, $A^\partial \partial \subset
\LND(A).$ Conclude that if $A$ is an algebra over $\C$ of $\trdg A \ge 2$ such
that $\LND(A) \neq \{0\},$ then the automorphism group $\Aut A$ is of infinite
dimension. In particular, if an affine variety $X$ with $\dim X \ge 2$
admits a non-trivial regular $\C_+-$action, then\footnote{This observation was
communicated to us by L. Makar-Limanov.} $\dim Aut X = \infty.$

\smallskip

(\nexr) \label{exr9.5}
Let $\Gamma$ be an irreducible affine algebraic curve. Show that it
admits a non-trivial regular $\C_{+}-$action iff $\Gamma \simeq \C.$

\bigskip

\defin
Let $A$ be an algebra over $\C$. The {\it Makar-Limanov invariant}
$\ML(A)$ of the algebra $A$ is the subalgebra
$\ML(A):=\bigcap_{\partial \in \LND(A)} A^\partial \subset A.$

\smallskip

The {\it Derksen invariant} $\Dk(A)$ of the algebra $A$
is the smallest subalgebra of $A$
which contains  $A^\partial$ for all $\partial\in LND(A)\setminus \{0\}.$

\label{defIN}
\bigskip

Clearly, $\ML(\C^{[n]})=\C,$ and $\Dk(\C^{[n]})=\C^{[n]}.$

\bigskip

\nthm{({\bf Derksen} \cite{De}).} {\it Let, as above,
$A_0$ denotes the algebra of regular functions on the Russell cubic threefold
$X_0=\{x + x^2y + z^2 + t^3 = 0\} \subset \C^4.$
Then $\Dk(A_0)\ne A_0.$ Hence, the algebra $A_0$ is not
isomorphic to $\C^{[3]}.$ In turn, the Russell cubic
$X_0$ is not isomorphic to
$\C^3,$ i.e. $X_0$ is an exotic $\C^3.$}

\label{derks}
\medskip

Before proceeding with the proof, we recall the following notions (see e.g.
\cite{Bou}).

\subsection{Degree functions, filtrations and the associated
graded algebras}

Let $A$ be an integral domain (usually, it will be also an algebra over $\C$).

\medskip

\defin A {\it degree function}\footnote{In a similar way, one may define a
degree function with values in
arbitrary ordered semigroup.} $\dg:A\longrightarrow \Z\bigcup \{-\infty\}$
on $A$ is a mapping which satisfies the following
axioms:

\medskip

(d1) $\dg 0=-\infty,$ and $\dg a\in \Z$ for all $a\ne 0;$ $\dg 1=0.$

(d2) $\dg fg=\dg f+\dg g$ for all $f,g\in A.$

(d3) $\dg (f+g)\le \max\{\dg f,\dg g\}$ for all $f,g\in A.$

\bigskip

\defin A degree function $\deg$ determines an {\it ascending filtration}
$F=\{\fr{i}\}$ on
$A$ where $\fr{i}:=\{a\in A\;\vert\; \dg a\le i\}.$ This filtration satisfies
the following conditions:

\medskip

(f1) $\fr{i}$ is a $\C-$linear subspace of $A,$ and $\fr{i}\subset \fr{i+1}$
({\it ascending}).

(f2) $A=\bigcup_{i\in \Z}\fr{i}$ ({\it exhaustive});
$\bigcap_{i\in \Z}\fr{i}=\{0\}$
({\it separated}); $1\in \fr{0} \setminus \fr{-1}.$

(f3) $(\fr{i}\setminus \fr{i-1})(\fr{j}\setminus \fr{j-1})\subset
(\fr{i+j}\setminus \fr{i+j-1}).$

\smallskip

Clearly, $F^0A \subset A$ is a subring (resp. subalgebra),
and the ring $A$ represents as an $F^0A-$module.

\smallskip

Vice versa, given a filtered domain $(A,\,F)$ which satisfies
the conditions (f1)-(f3) one can define a degree function $d_F$ on $A$ as
follows: $d_F(0) = -\infty$ and $d_F$(a)$ = i$ iff $a \in \fr{i}\setminus
\fr{i-1}.$

\bigskip

\defin
The {\it associated graded algebra} $\Gr A=\oplus_{i\in \Z} \Gr^i A$
of a filtered algebra $(A,\,F)$ where $\Gr^i A:=\fr{i}/\fr{i-1}$ can be
identified with the algebra of the Laurent
polynomials
$\{\sum_{i=k}^{k+l}\baf_iu^i\}$ where $\baf_i$ is either zero or
is equal to $\gr f_i:=f_i+\fr{i-1}\in \Gr^iA$ for some
$f_i\in \fr{i}\setminus \fr{i-1}.$
Due to the property (f3) of filtrations,
the mapping $\gr:A\longrightarrow \Gr A,$ $\gr f=\baf,$
is a homomorphism of multiplicative semigroups.

\bigskip

\defin \label{defWDF}
A {\it weight degree function} on the polynomial algebra $\C^{[n]}$ is a degree
function $d$ such that $d(p) = \max_i \{d(m_i)\}$ where $p \in \C^{[n]}$ is
a non-zero polynomial, and $m_i$ runs over the set
$M(p)$ of all the monomials of $p.$  Clearly,
$d$ is uniquely determined by the {\it weights}
$d_i := d(x_i),\,\,i=1,\dots,n.$ A weight degree function $d$ defines a grading
$\C^{[n]} = \oplus_{j \in \Z} \C^{[n]}_{d,\,j}$ where
$\C^{[n]}_{d,\,j} \setminus \{0\}$ consists of all the $d-$quasihomogeneous
polynomials of $d-$degree $j$. Accordingly, for any $p \in \C^{[n]} \setminus
\{0\}$ we have a unique decomposition
$p = \sum\limits_{i=m(p)}^{d(p)} p_i$ into a sum of $d-$quasihomogeneous
components; here $p_d:=p_{d(p)}$ is called the {\it principal
$d-$quasihomogeneous component} of $p.$ It is clear that $(pq)_d = p_dq_d.$

Let $X = (I) \subset \C^n$ be a reduced irreducible affine variety defined
by a prime ideal $I \subset \C^{[n]}.$ Denote $A = \C[X] = \C^{[n]}/I,$ and
let $\hI$ be the (graded) ideal in $\C^{[n]}$ generated by the principal
$d-$quasihomogeneous components $p_d$ where $p$ runs over $I.$
We say that the weight degree function $d$ is {\it appropriate for the ideal}
$I$ if
the following conditions hold:

\smallskip

\noindent $(*)$ ${\overline 0} \in X,$ i.e. $I \subset \alpha
:=(x_1,\dots,x_n);$

\smallskip

\noindent $(**)$ the ideal  $\hI$ is also prime, and
$\,\,x_i \notin \hI \,\,\,\forall\,\, i=1,\dots,n.$

\smallskip

\noindent For $f \in A \setminus \{0\}$ set $$d_A(f) = \min_{p \in [f]}
\{d(p)\}\,\,\,\,\,\,\,\,\,\,
{\rm where}\,\,\,\,\,\,\,\,\,\, [f]:= \{p \in \C^{[n]}\,|\,\,\,p\,|\,X =
f\}\,.$$
\exrs (\nexr) \label{exr9.6} Show that $d_A(f) = d(p)$ for a polynomial
$p \in [f]$ iff
$p_d \notin {\hI}.$

\medskip

(\nexr) \label{exr9.7} Assume that a weight degree function $d$ on the
polynomial algebra $\C^{[n]}$ is appropriate for an ideal $I\subset \C^{[n]}.$
Deduce that
$d_A$ is a degree function on $A,$ and that $d_A (\tx_i) = d(x_i) = d_i$
where $\tx_i:= x_i \,|\,X = x_i + I \in A,\,\,i=1,\dots,n.$
Hence, due to the property (d2) of a degree function,  $d_A(m\,|\,X) = d(m)$
for any monomial $m \in \C^{[n]}.$

\smallskip

{\it Hint.} Suppose that $f \in A$ and $d_A(f) = -\infty,$
that is, there exists a sequence of polynomials
$p_j \in \C^{[n]},\,\,j=1,\dots,$
such that $p_j\,\left.\right | \,X = f$ and
$\lim\limits_{j\to \infty} d(p_j) = -\infty.$
For $p \in \C^{[n]}$ set $\mu(p) = \min_{m \in M(p)}\{\deg m\}$
where deg is the usual degree. Then $p \in \alpha^{\mu(p)}$ where, as above,
$\alpha \subset \C^{[n]}$ denotes the maximal ideal
which corresponds to the origin of $\C^n.$  By the condition $(*)$ from
Definition \ref{defWDF},
$\widetilde{\alpha}:=  (\tx_1,\dots,\tx_n) \subset A$ is a proper ideal,
and we have $f = p_j\, | \,X \in \widetilde{\alpha}^{\mu(p_j)},\,\,j=1,\dots.$
Thus, by the Krull Theorem, $f \in \bigcap_{n \in \N} \widetilde{\alpha}^n =
\{0\},$
and so, $f = 0.$
Hence, $d_A(f) > -\infty$ for any
$f \in A \setminus \{0\}.$

The rest of the exercise, including checking of the other properties of
a degree function, can be done without difficulty.

\medskip

(\nexr) \label{exr9.8}
let $F = \{\fr{i}\}$ be the filtration on
$A$ determined by the above degree function $d_A,$ and let $\hA = \Gr A$ be the
associated graded algebra. Verify that the elements
$\hx_1,\dots,\hx_n\in \hA$ where $\hx_i := \gr \tx_i\in \hA,$
generate the graded algebra $\hA.$

\bigskip

\nlemma{({\bf Kaliman, Makar-Limanov}\footnote{We place
here this lemma and the preceding definition and exercises with a kind
permission
of Sh. Kaliman and L. Makar-Limanov.}\cite[Prop. 4.1]{KaML 4}).}
{\it Keeping the same
notation and assumptions as in the above exercises we have
$${\hA} \simeq \C^{[n]}/{\hI} = \C[{\hX}]\,$$
where ${\hX} = ({\hI}) \subset \C^n$ is the affine variety defined
by the prime ideal ${\hI}.$}

\label{lm9.1}
\medskip

\proof
According to Exercise \ref{exr9.8}, the elements $\hx_1,\dots,\hx_n\in \hA$
generate the graded algebra $\hA.$ Henceforth,
$ \hA =  \C^{[n]}/ J$ where $J \subset \C^{[n]} = \C[\hx_1,\dots,\hx_n]$ is
the ideal of relations between the generators $\hx_1,\dots,\hx_n$ in $ \hA .$
Thus, we must show that $J = {\hI}.$

Fix an arbitrary polynomial $p = \sum\limits_{i=m(p)}^{d(p)} p_i \in I.$ Then
$p \equiv 0 \quad \mod I,$ i.e. $p_d  \equiv -\sum\limits_{i=m(p)}^{d(p)-1}
p_i\quad \mod I,$ and hence
$$d_A(p_d\,|\,X) \le \max_{m(p) \le i \le d(p)-1} \{d_A(p_i\,|\,X)\} \le
\max_{m(p) \le i \le d(p)-1} \{ d(p_i)\} < d(p) = d(p_d)\,.$$
Therefore, $p_d\,|\,X \in \fr{d(p)-1}.$

Since the weight degree function $d$ is appropriate for the ideal
$I,$ by Exercise
\ref{exr9.7}, we have $d_A(m_j\,|\,X) = d(m_j) = d(p)$ for any monomial
$m_j \in M(p_d).$ Thus, $m_j\,|\,X \in \fr{d(p)} \setminus
\fr{d(p)-1}$ for any $m_j \in M(p_d),$
and $p_d\,|\,X = \sum_{m_j \in M(p_d)} (m_j\,|\,X) \in \fr{d(p)-1}.$
It follows that $ (m_j\,|\,X)^{\wedge} := \gr (m_j\,|\,X) =
m_j(\hx_1,\dots,\hx_n) \in
\hA^{d(p)},$ and $\sum_{m_j \in M(p_d)} (m_j\,|\,X)^{\wedge} = 0$ in
$\hA^{d(p)},$ i.e. $p_d(\hx_1,\dots,\hx_n) = 0$ in $\hA^{d(p)}.$
Whence, $p_d \in J,$ and so, $\hI \subset J.$

Vice versa, fix an element $f = \sum\limits_{i=m(f)}^{d(f)} f_i\in J.$  It is
clear that
$f_i(\hx_1,\dots,\hx_n) \in \hA^i$ (indeed, as above, this is true for any
monomial
$m \in M(f_i)$). Since $$\sum\limits_{i=m(f)}^{d(f)} f_i(\hx_1,\dots,\hx_n)
= f(\hx_1,\dots,\hx_n) = 0$$ we have $f_i(\hx_1,\dots,\hx_n) = 0$ for each
$i = m(f),\dots,d(f).$ Thus, $J$ is a homogeneous ideal of the $d-$graded
algebra $\C^{[n]}$ (see Definition \ref{defWDF} above). Hence, it is enough to
show that
$J_r \subset \hI$ for any $d-$homogeneous component $J_r$ of $J.$

Let $f \in J_r$ be a $d-$quasihomogeneous polynomial of $d-$degree $r = d(f).$
For any monomial $m \in M(f)$ we have, as above, that
$m\,|\,X \in \fr{r} \setminus \fr{r-1},$
and so, $m(\hx_1,\dots,\hx_n) \in \hA^r.$ Since $\sum_{m \in M(f)}
m(\hx_1,\dots,\hx_n) = f(\hx_1,\dots,\hx_n) = 0,$ it follows that
$f\,|\,X \in  \fr{r-1},$ i.e. $d_A(f\,|\,X) < r =d(f).$
By Exercise \ref{exr9.6},
this implies that $f_d \in \hI.$ But $f = f_d,$ and so, we are done.
\qed

\bigskip

{\bf Gradings and $\C^{*}-$actions} (see e.g. \cite{KamRu, Ru 3}).
Let $\hX$ be an affine variety endowed with a $\C^{*}-$action
$t.$ Then $t$ induces a grading $\hA=\oplus_{n\in \Z}\hA^n$
on the algebra
$\hA=\C[\hX]$ of regular functions on $\hX$ where $\hA^n:=
\{f\in \hA\;\vert\; f\circ t_\lambda=\lambda^nf\}$ consists of
the quasi-invariants of weight $n$ of $t.$

Vice versa, given a grading $\hA=\oplus_{n\in \Z}\hA^n$ of $\hA=\C[\hX]$,
one can define a $\C^{*}-$action on $\hA$ by setting
$t_\lambda(f_n)=\lambda^nf_n$ for $f_n\in \hA^n,$ $n\in \Z$,
and extending it to the whole $\hA$ in a natural way. If $\hA$
is finitely generated,
then it also has a finite system of homogeneous generators
$(f_{n_1},\ldots,f_{n_k}),$
$\,\,f_{n_i}\in \hA^{n_i}.$ The morphism
$F=(f_{n_1},\ldots,f_{n_k}):\hX\hookrightarrow \C^k$ is an embedding
equivariant
with respect to the linear $\C^{*}-$action $t_\lambda(x_1,\ldots,x_k)=
(\lambda^{n_1}x_{n_1},\ldots,\lambda^{n_k}x_{n_k})$ on $\C^k$
and the induced $\C^{*}-$action on $\hX.$

\bigskip

{\bf Gradings and locally nilpotent derivations}
(see e.g. \cite{ML 2, KaML 3, De}).

\medskip

\defin
\label{dfpd}
Let $\partial \in \LND(A)\setminus \{0\}$ where $(A,\,F)$ is a filtered
domain. Suppose that

\smallskip

\noindent ${\bf (*)}$ {\it there exists $k\in \Z$ such that
$\partial \fr{i}\subset \fr{i+k}$ for all $i\in \Z.$}

\smallskip

\noindent Denote by deg$\,\partial = k_0$ the minimal such $k$.  Define
$\widehatpar=\gr \partial : \Gr A\longrightarrow \Gr A$ as follows:
for $f\in \fr{i}\setminus \fr{i-1},$
set $\widehatpar \baf=\partial f+\fr{i+k_0-1},$ and then naturally
extend $\widehatpar$ to the whole algebra $\Gr A.$

\medskip

\exrs (\nexr) \label{exr9.9} Given $\partial \in LND(A) \setminus \{0\},$
verify that $\widehatpar \in \LNDG(\Gr A)\setminus \{0\}$ where $\LNDG(\hA)$
denotes the set of all homogeneous locally nilpotent
derivations of a graded algebra $\hA=\oplus_{n\in \Z}\hA^n.$

\smallskip

(\nexr) \label{exr9.10}
Suppose that a filtered domain $(A,\,F)$ is finitely generated.
Show that, given
$\partial \in \LND(A)\setminus \{0\},$ the condition $(*)$ of Definition
\ref{dfpd} above is fulfilled.

\smallskip

(\nexr) \label{exr9.11}
Let $\hA=\oplus_{n\in \Z}\hA^n$ be a graded algebra.
Show that, given any locally nilpotent derivation $\widehatpar\in
\LNDG(\hA),$ there exists $k_0=k_0(\widehatpar)\in \Z$ called the {\it degree}
of $\widehatpar$ such that $\widehatpar(\hA^n)\subset \hA^{n+k_0}.$
Furthermore, show that, if
$a=\sum_{i=k}^{k+l}a_i\in \Ker \widehatpar=\hA^{\widehatpar}$
where $a_i\in \hA^i,$ then
$a_i\in \hA^{\widehatpar},$ $i=k,\ldots,k+l.$
Therefore, $\hA^{\widehatpar}$ is
a graded subalgebra of the graded algebra $\hA.$

\smallskip

(\nexr) \label{exr9.12} Let $A$ be an integral domain. Given $\partial \in
\LND(A),$ set $\dgp a=n$ if $\partial^{n+1}a=0$ and $\partial^na\ne 0;$
$\dgp 0=-\infty.$ Verify that:

\smallskip

\noindent $(a)$ $\dgp$ is a degree function on $A$ over $\N;$

\smallskip

\noindent $(b)$ if $A$ is a $\C-$algebra, then
$\partial \lambda=0$ for any $\lambda \in \C;$

\smallskip

\noindent $(c)$ the equality $\partial (ab)=0$ where $a,\,b\in A\setminus\{0\}$
 implies that $\partial a=\partial b=0;$

\smallskip

\noindent $(d)$ the equality $\partial (a^k + b^l) = 0$ where\footnote{This
observation is due to L. Makar-Limanov; it is essential in his new approach to
the Theorem \ref{thmLML} \cite{ML 3}.}
$a^k + b^l \neq 0$ and $k,\,l \ge 2$ implies that $\partial a=\partial b=0.$

\medskip

(\nexr) \label{exr9.13}
Let $A=\C[x,\,y].$ Consider the $\C_{+}-$action
$\varphi_{\lambda}\,:\,(x,y)\longmapsto (x,y+\lambda x^2)$
on $\C^2.$ Let $\partial_{\varphi}$ be the locally nilpotent derivation
which corresponds to $\varphi$. Show that
$\partial_{\varphi}(x)=0,$ $\,\partial_{\varphi}(y)=x^2$ and
$\partial^2_{\varphi}(y)=0.$ Deduce that $\dgv x=0,\quad \dgv y=1$
for the associated degree function $\dgv$ on $A,$ and so, that
$\dgv f=\dgy f$ for any polynomial $f=f(x,y) \in A.$

\bigskip

{\bf The Brody hyperbolicity.} Recall that a complex manifold
$M$ which does not contain entire curves, that is, does not
admit non-constant holomorphic mappings $\C \to M,$ is called {\it Brody
hyperbolic.} By analogy with this, let us introduce the following notion.

\medskip

\defin We say that a quasi-projective variety $X$ is
{\it algebraically Brody hyperbolic} if any morphism $\C \to X$ is constant.
A regular function $f$ on $X$ is {\it (algebraically) Brody hyperbolic} if
so are its generic fibres  $f^{-1}(c),\,\,c\in \C.$

\medskip

\exrs (\nexr) \label{exr9.14}
Show that any non-constant polynomial of one variable
is Brody hyperbolic.
Verify that the polynomials $xy$ and
$x^k + y^l \in \C^{[2]}$ where $k,\,l \ge 2$ are Brody hyperbolic.

\smallskip

(\nexr) \label{exr9.15}
More generally, show that a polynomial $p \in \C^{[2]}$ fails to be
Brody hyperbolic iff it can be linearized, that is, $p \circ \alpha = x$
for some automorphism $\alpha \in \Aut \C^2.$

\bigskip

The next lemma is inspired by Exercises \ref{exr9.5} and
\ref{exr9.12}$(c),\,(d).$ It will be used in the proof of Derksen's
Theorem \ref{derks} in section 7.3 below.

\medskip

\lemma {\it $(a)$ Let $X$ be a quasi-projective variety,
and let $f \in A:=\C[X]$ be a Brody hyperbolic regular
function on $X.$ Then $f$
cannot be a constant of a locally nilpotent derivation
$\partial \in LND(A).$

\smallskip

\noindent $(b)$ Let $A$ be an integral domain over
$\C$ with the unit element ${\bf e} \in A,$ and let $p \in \C^{[k]}$ be a
Brody hyperbolic polynomial. If for some elements $a_1,\dots,a_k \in A$
we have $b:=p(a_1,\dots,a_k) \in A^{\partial}\setminus \C {\bf e},$
then $a_1,\dots,a_k \in A^{\partial}.$}

\label{lmBH}
\medskip

\proof $(a)$ Assume the contrary. Being an invariant of the associated
$\C_+-$action $\varphi_{\partial}$ on $X$ (see Exercise \ref{exr9.1} above)
the function $f$ is constant on the orbits of this action. Thus, generic
$\varphi_{\partial}-$orbits are $\C-$curves contained in the fibres of $f,$
in contradiction with the assumption that the function $f$ is Brody hyperbolic.

\smallskip

$(b)$ Let $A'$ be the subalgebra of the algebra $A$ generated by the elements
$a_1,\dots,a_k$ and all its successive $\partial-$derivatives. It is finitely
generated (henceforth, this is an affine domain) and $\partial-$invariant.
Consider the associated $\C_+-$action
$\varphi_{\partial}$ on the affine variety  $X':= \spec A'.$ By the assumption,
the nonconstant regular function
$b=p(a_1,\dots,a_k) \in A'=\C[X']$ is $\varphi_{\partial}-$invariant,
and so, it is constant along the $\varphi_{\partial}-$orbits.
Therefore, the image $\alpha(O_x)$ of the $\varphi_{\partial}-$orbit of a
generic point $x \in X'$ under the regular mapping
$\alpha\,:\,X' \to \C^k,\,\,x
\longmapsto (a_1(x),\dots,a_k(x)),$ is contained in a
generic fibre $p^{-1}(c)$
of the polynomial $p.$ Since this polynomial is assumed
to be Brody hyperbolic, the mapping $\alpha$ is constant
along the generic $\varphi_{\partial}-$orbits. It follows that its coordinate
functions $a_1,\dots,a_k$ are $\varphi_{\partial}-$invariants. Whence,
$a_1,\dots,a_k \in A^{\partial}$ (see again Exercise \ref{exr9.1}).
\qed

\medskip

\rems \nrem Suppose that under assumptions of Lemma \ref{lmBH}$(b),$
all the fibres of the polynomial $p$ except the zero one $p^{-1}(0)$ are Brody
hyperbolic. Then the conclusion of the lemma holds if one
replaces the condition $b \in A^{\partial}\setminus \C {\bf e}$ by the weaker
one $b\in A^{\partial}\setminus \{0\}.$

\smallskip

\nrem In Exercises \ref{exr9.1}, \ref{exr9.12}$(c),\,(d)$ one can make use of
Lemma \ref{lmBH}$(b),$
of the previous remark and of Exercise \ref{exr9.14}.

\subsection{Gradings and LND's on Russell's cubic}

In this section we provide a proof of Theorem \ref{derks}.
We use the following notation.

\medskip

{\bf Notation.} From now on $A=A_0=\C[x,y,z,t]/(p_0)$
where $p_0=x+x^2y+z^3+t^2,$
will be the algebra of regular functions on the Russel cubic
threefold $X_0\subset \C^4.$ Consider the weight degree function
$\dg x=-1, \dg y=2, \dg z=\dg t=0$
on the polynomial ring $\C^{[4]}.$ It is easily
seen that $\dg$ is appropriate for the principal
ideal $I:=(p_0)$ (see Definition \ref{defWDF}).
Hence, by Exercise \ref{exr9.7}, it induces a
degree
function $d_A$ on $A$ which, in turn, defines a filtration on $A$
(in what follows we use the notation $\dg$ instead of $d_A$). Let
$\hA:= \Gr A$ be the associated graded algebra.
{}From Lemma \ref{lm9.1} we obtain such a corollary.

\bigskip

\cor
{\it $\hA\simeq \C[\hx,\hy,\hz,\htt]/(q_0)$ where
$q_0(\hx,\hy,\hz,\htt)=\hx^2\hy+\hz^3+\htt^2,$ i.e. $\hA=\C[\hX_0],$
where $\hX_0=\{\hx^2\hy+\hz^3+\htt^2=0\}.$}

\bigskip

The next lemma exploits the specific equation of the Russell cubic.

\medskip

\nlemma{\cite{De}.} $(a)$ {\it Any element $f\in A$ has a unique presentation
of the form
$$f=a(x,z,t)+yb(y,z,t)+xyc(y,z,t)\,|\,X\,$$
where $a,b,c$ are polynomials.}

\noindent $(b)$ {\it We have:

\noindent $(i)$ $\hf=\hx^rh(\hz,\htt)$ iff $\dg f=r\le 0;$

\noindent $(ii)$ $\hf=\hy^rh(\hz,\htt)$ iff $\dg f=2r>0;$

\noindent $(iii)$ $\hf=\hx\hy^rh(\hz,\htt)$ iff $\dg f=2r-1>0,$

\noindent where $h(\hz,\htt)\in \C[\hz,\htt].$}

\label{lm9.2}
\medskip

{\bf Proof of $(b).$} Since the degree $\dg[yb(y,z,t)]$ is even and positive
when $b\ne 0,$
$\dg[xyc(y,z,t)]$ is odd and positive when $c\ne 0,$
and $\dg a(x,z,t) \le 0,$ we have that
in the case $(i)$ $f=a(x,z,t)$ and hence $\hf=\hx^r h(\hz,\htt);$
in the case $(ii)$ $\,\,\hf= {\rm gr}\,[yb(y,z,t)]=\hy^rh(\hz,\htt);$
finally, in the case $(iii)$
$\,\,\hf= {\rm gr}\,[xyc(y,z,t)]=\hx\hy^rh(\hz,\htt)$
for some $h(\hz,\htt)\in \C[\hz,\htt].$
\qed

\medskip

\cor {\it We have:

\smallskip

\noindent $\hA^0=\C[\hz,\htt],$ and thus $\hA$ is a
$\C[\hz,\htt]-$module);

\smallskip

\noindent $\hA^i=\hx^{-i}\C[\hz,\htt]$ for $i\le 0;$

\smallskip

\noindent $\hA^{2r}=\hy^r\C[\hz,\htt]$ for $r>0,$ and

\smallskip

\noindent $\hA^{2r-1}=\hx\hy^r\C[\hz,\htt]$ for $r>0.$}

\label{cor9.2}
\bigskip

\lemma {\it Consider a surface
$S_c=\{c\hx^2+\hz^2+\htt^3=0\}\subset \C^3$ where $c \in \C^*.$
Then for any nonconstant polynomial $h \in \C[z,\,t]$
the regular function
$f:=h\,|\,S_c \in \C[S_c]$ on the surface $S_c$ is Brody hyperbolic.}

\label{lm9.3}
\medskip

\proof\footnote{We give a simplified proof suggested by Sh.
Kaliman.} We may assume that the generic fibres
of the polynomial $h$ are
irreducible plane curves. Indeed, otherwise the polynomial $h$ has
the Stein factorization $h = \varphi\circ h_1$ where
$\varphi \in \C[z],\,\,\deg \varphi\ge 2,$ and
$h_1 \in \C[z,\,t]$ is a
polynomial with irreducible generic fibres, and we may replace $h$ by the
polynomial $h_1.$

Assume, on the contrary, that for some
$c\ne 0$ the regular function $f=h\,|\,S_c$ is not
Brody hyperbolic, that is, its generic fiber
$F_{\lambda}  := f^{-1}(\lambda) \subset S_c$
admits a nonconstant morphism $\C \to F_{\lambda}.$
It follows that an irreducible component, say, $F_{\lambda}'$ of the curve
$F_{\lambda}$
is a smooth affine curve isomorphic to $\C$
(see Exercise \ref{exr9.5} above).

The curve $F_{\lambda}$ is represented as a two-sheeted ramified covering
of $C_{\lambda}$ under the projection
$\pi\,:\,S_c \to \C^2,\,\,\,(\hx,\,\hz,\,\htt) \mapsto (\hz,\,\htt).$ Thus,
the irreducible affine plane curve $C_{\lambda}$
admits a dominant morphism $\C \simeq F_{\lambda}'\to C_{\lambda},$
and hence, it is
isomorphic to $\C,$ too.

We have the following alternative: either

$\bullet$ The generic fibre $F_{\lambda}$
of the regular function $f$ is irreducible, i.e. $F_{\lambda}=F_{\lambda}'
\simeq \C,$ and then
$p:=\pi\,|\,F_{\lambda}'\,:\,F_{\lambda}'  \simeq \C \to C_{\lambda} \simeq \C$
is a two-sheeted branched covering;

or

$\bullet$  The generic fibre $F_{\lambda}$ is reducible, and then the
mapping $p$ as above is univalent, and hence, isomorphic.

In the first case a generic curve $C_{\lambda}$
meets the ramification locus
$\Gamma_{2,\,3} = \{\hz^2+\htt^3=0\} \subset {\bf C}^2$ of
the projection $\pi$ at one point which corresponds to the
only critical value of the quadratic polynomial $p \in \C[u].$
Therefore, the restriction
$h_1\,|\,\Gamma_{2,\,3}\,:\,\Gamma_{2,\,3} \to \C$
is generically one-to-one, which is impossible.

In the second case a generic curve $C_{\lambda}$ does not meet
the ramification locus $\Gamma_{2,\,3}$ at all, and
then it should be contained in an elliptic curve
$\hz^2+\htt^3=$const$\neq 0,$ which is impossible either.
This completes the proof.
\qed

\bigskip

Lemmas \ref{lmBH}$(a)$ and \ref{lm9.3} yield the following corollary.

\medskip

\cor {\it The surface $S_c$ where $c \neq 0$
does not admit a non-trivial $\C_{+}-$action with a
non-constant invariant function $h(\hz,\htt).$}

\label{cor9.3}
\bigskip

\lemma {\it For any nonzero locally nilpotent derivation
$\widehatpar \in \LNDG(\hA)$
and for any homogeneous element $\hf\in \hA^n$
the equality $\widehatpar \hf=0$ implies that $n:=\dg \hf\le 0.$}

\label{lm9.4}
\medskip

\proof Assume, on the contrary, that $n>0.$ Suppose first
that $n$ is odd, i.e. $n=2r-1$ for some $r\in \N.$
Then, by Lemma \ref{lm9.2}$(b),$ we have that $\,\,\hf=\hx\hy^rh(\hz,\htt)$
for some non-zero polynomial
$h(\hz,\htt)\in \C[\hz,\htt],$
and the equality $\widehatpar \hf=0$
implies that $\widehatpar \hx=\widehatpar \hy=\widehatpar h=0$
(see Exercise \ref{exr9.12}$(c)$ above).
Thus, $\hx,\,\hy\in \hA^\widehatpar$ are invariants of the
associated $\C_{+}-$action $\varphi=\varphi_{\widehatpar}$ on
$\hX_0=\{\hx^2\hy+\hz^3+\htt^2=0\}$ (see Exercise \ref{exr9.1} above).
Therefore, each orbit of the $\C_{+}-$action $\varphi$ is contained
in a curve
$\Gamma_{c_1,c_2}=\{\hx=c_1,\,\hy=c_2\} \subset \hX_0.$
Conversely, any such curve in $\hX_0$
consists of $\C_{+}-$orbits. Since a generic curve
$\Gamma_{c_1,c_2} = \{\hz^3+\htt^2=-c_1^2c_2\}$ is elliptic,
it does not admit
an embedding of $\C$ (see Exercise \ref{exr9.5} above),
and hence all the points of
$\Gamma_{c_1,c_2}$ are fixed by the $\C_{+}-$action $\varphi.$
It follows that the $\C_{+}-$action $\varphi$ on $\hX_0$
is trivial, i.e. $\widehatpar =0,$ a contradiction.

Now consider the case when $n > 0$ is even,
i.e. $n=2r$ for some $r\in \N.$ In this case by Lemma \ref{lm9.2}$(b),$ we have
$\hf=\hy^rh(\hz,\htt)$ where $h(\hz,\htt)\in \C[\hz,\htt]$ and $h\ne 0.$
Hence, $\widehatpar \hy=0,$ and so $\hy\in \hA^\widehatpar$ is an invariant of
the $\C_{+}-$action $\varphi$ on $\hX_0$ (see Exercise \ref{exr9.1}).
Thus, $\hX_0$ is foliated
by the $\varphi-$invariant surfaces $S_c=\{\hy=c\}, \,\,c\in \C.$
Denote by $\partial_c\in \LND(A_c)$ the corresponding locally nilpotent
derivation on
$A_c=\C[S_c]$, that is, the infinitesimal
generator of the $\C_{+}-$action $\varphi_c:= \varphi \,\vert \,{S_c}$
on ${S_c}.$
For a generic $c\in \C,$ $\varphi_c$ is a non-trivial
$\C_{+}-$action on $S_c,$ whence $\partial_c\ne 0.$

Next we show that there exists a non-constant $\varphi-$invariant function
$h_1(\hz,\htt)\in \hA^\widehatpar.$
Indeed, since $\trdg \hA^\widehatpar=
\trdg \hA-1=2$ (see Exercise \ref{exr9.2} above),
the subalgebra $\hA^\widehatpar$ contains a function
$g$ such that $\hy$ and $\hg$ are algebraically independent;
in particular,
$\hg\notin \C[\hy].$ Furthermore, $\hg$ and $\hy$ are both
$\varphi-$invariants,
and so, for $s \in \N$ sufficiently large $\hg\hy^s$ is a
$\varphi-$invariant of a positive degree. We have proven above that
the equality
$\hg\hy^s=\hx\hy^rh_1(\hz,\htt)$ is impossible. Hence, we get
$\hA^\widehatpar \ni \hg\hy^s = \hy^rh_1(\hz,\htt)$ for some $r > 0$ where
$h_1 \in \hA^\widehatpar$ is non-constant.
Finally, the restriction of the polynomial $h_1$ onto a generic surface
$S_c$ is a nonconstant invariant of the non-trivial $\C_{+}-$action
$\varphi_c$ on $S_c,$ which contradicts to Corollary \ref{cor9.3}.
The proof is completed.
\qed

\bigskip

\cor {\it $(a)$ $A^\partial\subset \fr{0}$ for any non-zero locally nilpotent
derivation $\partial\in \LND(A).$

\smallskip

\noindent $(b)$ $\Dk (A)\subset \fr{0} \ne A.$}

\label{cor9.4}
\medskip

\proof The statement of $(b)$ follows from $(a)$
in virtue of Definition \ref{defIN}. To prove $(a)$ assume
the contrary, i.e. that for some $f\in A$ where $\dg f>0$
and for some $\partial\in \LND\setminus \{0\}$ we have $\partial f=0.$
Then
$\dg \hf >0$ as well, and $\widehatpar \hf=0$
(see Definition \ref{dfpd} and Exercise \ref{exr9.9}). It remains to apply
Lemma \ref{lm9.4}.
\qed

\bigskip

Now the proof of Theorem \ref{derks} is completed.

\section{Concluding remarks}

In this section as in section 4.2 above we mainly follow \cite{KaZa}.
To begin with, remind the following problem.

\bigskip

\noindent {\bf Generalized Serre Problem}:

\smallskip

\noindent {\it Let $X$ be a smooth contractible affine variety.
Is any algebraic vector bundle on $X$ trivial?}

\smallskip

Due to the Quillen-Suslin Theorem,
this is true for the affine spaces.

\medskip

The next question suggested by A. Beilinson generalizes
the Zariski Cancellation Problem. Following \cite[sect. 4]{BaWr},
we say that a closed subvariety $Y$ of an
affine variety $X$ is a {\it retract} of $X$ if there exists a morphism
$f\,:\,X \to Y$ such that $f\,|\,Y = {\rm id}_Y.$ An abstract variety $Y$
is called a {\it retract} of $X$ if the image of $Y$ under some proper
embedding $Y \hookrightarrow X$ is a retract of $X.$

\smallskip

\noindent {\it Is it true that an affine variety which is a retract of
an affine space is isomorphic to an affine space?

More
generally, let $R$ be a commutative ring and $A$ be an
$R-$algebra such that any surjective map of $R-$algebras
$B \longrightarrow\!\!\!\to A$ admits a right inverse.
Is it true that $A$ is isomorphic to ${\rm Sym}\,P$ where
$P$ is a projective $R-$module?}

\medskip

\exrs (\nexr) \label{exr0.1} Show that if $X$ is a retract of an affine space,
then it is a smooth contractible affine variety, $\ok(X)=-\infty,$
and, moreover, any  algebraic vector bundle on $X$ is trivial.

\smallskip

(\nexr) \label{exr0.2} Assume that a variety $X$ is a retract of an affine
space.
Show that then $X$ is a retract of any affine variety which contains a copy of
$X.$ Deduce that an affine space is a retract of any affine variety which
contains a copy of it. In particular, if a subvariety $X \subset \C^n$ is
isomorphic to $\C^k,$ then $X$ is a retract of $\C^n$ (under the given
embedding).

\smallskip

(\nexr) \label{exr0.3} Let $f\,:\,X \hookrightarrow \C^n$ and
$g\,:\,X \hookrightarrow \C^k$ be two closed embeddings. Verify that if
$f(X)$ is a retract of $\C^n,$ then also $g(X)$ is a retract of $\C^k.$

\smallskip

\noindent {\it Hint.} Apply the Jelonek-Kaliman-Nori-Srinivas Theorem on
equivalence of embeddings cited in the Introduction.

\smallskip

(\nexr) \label{exr0.4} Find an example of a pair of affine varieties $X$ and
$Y$ and a pair of proper embeddings $f,\,g\,:\,Y \hookrightarrow X$ such that
the image $f(Y)$ is a retract of $X$ whereas the image $g(Y)$ is not.

\bigskip

{\bf Notation.} Let $X$ be a smooth affine $n-$fold,
and $A = \C[X]$ be the algebra of regular functions on $X.$
Set $ML(X)=ML(A)$ and $Dk(X) = Dk(A)$ (see Definition \ref{defIN}).

\smallskip

We would like to
know when these invariants are nontrivial. For instance,

\smallskip

\noindent {\it Is it true that $ML(X_0 \times \C) \neq \C$ where $X_0$ is the
Russell cubic?}

\smallskip

\noindent It is unknown {\it whether this product is, indeed,
an exotic $\C^4.$}

\medskip

\rem In general, the Makar-Limanov invariant is not invariant under
cancellation\footnote{This observation is due to T. Bandman and L.
Makar-Limanov.}. Indeed, consider the {\it Danielewski surfaces}
$D_n=\{x^ny+z^2=1\}$ in $\C^3$ (see Introduction), and set $X_n = D_n \times
\C.$ Then we have:
$X_1 \simeq X_2$ \cite{Dan}, and so, $ML(X_1) = ML(X_2),$ whereas
$ML(D_2) \neq ML(D_1)=\C$ (see \cite{ML 1, ML 4}, and also
Corollary \ref{corTR} below).

\medskip

\exrs (\nexr) Verify that $D_1 \simeq (\Pp^1 \times \Pp^1) \setminus \Delta$
where $\Delta \subset \Pp^1 \times \Pp^1$ is the diagonal. Find the minimal
SNC-compactifications for the other Danielewski surfaces $D_n$ and their dual
graphs.

\smallskip

(\nexr) Deduce that the fundamental groups at infinity $\pi_1^{\infty}(D_n)$
are pairwise non-isomorphic, and hence the surfaces $D_n$ are pairwise
non-homeomorphic. What are the homology groups $H_*(D_n;\,\Z)?$ (see
Fieseler \cite{Fi}, tom Dieck \cite{tD 3}).

\smallskip

(\nexr) Show that $\ok(X) = -\infty$ for any smooth quasiprojective variety $X$
which admits a non-trivial regular $\C_+-$action. Deduce that
$\ok(D_n) = -\infty$ for all $n.$ Observe that the latter conclusion also
follows by the Iitaka-Fujita Strong Cancellation Theorem \ref{thif}.

\bigskip

{\bf Notation.} To any regular $\C_+-$action $g$ on $X$
there corresponds a one-parameter algebraic subgroup $G$ of the automorphism
group
$\Aut X.$ Denote $\Aut_+ X$ the subgroup of the group $\Aut X$ generated
by all such $\C_+-$subgroups.

\smallskip

The answer to the following question is unknown:

\smallskip

\noindent {\it Wherther or not $\Aut_+ \C^n = \Aut \C^n?$}

\smallskip

\noindent (see e.g. \cite{For} on some density properties of this subgroup).

\smallskip

It is easily seen that for any
$n \ge 2$ the group $\Aut_+ \C^n$ acts
$m-$transitively on $\C^n,$
i.e. it is transitive on the set of all $m-$tuples of distinct points
in $\C^n,$ whatever $m \in \N$ is. In such a case
we say that the action is {\it infinitely transitive}.

Clearly, the Makar-Limanov invariant is trivial, i.e. $ML(X)=\C,$
if the group $\Aut_+ X$ is transitive on a variety
$X$ or, at least, has a dense orbit.
The next surprising examples were found in \cite{ML 1}
for $k=1$ and in  \cite{KaZa} for $k>1.$

\medskip

\nthm{\cite[Thm. 6.1]{ML 1, KaZa}.} {\it Let
$p \in \C^{[k]},\,\,\,k\ge 1,$
be a non-constant polynomial, and let $X_p \subset \C^k$ resp. $Y_p\subset
\C^{k+2}$ be the hypersurface given by the equation $p(\bx)=0$ resp.
$uv-p(\bx)=0$ where $\bx = (x_1,\dots,x_k).$
Then the following statements hold.

\smallskip

$(a)$ The hypersurface $Y_p$ is smooth iff
$X_p$ is a smooth reduced hypersurface. If so, then the variety
$Y_p$ is simply connected, and there is
an isomorphism of the reduced homology groups
${\tilde H}_*(Y_p;\,\Z) \simeq {\tilde H}_{*-2}(X_p;\,\Z).$

\smallskip

$(b)$  The group $\Aut_+ Y_p$ is transitive on the smooth part $\reg Y_p$
of the variety $Y_p;$ moreover, for $k > 1$ this group is infinitely
transitive on $\reg Y_p.$}

\label{thmTR}
\bigskip

\cor {\it $(a)$ Suppose that the hypersurface $X_p\subset \C^k$ is
smooth and reduced. Then the hypersurface
$Y_p\subset \C^{k+2}$ is contractible iff the
variety $X_p$ is acyclic.

\smallskip

$(b)$ We have\footnote{At the same time, by Theorem \ref{thmTR}$(b)$,
the variety $Y_p$ may have rather non-trivial topology.}
$ML(Y_p) = \C.$}

\label{corTR}
\bigskip

\rems \nrem
There is such a question:

\smallskip

\noindent {\it Suppose that an affine algebra $A$ over $\C$ possesses
$n=\trdg A$ locally nilpotent derivations
which are linearly independent over $A.$
Does it follow that $A \simeq \C^{[n]}?$}

\smallskip

\noindent The answer is negative. Indeed,
following the lines of the proof of Theorem 6.1 in \cite{KaZa}
it is easily seen that smooth non-contractible varieties $Y_p$ as above
provide corresponding counterexamples.

However, we do not know what is the answer under a stronger assumption
that the corresponding regular vector fields on the variety
$X:={\rm spec}\, A$
are linearly independent at any point of $X.$
It is positive provided that, in addition, the
given locally nilpotent derivations pairwise commute.
Indeed, then the variety $X$ is the unique
orbit of the associated regular  free
action of the additive group $\C_+^n$ on $X.$

\medskip

\nrem
For $k=1$ the only example
of a smooth contractible surface of type $Y_p$ is
the affine plane $\C^2.$ Furthermore, as follows from
the Abhyankar-Moh and Suzuki Theorem \cite{AM, Suz 1}
(see Theorem \ref{cusuthm}$(a)$ above), for $k=2$
every smooth contractible $3-$fold in $\C^4$ of
type $Y_p$  is isomorphic to the affine space $\C^3,$
and the embedding $\C^3 \simeq Y_p \hookrightarrow \C^4$
is rectifiable (see \cite[Prop. 5.2$(b)$]{KaZa}).

Starting with $k=3$ one can obtain a number of new examples of
smooth contractible hypersurfaces of type $Y_p$ in $\C^{k+2}$ choosing for
$X_p$ smooth acyclic (not necessarily contractible) hypersurfaces in
$\C^k$
(see  e.g. Example \ref{ex7.2} above of such a surface
$X_p=Y_{k,l,s}\subset \C^3$).
Presumably, among these hypersurfaces $Y_p$ there are examples of exotic
$\C^{n},\,\,n \ge 4,$ with infinitely transitive automorphism groups
(for instance, the hypersurface $uv - p_{k,l}(x,\,y,\,z) =0$ in $\C^5$ where
$p_{k,l} \in \C^{[3]}$ is a tom Dieck-Petrie polynomial should be such one).
This would provide a negative answer to a question proposed by the author in
\cite[(2.8)]{OPOV}.
Thus, it would be important to know

\smallskip

\noindent {\it For which polynomials $p \in \C^{[k]}$ the variety $Y_p$ is an
exotic $\C^{k+1}?$}

\smallskip

\noindent Corollary \ref{corTR}$(b)$ above shows that
the Makar-Limanov invariant does not supply an
answer to this question.

\medskip

\nrem Furthermore, these examples are interesting in connection with
the Abhyankar-Sathaye Embedding Problem (see Introduction).
Indeed, consider the regular function
$f := u\,|\,Y_p \in \C[Y_p].$
For $c \neq 0$ the level hypersurface $U_c = f^{-1}(c)$ is isomorphic to
$\C^k;$
further, the complement $Y_p \setminus U_0$ is isomorphic to the cylinder
$\C^* \times \C^k.$
At last, the hypersurface $U_0$ is isomorphic to the product $X_p \times \C.$
Actually, the hypersurface $U_0$ coincides with the exceptional divisor $E'$
of the affine modification
$$\si\,:\,Y_p \to \C^{k+1},\,\,(\bx,\,u,\,v) \longmapsto (\bx,\,v),$$
of the affine space $\C^{k+1}$ along the hyperplane $D:=\{v=0\}$ with center
$X_p\subset D$ (see Example \ref{exCI}).
The following question arises:

\smallskip

\noindent {\it Is it true that $Y_p \simeq \C^{k+1}$ iff
$X_p \simeq \C^{k-1}?$}

\smallskip

\noindent It has a direct relation to the Zariski Cancellation Problem.
Namely, it would be useful to know

\smallskip

\noindent {\it Whether the implication $X \times \C \simeq \C^k$
$\Longrightarrow$
$X  \simeq \C^{k-1}$ holds at least for hypersurfaces $X$ in $\C^k?$}

\medskip

\nrem Assume that for some polynomial $p \in \C^{[k]}$ we have
$U_0 = X_p \times \C \not\simeq \C^k$ whereas $Y_p \simeq \C^{k+1}.$
Then, evidently, the embedding
$\C^k \simeq U_1 \hookrightarrow Y_p  \simeq \C^{k+1}$
is not rectifiable, thus providing a counterexample to
the Abhyankar-Sathaye Embedding Problem.
Otherwise, i.e. in the case where $U_0  \not\simeq \C^k$ and
$Y_p  \not\simeq \C^{k+1},$ we obtain an example showing that
the analog of the Miyanishi Theorem \ref{thmMI}
does not hold in dimension $n = k+1$ (clearly, here $n\ge 4$).

\bigskip

The following fact, related to
the Abhyankar-Sathaye Embedding Problem,
was established in Sathaye \cite{Sat} and Wright \cite{Wr 2}
in the case where $X \simeq \C^2,$ and in
\cite{KaZa} in the present more general form.

\medskip

\nthm{\cite[Theorem 4.2]{Sat, Wr 2, KaZa}.} {\it Let
$X=X_{n,\,f,\,g}$ be an irreducible smooth surface in
$\C^3$ given by the equation
$f(x,\,y)z^n + g(x,\,y) = 0$ where
$f,\,g \in \C[x,\,y],\,\,n \in \N.$
Then the following conditions $(i)\, -\, (iv)$
are equivalent:

\smallskip

\noindent $(i)$ $e(X)=1$ and $H_1(X;\,\Z) = 0$
where $e(X)$ denotes the Euler characteristic of the variety $X.$

\smallskip

\noindent $(ii)$ The surface $X$ is acyclic, i.e. ${\tilde H}_*(X;\,\Z) = 0.$

\smallskip

\noindent $(iii)$ $X \simeq \C^2.$

\smallskip

\noindent $(iv)$ The surface $X$ is rectifiable, i.e. it can be
transformed into a plane by an automorphism of $\C^3.$

\smallskip

\noindent For $n > 1$ these conditions are equivalent
also to the following one:

\smallskip

\noindent $(v)$ The pair $(f,\,g)$ is rectifiable, i.e. it can be
transformed into a pair $(\alpha(f),\,\alpha(g)) = (p(x),\,y)$
by an automorphism $\alpha \in \Aut\C^{[2]}.$}

\label{thmSW}
\smallskip

\bigskip

Quite recently, following the same idea as in the proof of
Theorem \ref{thmSW} in \cite{KaZa}, S. Venereau has obtained such a
result\footnote{We place it here with his kind permission.}.

\medskip

\nprop{({\bf Venereau}).} {\it Let
$q \in \C[\bx,\,u]$ (where $\bx = (x_1,\dots,x_k)$)
be a polynomial such that
$p(\bx) := q(\bx,\,0) \in  \C[\bx]$ is a non-constant polynomial.
Consider the hypersurfaces $X_p := \{p(\bx) =0\} \subset \C^k$ and
$Y_{q,\,n} := \{u^nv -  q(\bx,\,u)=0\}\subset \C^{k+2}.$
If the hypersurface $X_p$ is rectifiable in $\C^k,$ then also the hypersurface
$Y_{q,\,n}$ is rectifiable in $\C^{k+2}.$}

\label{prV}
\bigskip

We do not know whether the converse is true.
Notice that under the conditions
of Proposition \ref{prV} an analog of
Theorem \ref{thmTR}$(a)$ above holds.
See also \cite[ Prop. 5.3, 5.4]{KaZa} on extensions of
Theorem \ref{thmTR}$(a)$ to complete intersections of higher codimensions.

\bigskip

\rem Every tom Dieck-Petrie surface $X_{k,l}$ (see Examples \ref{ex3.1},
\ref{exDP} and \ref{ex6.1} above) admits at least two non-equivalent embeddings
into $\C^3$ \cite[Example 6.3]{KaZa}.
The question arises:

\smallskip

\noindent {\it What is the cardinality of the set of all pairwise
non-equivalent proper embeddings $X_{k,l} \hookrightarrow \C^3?$ Are those
embeddings rigid or, conversely, do admit non-trivial deformations up
to the action of the automorphism group $\Aut \C^3$ on $\C^3?$}

\smallskip

\noindent The same question has sense for any contractible or acyclic surface
embeddable into $\C^3.$

\noindent  Mikhail Zaidenberg

\noindent  Universit{\'e} Grenoble I

\noindent  Institut Fourier

\noindent  UMR 5582 CNRS-UJF

\noindent  BP 74

\noindent  38402 St. Martin d'H{\'e}res--cedex

\noindent  France

\noindent  e-mail: zaidenbe@mozart.ujf-grenoble.fr


\bigskip

\begin{thebibliography}{999}
\addcontentsline{toc}{a}{References}



\bibitem[AEH]{AEH} S. Abhyankar, P. Eakin, W. Heinzer,
{\em On the uniqueness of
the coefficient ring in a polynomial ring}, J. Algebra, {\bf 23} (1972),
310--342.
\bibitem[AM]{AM} S.S. Abhyankar, T.T. Moh, {\em Embedding of the line in the
plane}, J. Reine Angew. Math. {\bf 276} (1975), 148--166.
 \bibitem[AF]{AF} A. Andreotti, T. Frenkel, {\em The Lefschetz theorem on
hyperplane sections}, Ann. Math. {\bf 69} (1959), 713--717.
 \bibitem[An]{An} G. Angerm{\"u}ller, {\em Connectedness properties of
polynomial maps between affine spaces}, Manuscr. Math. {\bf 54} (1986),
349-359.
 \bibitem[AAS]{AAS} {\em Automorphisms of affine spaces}, Proc. Conf., July
4-8, 1994, Cura\c cao, Van den Essen (ed.), Kluwer Acad. Publ., Dordrecht e.a.,
1995.
 \bibitem[BD]{BD} G. Barthel, A. Dimca, {\em
On complex projective hypersurfaces which are homology-$P\sb n$'s}, In:
Singularities, Proc Conf.
`Singularities in geometry and topology', Lille (France), 3-8 June, 1991,
J.-P. Brasselet (ed.), Cambridge: Cambridge University Press, Lond. Math.
Soc. Lect. Note Ser. {\bf 201} (1994), 1-27.
 \bibitem[Ba]{Ba} H. Bass, {\em A non-triangular action of $G_a$ on
${\AA}^3$}, J. Pure Appl. Algebra {\bf 33} (1984), 1--5.
 \bibitem[BCW]{BCW} H. Bass, E. Connell, D. Wright, {\em The Jacobian
Conjecture: reduction of degree and formal expansion of the inverse}, Bull.
Amer. Math. Soc. {\bf 7} (1982), 287-330.
 \bibitem[BaHa]{BaHa} H. Bass, W. Haboush, {\em  Linearizing certain reductive
group actions}, Trans. Amer. Math. Soc. {\bf 292} (1985), 463--482.
 \bibitem[BaWr]{BaWr} H. Bass, D. Wright, {\em Localization in the K-theory of
invertible algebras}, J. Pure Appl. Algebra  {\bf 9} (1976), 89-105.
 \bibitem[BCTSSD]{BCTSSD} A. Beauville, J.-L. Colliot-Thelene, J.-J. Sansuc,
P. Swinnerton-Dyer, {\em Vari{\'e}t{\'e}s stablement rationnelles non
rationnelles},
Ann. Math. {\bf 121} (1985), 283-318.
 \bibitem[Bia 1]{Bia 1} A. Bia{\l}ynicki-Birula,
{\em  Remarks on the action of an
algebraic torus on $k^n$, I, II}, Bull. Acad. Polon. Sci. S{\'e}r. Sci. Math.
{\bf 14} (1966), 177--181; {\bf 15} (1967), 123--125.
 \bibitem[Bia 2]{Bia 2} A. Bia{\l}ynicki-Birula,
{\em Some theorems on action of algebraic groups}, Ann. Math. {\bf 98} (1973),
480--497.
\bibitem[Bou]{Bou} N. Bourbaki, {\em Alg{\`e}bre Commutative, Ch. 3}, Hermann,
Paris, 1961.
 \bibitem[Bre]{Bre} G.E. Bredon,
{\em Introduction to compact transformation groups}, Ac. Press, N.Y., 1972.
\bibitem[Bri]{Bri} E. Brieskorn,
{\em Beispiele zur Differentialtopologie von Singularit\"aten},
Invent. Math. {\bf 2} (1966), 1-14.
 \bibitem[ChoDi]{ChoDi} A.D.R. Choudary, A. Dimca, {\em Complex hypersurfaces
diffeomorphic to affine spaces}, Kodai Math. J. {\bf 17} (1994), 171--178.
 \bibitem[Co]{Co} P.H. Cohn, {\em Free rings and their relations},
Second edition, Acad. Press, London e.a., 1985.
 \bibitem[Dan]{Dan} W. Danielewski, {\em On the cancellation problem and
automorphism group of affine algebraic varieties}, preprint, 1989.
 \bibitem[Dav]{Dav} E.D. Davis,  {\em Ideals of the principal class,
$R$-sequences and a certain monoidal transformation}, Pacific J. Math.
{\bf 20} (1967), 197--205.
 \bibitem[De]{De} H. Derksen, {\em
Constructive Invariant Theory and the Linearization Problem},
Ph.D. thesis, Basel, 1997.
\bibitem[tDP 1]{tDP 1} T. tom Dieck, T. Petrie, {\em Contractible affine
surfaces
of Kodaira dimension one}, Japan J. Math. {\bf 16} (1990), 147--169.
 \bibitem[tDP 2]{tDP 2} T. tom Dieck, T. Petrie,
{\em The Abhyankar--Moh problem in dimension 3}, Lect. Notes Math. {\bf 1375},
1989, 48--59.
 \bibitem[tDP 3]{tDP 3} T. tom Dieck, T. Petrie, {\em Homology planes. An
announcement and survey}, Topological Methods in Algebraic Transformation
Groups, Progress in Mathem., {\bf 80}, Birkhauser, Boston, 1989, 27--48.
\bibitem[Di 1]{Di 1} A. Dimca, {\em Hypersurfaces in ${\C}^{2n}$ diffeomorphic
to ${\bf R}^{4n - 2} \,(n \geq 2)$}, Max-Planck Institute, preprint, 1991.
\bibitem[Di 2]{Di 2} A. Dimca, {\em Singularities and Topology of
Hypersurfaces}, Universitext, Springer, 1992.
 \bibitem[Do]{Do} A. Dold, {\em Lectures on algebraic topology}, Springer,
Berlin e.a., 1974.
 \bibitem[Dr]{Dr} L. M. Druzkowski, {\em The Jacobian Conjecture:
survey of some results}, in: Topics in Complex Analysis, Banach Center
Publications, {\bf 31}, Warszawa 1995, 163--171.
 \bibitem[EH]{EH} P. Eakin, W. Heinzer,  {\em A cancellation problem for
rings}, In: Conference on Commutative Algebra (J.W. Brewer, E.A. Rutter, eds.),
Lect. Notes in Mathematics, Springer, Berlin e.a. {\bf 311} (1973), 61--77.
 \bibitem[Ei]{Ei} D. Eisenbud, {\em Commutative algebra with a view towards
 algebraic geometry}, Graduate Texts in Math., Springer, N.Y. e.a., 1994
 \bibitem[Fi]{Fi} K.-H. Fieseler, {\em On complex affine surfaces with
$\C_+-$action}, Comment. Math. Helv. {\bf 69:1} (1994), 5-27.
 \bibitem[FlZa 1]{FlZa 1} H. Flenner, M. Zaidenberg,
{\sl Q-acyclic surfaces and their deformations}, Proc. Conf.
"Classification of Algebraic Varieties", Mai 22--30, 1992, Univ. of l'Aquila,
L'Aquila, Italy, Livorni (ed.) Contempor. Mathem. {\bf 162}, Providence, RI,
1994, 143--208.
 \bibitem[FlZa 2]{FlZa 2} H. Flenner, M. Zaidenberg,
{\em On a class of rational cuspidal plane curves}, Manuscr. Mathem. {\bf 89}
(1996), 439-460.
 \bibitem[For]{For} F. Forstneri\v c, {\em Holomorphic automorphisms of $\C^n:$
A survey}, Complex Analysis and Geometry (Trento, 1993), Ancona e.a.
(eds.) Lect. Notes in Pure and Applied Math. {\bf 173} Marcel Dekker, N.Y.,
1996, 173-200.
 \bibitem[Fu 1]{Fu 1} T. Fujita, {\em On Zariski problem}, Proc. Japan Acad.
Ser. A
Math. Sci. {\bf 55} (1979), 106--110.
 \bibitem[Fu 2]{Fu 2} T. Fujita, {\em On the topology of non complete algebraic
surfaces}, J. Fac. Sci. Univ. Tokyo, Sect.IA, {\bf 29} (1982), 503--566.
 \bibitem[Fur]{Fur} M. Furushima, {\em The complete classification of
compactifications of $\C^3$ which are projective manifolds with second Betti
number equal to one}, Math. Ann. {\bf 297} (1993), 627--662.
 \bibitem[GuMiy 1]{GuMiy 1} R.V. Gurjar, M. Miyanishi,
{\em Affine surfaces with
$\k \le 1$}, Algebraic Geometry and Commutative Algebra, in honor of
M. Nagata, 1987, 99--124.
 \bibitem[GuMiy 2]{GuMiy 2} R.V. Gurjar, M. Miyanishi,
{\em Affine lines on logarithmic $\Q-$homology planes}, Math. Ann. {\bf 294}
(1992), 463--482.
 \bibitem[GuPa]{GuPa} R.V. Gurjar, A.J. Parameswaran, {\em Affine lines on
$\Q-$homology planes}, J. Math. Kyoto Univ. {\bf 35:1} (1995), 63--77.
 \bibitem[GuPraSha]{GuPraSha} R.V. Gurjar, C.R. Pradeep, A.R. Shastri,
{\em On rationality of logarithmic $\Q-$homology planes: I, II, III},
preprints, 1997.
 \bibitem[GuSha]{GuSha} R.V. Gurjar, A.R. Shastri, {\em On rationality of
complex homology 2--cells: I, II}, J. Math. Soc. Japan {\bf 41} (1989), 37--56,
175--212.
 \bibitem[Gut]{Gut} A. Gutwirth, {\em The action of an algebraic torus on the
affine plane}, Trans. Amer. Math. Soc. {\bf 105} (1962), 407 - 414.
 \bibitem[Ha]{Ha} H. A. Hamm, {\em Lefschetz theorems for singular
varieties}, Proceedings of Symposia in Pure Mathematics, Part I
(Arcata Singularities Conference), {\bf 40}, 1983, 547--557.
 \bibitem[Hir]{Hir} F. Hirzebruch, {\em The topology of normal singularities
of an algebraic surface}, S{\'e}minaire Bourbaki
{\bf 15} (1962/1963), No. 250,
9p. (1964).
 \bibitem[HNK]{HNK} F. Hirzebruch, W.D. Neumann, and S.S. Koh, {\em
Differentiable manifolds and quadratic forms}, Lect. Notes Pure Appl. Math.
{\bf 4}, M. Dekker Inc., New York, 1971
 \bibitem[Ho]{Ho} M. Hochster, {\em Non-uniqueness of coefficient rings in a
polynomial ring}, Proc. Amer. Math. Soc. {\bf 34} (1972), 81--82.
 \bibitem[Ii 1]{Ii 1} S. Iitaka, {\em On logarithmic Kodaira dimensions of
algebraic varieties}, Complex Analysis and Algebraic geometry, Iwanami, Tokyo,
1977, 175--189.
 \bibitem[Ii 2]{Ii 2} S. Iitaka, {\em Some applications of logarithmic Kodaira
dimensions}, Algebraic Geometry (Proc. Intern. Symp. Kyoto 1977), Kinokuniya,
Tokyo, 1978, 185--206.
 \bibitem[Ii 3]{Ii 3} S. Iitaka, {\em Algebraic Geometry: An introduction to
birational geometry of algebraic varieties}, Graduate Texts in
Mathematics, {\bf 76}, Springer Verlag, Berlin--Heidelberg--New York,
1982.
 \bibitem[IiFu]{IiFu} S. Iitaka, T. Fujita, {\em Cancellation theorem for
algebraic varieties}, J. Fac. Sci. Univ. Tokyo, Sect.IA, {\bf 24} (1977),
123--127.
 \bibitem [Je]{Je} Z. Jelonek, {\em The extension of regular and rational
embeddings}, Math. Ann. {\bf 113} (1987) 113--120.
 \bibitem [Ju]{Ju} H. W. E. Jung, {\em {\"U}ber ganze birationale
Transformationen der Ebene}, J. reine und angew. Math., {\bf 184} (1942),
161--174.
 \bibitem[Ka 1]{Ka 1} S. Kaliman, {\em Smooth contractible hypersurfaces in
${\C}^{n}$ and exotic algebraic structures on ${\C}^{3}$}, Math. Zeitschrift
{\bf 214} (1993), 499--510.
 \bibitem[Ka 2]{Ka 2} S. Kaliman, {\em Exotic analytic structures and
Eisenman intrinsic measures}, Israel Math. J. {\bf 88} (1994), 411--423.
 \bibitem[Ka 3]{Ka 3} S. Kaliman, {\em Exotic structures on $\C^n$ and
$\C^*$-action on $\C^3$}, Proc. Conf. "Complex Analysis and Geometry",
Lect. Notes in Pure and Appl. Math., Marcel Dekker Inc. {\bf 173} (1996),
299--300.
 \bibitem[Ka 4]{Ka 4} S. Kaliman, {\em Isotopic embeddings of affine algebraic
varieties into $\C^n$}, Contempor. Mathem. {\bf 137} (1992), 291--295.
 \bibitem[Ka 5]{Ka 5} S. Kaliman, {\em Extensions of isomorphisms between
affine algebraic subvarieties of $k^n$ to automorphisms of $k^n$}, Proc.
Amer. Math. Soc. {\bf 113} (1991), 325--334.
 \bibitem[KaML 1]{KaML 1} S. Kaliman, L. Makar-Limanov,
{\em On some family of contractible hypersurfaces in ${\C}^4$},
S{\'e}minaire d'alg{\`e}bre. Journ{\'e}es Singuli{\`e}res et Jacobi{\'e}nnes,
26--28 mai 1993, Pr{\'e}publication de l'Institut Fourier, Grenoble, 1994,
57--75.
\bibitem[KaML 2]{KaML 2} S. Kaliman, L. Makar-Limanov, {\em Affine algebraic
manifolds without dominant morphisms from Euclidean spaces}, Rocky Mount. J.
Math. {\bf 27:2} (1997), 601 --609.
 \bibitem[KaML 3]{KaML 3} S. Kaliman, L. Makar-Limanov, {\em On
Russell--Koras contractible threefolds},
J. of Algebraic Geom. {\bf 6} (1997), 247-268.
 \bibitem[KaML 4]{KaML 4} S. Kaliman, L. Makar-Limanov,
{\em Locally nilpotent derivations of Jacobian type}, Preprint, 1998, 16p.
 \bibitem[KaKoMLRu]{KaKoMLRu} S. Kaliman, M. Koras, L. Makar-Limanov, P.
Russell,
{\em ${\bf C}^*$-actions on ${\bf C}^3$ are linearizable},
Electronic Research
Announcement Journal of the AMS, {\bf 3} (1997), 63--71.
\bibitem[KaZa]{KaZa} S. Kaliman, M. Zaidenberg,
{\em Affine modifications and affine varieties with a very transitive
automorphism group},
Pr{\'e}publication de l'Institut Fourier des Math{\'e}matiques,
406, Grenoble 1998, 46p. E-print: math.AG/9801076
 \bibitem[Kam 1]{Kam 1} T. Kambayashi,
{\em On Fujita's strong cancellation theorem for the affine space},
J. Fac. Sci. Univ. Tokyo {\bf 23} (1980), 535--548.
 \bibitem[Kam 2]{Kam 2} T. Kambayashi, {\em Pro-affine groups,
Ind-affine groups and the Jacobian Problem}, J. Algebra {\bf 185} (1996),
481-501.
 \bibitem[KamRu]{KamRu} T. Kambayashi, P. Russell, {\em On linearizing
algebraic torus actions}, J. Pure Appl. Algebra {\bf 23} (1982), 243--250.
 \bibitem[Kaw 1]{Kaw 1} Y. Kawamata, {\em Addition formula of logarithmic
Kodaira dimension for morphisms of relative dimension one}, Proc. Intern.
Sympos. Algebr. Geom., Kyoto, 1977. Kinokuniya, Tokyo, 1978, 207--217.
 \bibitem[Kaw 2]{Kaw 2} Y. Kawamata,
{\em On the classification of non-complete
algebraic surfaces}, Algebraic Geom., Proc.
Summer Meeting, Copenhagen, 1978,
Lect. Notes Math. {\bf 732}, 1979, 215--232.
 \bibitem[KMM]{KMM} Y. Kawamata, K. Matsuda, K. Matsuki,
{\em Introduction to
the minimal model problem}, Alg. Geom. Sendai, 1985,
Adv. St. in Pure Math. {\bf 10} (1987), 283-360.
 \bibitem[Kl]{Kl} S. Kleiman, {\em Toward a numerical theory of ampleness},
Ann. Math. {\bf 84} (1966), 293-344.
 \bibitem[Ko]{Ko} M. Koras, {\em A characterization of ${\bf A}^2 / \Z_a$},
Compositio Math. {\bf 87} (1993), 241--267.
 \bibitem[KoRu 1]{KoRu 1} M. Koras, P. Russell, {\em On linearizing "good"
${\C}^*$-action on ${\C}^3$}, Canadian Math. Society Conference Proceedings,
{\bf 10} (1989), 93--102.
 \bibitem[KoRu 2]{KoRu 2}  M. Koras, P. Russell, {\em Contractible threefolds
and ${\bf C}^*$-actions on ${\bf C}^3$}, J. of Algebraic Geom. {\bf 6} (1997),
671-695.
 \bibitem[KoRu 3]{KoRu 3}  M. Koras, P. Russell, {\em ${\C}^*$-actions on
${\C}^3$: the smooth locus of the quotient is not of hyperbolic type},
preprint, CICMA Reports, Concordia-Laval-McGill, 1996-{\bf 06}, 93p.
 \bibitem[Kr 1]{Kr 1} H. Kraft, {\em Algebraic automorphisms of affine space},
Topological Methods in Algebraic Transformation Groups, Birkh{\"a}user, Boston
e.a., 1989, 81--105.
 \bibitem[Kr 2]{Kr 2} H. Kraft, {\em $\C^*$--actions on affine space},
Operator
Algebras etc., Progress in Mathem. {\bf 92}, 1990, Birkh{\"a}user,
Boston e.a.,
561--579.
 \bibitem[Kr 3]{Kr 3} H. Kraft, {\em Challenging problems on affine $n$-space},
S{\'e}minaire Bourbaki {\bf 802} (1994/1995), Ast{\'e}risque {\bf 237} (1996),
295--318.
 \bibitem[KrPeRun]{KrPeRun} H. Kraft, T. Petrie, J.D. Rundall, {\em Quotient
varieties}, Advances in Math. {\bf 74} (1989), 145-162.
 \bibitem[KrPo]{KrPo} H. Kraft, V. L. Popov, {\em Semisimple group actions on
the three dimensional affine space are linear}, Comment. Math. Helv.
{\bf 60} (1985), 466--479.
 \bibitem[Lef]{Lef} S. Lefschetz,
{\em L'analysis situs et la g{\'e}om{\'e}trie
alg{\'e}brique}, Paris, 1924.
 \bibitem[Lib]{Lib} A. Libgober, {\em A geometric procedure for killing the
middle dimensional homology groups of algebraic hypersurfaces.} Proc. Amer.
Math. Soc. {\bf 63} (1977), 198--202.
 \bibitem[LiZa]{LiZa} V. Lin, M. Zaidenberg, {\sl An irreducible simply
connected curve in ${\C}^{2}$ is equivalent to a quasihomogeneous curve},
Soviet Math. Dokl., {\bf 28} (1983), 200-204.
 \bibitem[ML 1]{ML 1} L. Makar-Limanov,
{\em On groups of automorphisms of a class of surfaces},
Israel J. Math. {\bf 69} (1990), 250-256.
 \bibitem[ML 2]{ML 2} L. Makar-Limanov,
{\em On the hypersurface $x + x^2y + z^2 + t^3 = 0$
in ${\C}^{4}$ or a ${\C}^3$-like threefold which is not ${\C}^3$}, Israel J.
Math. {\bf 96} (1996), 419--429.
 \bibitem[ML 3]{ML 3} L. Makar-Limanov,
{\em Again $x + x^2y + z^2 + t^3 = 0$}, Preprint, 1998, 3p.
 \bibitem[ML 4]{ML 4} L. Makar-Limanov,
{\em On the group of automorphisms of a surface $x^ny = P(z)$},
Preprint, 1997, 11p.
 \bibitem[Mil 1]{Mil 1} J. Milnor, {\em Lectures on the h-cobordism Theorem},
Princeton Univ. Press, Princeton, NJ, 1965.
\bibitem[Mil 2]{Mil 2} J. Milnor, {\em Morse Theory}, Princeton Univ. Press,
Princeton, NJ, 1963.
\bibitem[Mil 3]{Mil 3} J. Milnor, {\em On the 3-dimensional Brieskorn manifolds
$M(p,\,q,\,r)$}, in: Knots, groups, and 3-manifolds, L. P. Neuwirth, ed.
Annals of Math. Stud., Princeton Univ. Press,
Princeton, NJ, 1975, 175--225.
 \bibitem[MilSta]{MilSta} J. Milnor, J. Stasheff, {\em Characteristic classes},
Annals of Mathem. Studies {\bf 76}, Princeton Univ. Press and Univ. of Tokyo
Press, Princeton, NJ, 1974.
 \bibitem[Miy 1]{Miy 1} M. Miyanishi, {\em Algebraic characterization of the
affine 3--space}, Proc. Algebraic Geom. Seminar, Singapore, World Scientific,
1987, 53--67.
 \bibitem[Miy 2]{Miy 2} M. Miyanishi, {\em Recent topics on open algebraic
surfaces}, Amer. Math. Soc. Transl. {\bf 172} (1996), 61-75.
 \bibitem[MiySu]{MiySug} M. Miyanishi, T. Sugie, {\em Affine surfaces
containing cylinderlike open sets}, J. Math. Kyoto Univ., {\bf 20} (1980),
11--42.
 \bibitem[Mo]{Mo} S. Mori, {\em Classification of higher-dimensional
varieties},
Proc. Symp. in Pure Math. {\bf 46} (1987), 269-331
 \bibitem[MS]{MS} S. M{\"u}ller--Stach, {\em Projective compactifications of
complex affine varieties}, London Math. Soc. Lect. Notes Ser. {\bf 179} (1991),
277--283.
 \bibitem[Mu]{Mu} D. Mumford, {\em The topology of normal singularities of an
algebraic surface and a criterion for simplicity}, Public. Math. IHES {\bf 9}
(1961), 229--246.
 \bibitem[Na 1]{Na 1} M. Nagata, {\em On an automorphism group of $k[x,y]$},
Kinokuniya, Tokyo, 1972.
 \bibitem[Na 2]{Na 2} M. Nagata, {\em Commutative algebra and algebraic
geometry},
Proc. Intern. Mathem. Conf. L.H.Y. Chen, T.B.Ng, M.J.Wicks (eds.),
North-Holland Publ. Co., 1982, 125-154.
 \bibitem[Ne]{Ne} A. N{\'e}methi, {\em Global Sebastiani--Thom theorem for
polynomial maps}, J. Math. Soc. Japan, {\bf 43} (1991), 213--218.
 \bibitem[OPOV]{OPOV} {\em Open problems on open varieties (Montreal 1994
problems)}, P. Russell (ed.),
Pr{\'e}publication de l'Institut Fourier des Math{\'e}matiques, {\bf 311},
Grenoble 1995, 23p. E-print alg-geom/9506006.
  \bibitem[Or 1]{Or 1} S. Orevkov, {\em On three-sheeted polynomial mappings of
$\C^2$}, Math. USSR Izvestiya, {\bf 29} (1987), 587-598.
  \bibitem[Or 2]{Or 2} S. Orevkov, {\em Acyclic algebraic surfaces bounded by
Seifert spheres}, Osaka J. Math. {\bf 34:2} (1997), 457--480.
 \bibitem[Pe]{Pe} T. Petrie, {\em Topology, representations and equivariant
algebraic geometry}, Contemporary Math. {\bf 158}, 1994, 203--215.
 \bibitem[Pi]{Pi} S. Pinchuk,
{\em A counterexample to the strong real Jacobian conjecture},
Math. Z. {217:1} (1994), 1-4.
 \bibitem[Po]{Po} V.L. Popov, {\em Algebraic actions of connected reductive
algebraic groups on ${\bf A}^3$ are linearizable}, preprint, 1996, 3p.
 \bibitem[Pro]{Pro} Y.G. Prokhorov,
{\em Compactifications of $\C^4$ of index 3},
Algebraic Geometry and its Applications, Proc. 8th Algebraic Geometry Conf.,
Yaroslavl' 1992, Tikhomirov, Tyurin (eds.)
Vieweg, Braunschweig/Wiesbaden, 1994, 159-169.
 \bibitem[Ram]{Ram} C.P. Ramanujam, {\sl A topological characterization of the
affine plane as an algebraic variety}, Ann. Math., {\bf 94} (1971), 69-88.
 \bibitem[Re]{Re} R. Rentschler, {\em Op{\'e}rations du groupe additif sur le
plane affine}, C.R. Acad. Sci. Paris, {\bf 267} (1968), 384--387.
 \bibitem[Ru 1]{Ru 1} P. Russell, {\em On a class of ${\C}^3$-like threefolds},
Preliminary Report, 1992.
 \bibitem[Ru 2]{Ru 2} P. Russell, {\em On affine-ruled rational surfaces},
Math. Ann., {\bf 255} (1981), 287--302.
 \bibitem[Ru 3]{Ru 3} P. Russell, {\em Gradings of polynomial rings},
Algebraic Geometry and its Applications (C. L. Bajaj ed.), Springer, 1994.
 \bibitem[Ru 4]{Ru 4} P. Russell, {\em Sufficiently homogeneous closed
embeddings of ${\bf A}^{n-1}$ into ${\bf A}^{n}$ are linear}, Preprint CICMA,
1997, 13p.
 \bibitem[Sak]{Sak} F. Sakai, {\em Kodaira dimension of complement of divisor},
Complex Analysis and Algebraic geometry, Iwanami, Tokyo, 1977, 239--257.
 \bibitem[Sat]{Sat} A. Sathaye, {\em On linear planes}, Proc. Amer. Math. Soc.
{\bf 56} (1976), 1--7.
 \bibitem[Sch]{Sch} G. W. Schwarz, {\em Exotic algebraic group actions},
C. R. Acad. Sci. Paris {\bf 309} (1989), 89--94.
 \bibitem[Sn]{Sn} D. Snow, {\em Unipotent actions on affine space},
Topological methods in algebraic transformation groups, Proc. Conf., New
Brunswick/NJ (USA) 1988, Prog. Math. {\bf 80} (1989), 165-176.
 \bibitem[Sr]{Sr} V. Srinivas, {\em On the embedding dimension of an affine
variety}, Math. Ann.  {\bf 289} (1991), 125-132.
 \bibitem[Sug]{Sug} T. Sugie, {\em On Petrie's problem concerning homology
planes},
J. Math. Kyoto Univ. {\bf 30} (1990), 317-324.
 \bibitem[Suz 1]{Suz 1} M. Suzuki, {\em Propi{\'e}t{\'e}s topologiques des
polyn{\^o}mes de deux variables complexes, et automorphismes alg{\'e}brigue de
l'espace
${\C}^{2}$}, J. Math. Soc. Japan, {\bf 26} (1974), 241-257.
 \bibitem[Suz 2]{Suz 2} M. Suzuki, {\em Sur les op{\'e}ration holomorphes
du groupe additif complexe sur l'espace de deux variables complexes},
Ann. Sci. {\'E}cole Norm.  Sup. {\bf 10} (1977), 517--546.
 \bibitem[tD 1]{tD 1} T. tom Dieck, {\em Hyperbolic modifications and acyclic
affine foliations}, preprint, Mathematica Gottingensis, G{\"o}ttingen,
H. {\bf 27} (1992), 1--19.
 \bibitem[tD 2]{tD 2} T. tom Dieck, {\em Ramified coverings of acyclic
varieties}, preprint, Mathematica Gottingensis, G{\"o}ttingen, H. {\bf 26}
(1992), 1--20.
 \bibitem[tD 3]{tD 3} T. tom Dieck, {\em Homology planes without cancellation
property}, Arch. Math. {\bf 59} (1992), 105--114.
 \bibitem[tD 4]{tD 4} T. tom Dieck, {\em Symmetric homology planes}, Math. Ann.
{\bf 286} (1990), 143-152.
 \bibitem[VdV]{VdV} A. Van de Ven, {\em Analytic compactifications of complex
homology cells}, Math. Ann. {\bf 147} (1962), 189--204.
 \bibitem[vdK]{vdK} W. van der Kulk,
{\em On polynomial rings in two variables}, Nieuw Arch. Wisk. (3) {\bf 1}
(1953), 33--41.
 \bibitem[Vie]{Vie} E. Viehweg, {\em Canonical divisors and the additivity
of the Kodaira dimension for morphisms of relative dimension one},
Compos. Math. {\bf 35} (1977). 197-223, {\em Correction}, ibid. 336.
  \bibitem[Wa]{Wa} S. S.-S. Wang, {\em Extension of derivations},
J. Alg. {\bf 69} (1981), 240-246.
  \bibitem[Wil]{Wil} P.M.H. Willson, {\em Towards birational classification of
algebraic varieties}, Bull. London Math. Soc. {\bf 19} (1987), 1-48.
 \bibitem[Win]{Win} J. Winkelmann, {\em On free holomorphic $\C_+$--actions on
$\C^n$ and homogeneous Stein manifolds}, Math. Ann. {\bf 286} (1990), 593--612.
 \bibitem[Wr 1]{Wr 1} D. Wright, {\em Abelian subgroups of ${\rm
Aut}_k(k[X,Y])$
and applications to actions on the affine plane}, Ill. J. Math. {\bf 23}
(1979), 579--634.
\bibitem[Wr 2]{Wr 2} D. Wright, {\em Cancellation of variables of the form
$bT^n - a$}, J. Algebra {\bf 52} (1978), 94--100.
 \bibitem[Za 1]{Za 1} M. Zaidenberg, {\sl Isotrivial families of curves on
affine surfaces and characterization of the affine plane}, Math. USSR Izvestiya
{\bf 30} (1988), 503-531. Addendum: ibid, {\bf 38} (1992), 435--437.
 \bibitem[Za 2]{Za 2}  M. Zaidenberg, {\sl Ramanujam surfaces and exotic
algebraic
structures on ${\C}^n$}, Soviet Math. Doklady {\bf 42} (1991), 636--640.
 \bibitem[Za 3]{Za 3}  M. Zaidenberg, {\sl An analytic cancellation theorem and
exotic algebraic structures on ${\C}^n$, $n \ge 3$},  Ast{\'e}risque {\bf 217}
(1993), 251--282.
\bibitem[Za 4]{Za 4} M. Zaidenberg, {\sl On Ramanujam surfaces,
$\C^{**}$-families
and exotic algebraic structures on $\C^n$, $n \ge 3$},
Trans. Moscow Math. Soc. {\bf 55} (1994), 1--56.
\bibitem[Za 5]{Za 5} M. Zaidenberg, {\sl On exotic algebraic structures on
affine spaces}, Geometric Complex Analysis, J. Noguchi e.a. eds. World
Scientific Publ. Co., Singapore 1996, 691--714; E-print alg-geom/9506005.


\end{thebibliography}
\end{document}